\documentclass[a4paper]{amsart}
\usepackage[utf8]{inputenc}

\usepackage{bbm}
\usepackage{amsmath}
\usepackage{amsfonts}
\usepackage{amssymb}
\usepackage{amsthm}
\usepackage[foot]{amsaddr}
\usepackage{enumerate}
\usepackage{xcolor}
\usepackage{mathrsfs}
\usepackage{stmaryrd}
\usepackage{graphicx}
\usepackage{caption}
\usepackage[a4paper,margin=1in,headsep=25pt,footskip=35pt]{geometry}
\usepackage{fancyhdr}
\usepackage{comment}
\usepackage{hyperref}
\usepackage{cleveref}
\usepackage{tikz}
\usetikzlibrary{shapes.geometric}
\usetikzlibrary{decorations.shapes}
\usetikzlibrary{decorations.pathmorphing, decorations.pathreplacing}

\SetSymbolFont{stmry}{bold}{U}{stmry}{m}{n}

\newtheorem{thm}{Theorem}

\newtheorem{lem}[thm]{Lemma}
\newtheorem{prop}[thm]{Proposition}

\newtheorem{defn}[thm]{Definition}
\newtheorem*{defn*}{Definition}
\newtheorem{assump}[thm]{Assumption}

\newtheorem{rem}[thm]{Remark}

\newcommand{\deq}{\mathrel{\mathop:}=}
\newcommand{\e}[1]{\mathrm{e}^{#1}}
\newcommand{\R} {\mathbb{R}}
\newcommand{\C} {\mathbb{C}}

\newcommand{\N} {\mathbb{N}}
\newcommand{\Z} {\mathbb{Z}}

\newcommand{\E} {\mathbb{E}}

\newcommand{\adj}{^{*}} 
\newcommand{\tp}{^{\intercal}}

\DeclareMathOperator{\diag}{diag}

\DeclareMathOperator{\Tr}{Tr}

\DeclareMathOperator{\supp}{supp}
\DeclareMathOperator{\spann}{span}
\DeclareMathOperator{\spec}{spec}

\DeclareMathOperator{\re}{\mathrm{Re}}
\DeclareMathOperator{\im}{\mathrm{Im}}

\newcommand{\caA}{{\mathcal A}}
\newcommand{\caB}{{\mathcal B}}

\newcommand{\caF}{{\mathcal F}}

\newcommand{\caH}{{\mathcal H}}
\newcommand{\caI}{{\mathcal I}}
\newcommand{\caJ}{{\mathcal J}}

\newcommand{\caL}{{\mathcal L}}
\newcommand{\caM}{{\mathcal M}}

\newcommand{\caP}{{\mathcal P}}

\newcommand{\caS}{{\mathcal S}}

\newcommand{\caU}{{\mathcal U}}
\newcommand{\caV}{{\mathcal V}}

\newcommand{\caY}{{\mathcal Y}}

\newcommand{\bbC}{{\mathbb C}}

\newcommand{\bbF}{{\mathbb F}}

\newcommand{\bbH}{{\mathbb H}}

\newcommand{\frc}{\mathfrak{c}}

\newcommand{\frg}{{\mathfrak g}}

\newcommand{\frC}{{\mathfrak C}}

\newcommand{\bsa}{{\boldsymbol a}}

\newcommand{\bsc}{{\boldsymbol c}}

\newcommand{\bsh}{{\boldsymbol h}}

\newcommand{\bsn}{{\boldsymbol n}}

\newcommand{\bss}{{\boldsymbol s}}

\newcommand{\bsu}{{\boldsymbol u}}

\newcommand{\bsw}{{\boldsymbol w}}
\newcommand{\bsx}{{\boldsymbol x}}
\newcommand{\bsy}{{\boldsymbol y}}
\newcommand{\bsz}{{\boldsymbol z}}

\newcommand{\bsL}{{\boldsymbol L}}

\newcommand{\bsU}{{\boldsymbol U}}

\newcommand{\rmd}{\mathrm{d}}

\newcommand{\wt}{\widetilde}
\newcommand{\ol}{\overline}

\newcommand{\mr}{\mathring}
\newcommand{\wh}{\widehat}
\newcommand{\beq}{ \begin{equation} }
	\newcommand{\eeq}{ \end{equation} }
\newcommand{\beqs}{\begin{equation*}}
	\newcommand{\eeqs}{\end{equation*}}

\newcommand{\lone}{\mathbbm{1}} 

\newcommand{\dd}{\mathrm{d}}
\newcommand{\ii}{\mathrm{i}}

\newcommand{\llbra}{\llbracket}
\newcommand{\rrbra}{\rrbracket}

\newcommand\norm[1]{\Vert#1\Vert}
\newcommand\Norm[1]{\left\Vert#1\right\Vert}

\newcommand\Absv[1]{\left\vert#1\right\vert}
\newcommand\absv[1]{\vert#1\vert}

\newcommand\brkt[1]{\langle#1\rangle}
\newcommand\Brkt[1]{\left\langle#1\right\rangle}
\newcommand\bbrktt[1]{\llbra #1\rrbra}

\numberwithin{equation}{section} 
\numberwithin{thm}{section}

\title{Non--Hermitian spectral universality at critical points}

\author[]{Giorgio Cipolloni$^{\ddagger}$}
\email{gcipolloni@arizona.edu}

\author[]{L\'{a}szl\'{o} Erd\H{o}s$^{\dagger}$}
\email{lerdos@ist.ac.at}

\author[]{Hong Chang Ji$^*$}
\email{hji56@wisc.edu}

\address{$^{\dagger}$Institute of Science and Technology Austria}
\address{$^{\ddagger}$ Mathematics department, University of Arizona}
\address{$^*$ Mathematics department, University of Wisconsin}
\thanks{$^\dagger$ Supported by ERC Advanced Grant ``RMTBeyond" No. 101020331}

\begin{document}

\begin{abstract}
For general large non--Hermitian random matrices $X$ and deterministic normal
deformations $A$, we prove that the local eigenvalue statistics of $A+X$ close to the critical edge points of its spectrum are universal. This concludes the proof of the third and last remaining typical universality class for non--Hermitian random matrices (for normal deformations), after bulk and sharp edge universalities have been established in recent years.
\end{abstract}
\maketitle
	\section{Introduction}
		
The Wigner--Dyson--Mehta conjecture states that the local statistics of the eigenvalues of large Hermitian matrices are universal, i.e. they do not depend on the details of the system but only on its symmetry class. Moreover, different statistics emerge in the \emph{bulk} of the spectrum, at the \emph{regular spectral edges} with a square root vanishing behavior of the eigenvalue density, and at the \emph{cusps} where the density vanishes as a cubic root. The fact that for a very large class of mean-field random matrices  these are the only possible emerging statistics follows by the classification theorem \cite{Ajanki-Erdos-Kruger2017a}. By now the universality of these statistics is  well understood for very general Hermitian random matrix ensembles in the bulk \cite{Tao-Vu2011, Erdos-Schlein-Yau2011}, regular edge \cite{Soshnikov1999, Tao-Vu2010CMP,  Bourgade-Erdos-Yau2014}, and cusp \cite{Erdos-Kruger-Schroder2020, Cipolloni-Erdos-Kruger-Schroder2019} regimes. Since the literature about this topic is very vast, here we refer only to the first results proving universality in each one of these regimes, typically for the simplest ensemble. Many subsequent works strengthened  these universality results and also extended them  to more general ensembles. We refer the interested reader to \cite[Section 1.2]{Campbell-Cipolloni-Erdos-Ji2024} for further references and historical notes.

\subsection{Main result}
	
Motivated by these Hermitian results, we aim to study a similar universality question for 
non--Hermitian matrices. They exhibit a more intricate scenario since their spectrum is genuinely in the
two--dimensional plane,
moreover their analysis  is technically more challenging. 
In this paper we consider $N\times N$ \emph{deformed i.i.d. matrices}, i.e. matrices of the form $A+X$, where $X$ has independent, identically distributed centered complex entries (\emph{i.i.d. matrices}) with $\E X_{ij}^2=0$ and $\E|X_{ij}|^2= N^{-1}$, and $A$ is a normal deterministic  matrix, called \emph{deformation}. Under mild assumptions on $A$, 
a non-Hermitian  classification theorem  was established in \cite{Erdos-Ji2023circ} 
asserting that for large $N$ the limiting eigenvalue density of $A+X$ exhibits only two 
possible  behaviors close to the edges of its support: i) \emph{sharp edge} points, where the density has a jump of order one, ii) \emph{critical edge} points, where the density vanishes quadratically (possibly except for one direction
where vanishing of all polynomial order may occur, see  \cite{Alt-Kruger2024Brown}  for a complete classification
of this special situation).\footnote{See also \cite{Alt-Kruger2024} for an extension of \cite{Erdos-Ji2023circ} for the case when $X$ has a variance profile but $A$ is diagonal.}
Since local eigenvalue statistics are expected to be determined by the local behaviour of the limiting density, 
the result of \cite{Erdos-Ji2023circ} indicates that
also for non--Hermitian matrices only three possible different statistics can emerge: \emph{bulk}, \emph{sharp edge}, and \emph{critical edge}.  Explicit form of the correlation functions at the bulk and at the sharp edge
have been computed in \cite{Borodin-Sinclair2009, Forrester-Nagao2007, Ginibre1965, Mehta2004} based upon the simplest Gaussian i.i.d. matrix ensemble (also called \emph{Ginibre} matrices).  The statistics close to critical edge points,
 computed for certain concrete examples by Liu and Zhang \cite{Liu-Zhang2023} only very recently, are much richer than the ones appearing in the bulk and close to sharp edges.
In fact, unlike in these regimes, the eigenvalue fluctuation close to critical edges are not described by a single limiting distribution but by a two--parameter family of statistics parametrized by the values of $\alpha, \beta$ defined later in \eqref{eq:defangle} and \eqref{eq:defbneta}, respectively. The parameter $\alpha$ is determined by the shape
of the quadratic behavior of the density (it is the ratio of the two eigenvalues of the Hessian), while 
the parameter $\beta$ describes a small deviation from the exact criticality.
This $\beta$ parameter creates an analogy with the statistics close to cusps in the Hermitian case, 
which are described by a one--parameter family of statistics given by the \emph{Pearcey kernel} \cite{Tracy-Widom2006}. 
The parameter $\alpha$ has no Hermitian analogue since in one--dimension there is only 
one cubic singularity (modulo trivial rescaling), while a quadratic density function in two dimensions
has one intrinsic parameter.

The main result of this paper (Theorem~\ref{theo:mainthm} below)  is the proof of universality
of local eigenvalue statistics close to the critical edges of the limiting eigenvalue density of $A+X$. This result, together with the recent proofs of bulk and edge universality in \cite{Dubova-Yang2024, Maltsev-Osman2023, Osman2024} and \cite{Campbell-Cipolloni-Erdos-Ji2024} (see also for previous \cite{Cipolloni-Erdos-Schroder2021, Tao-Vu2015} and partial \cite{Liu-Zhang2024} results), respectively, proves the last remaining case of universality of local spectral correlation functions for non--Hermitian matrices.  We also mention recent progress in the proof of universality of the spectral radius of i.i.d. matrices \cite{Cipolloni-Erdos-Xu2023}, the universality of the eigenvector entries \cite{Dubova-Yang-Yau-Yin2024,Osman2024vector}, and the universality of the eigenvector overlap \cite{Osman2024overlap}.

\begin{figure}
	\centering
	\includegraphics[width=0.32\textwidth]{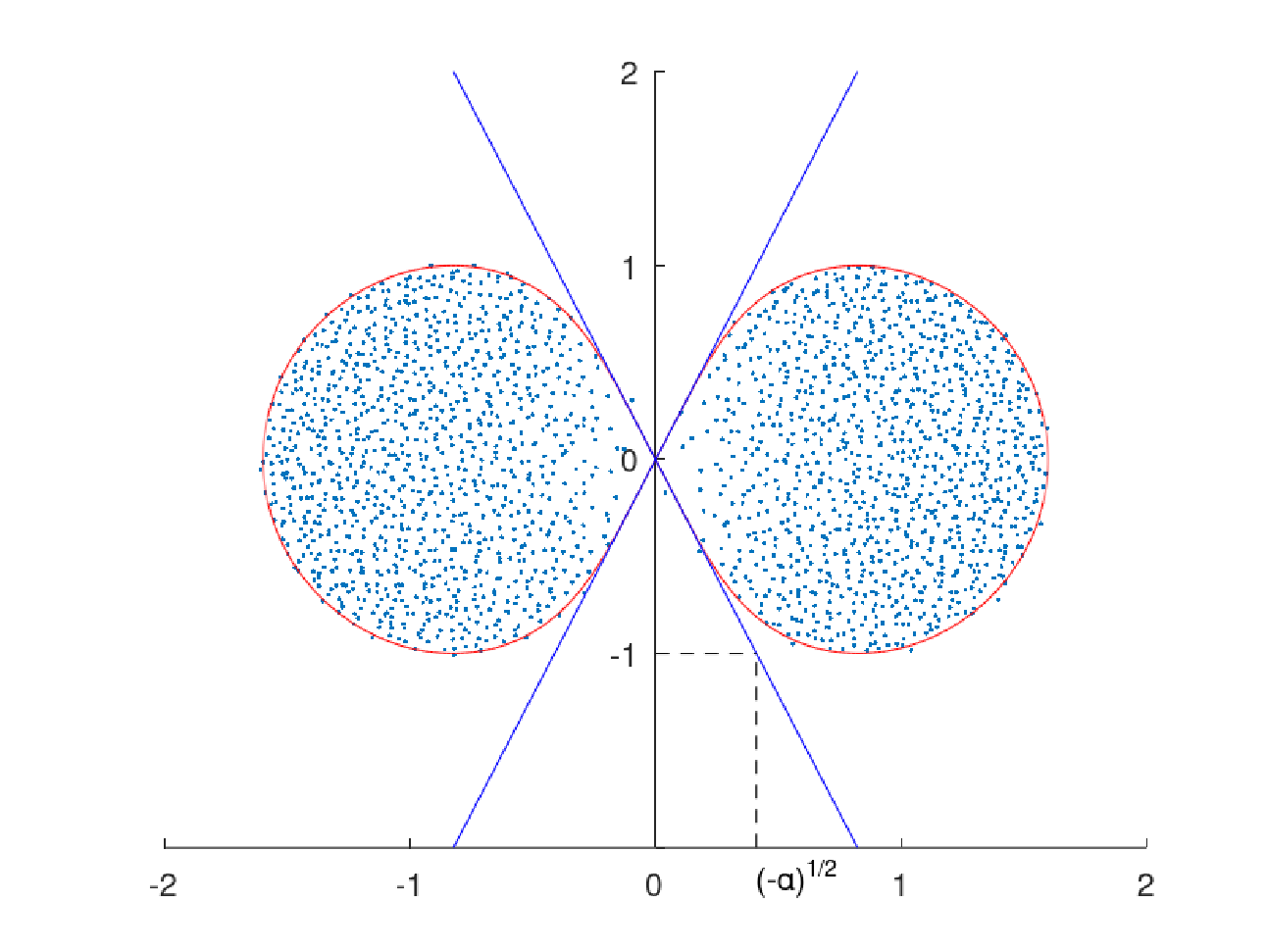}
	\includegraphics[width=0.32\textwidth]{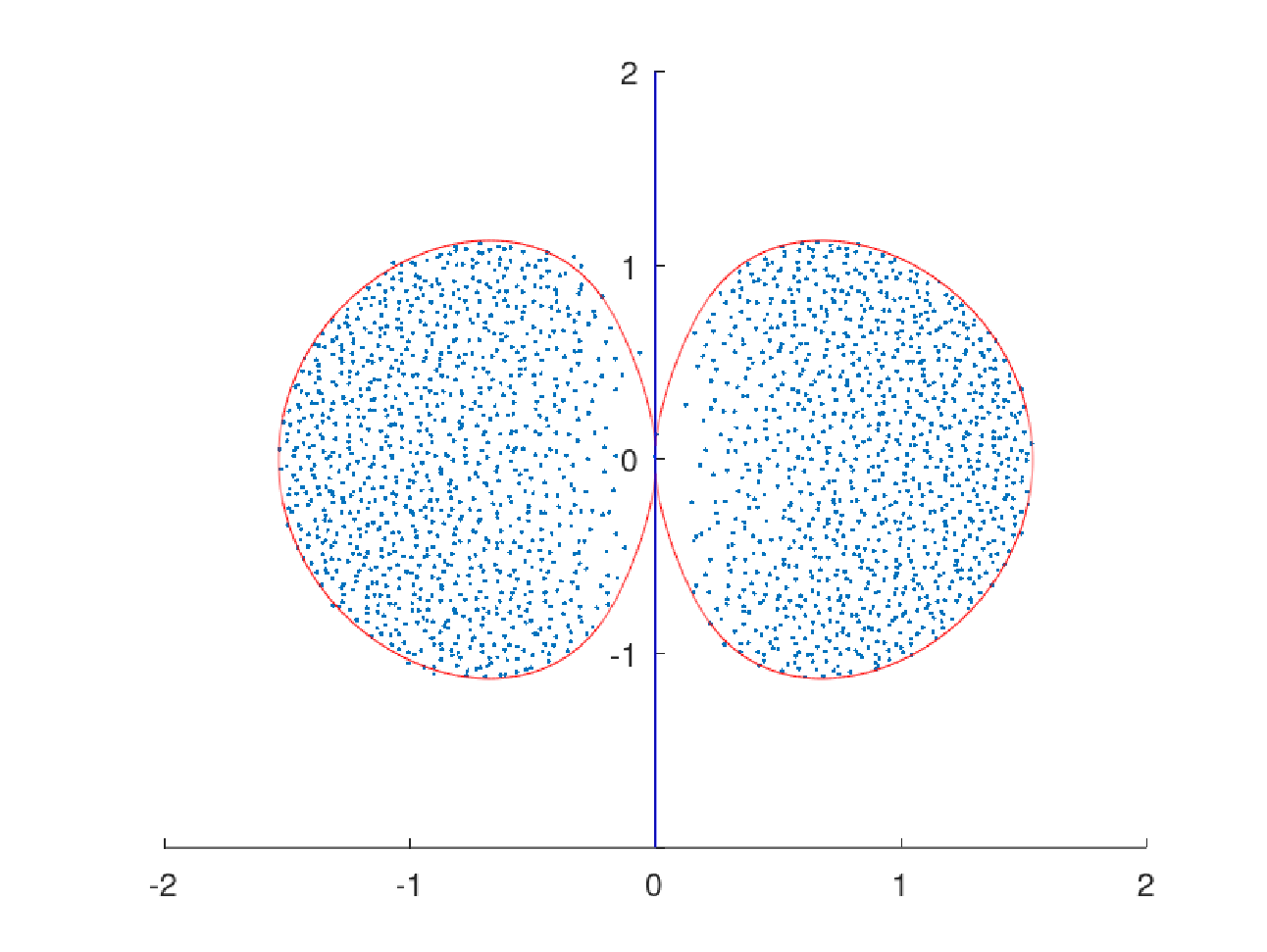}
	\includegraphics[width=0.32\textwidth]{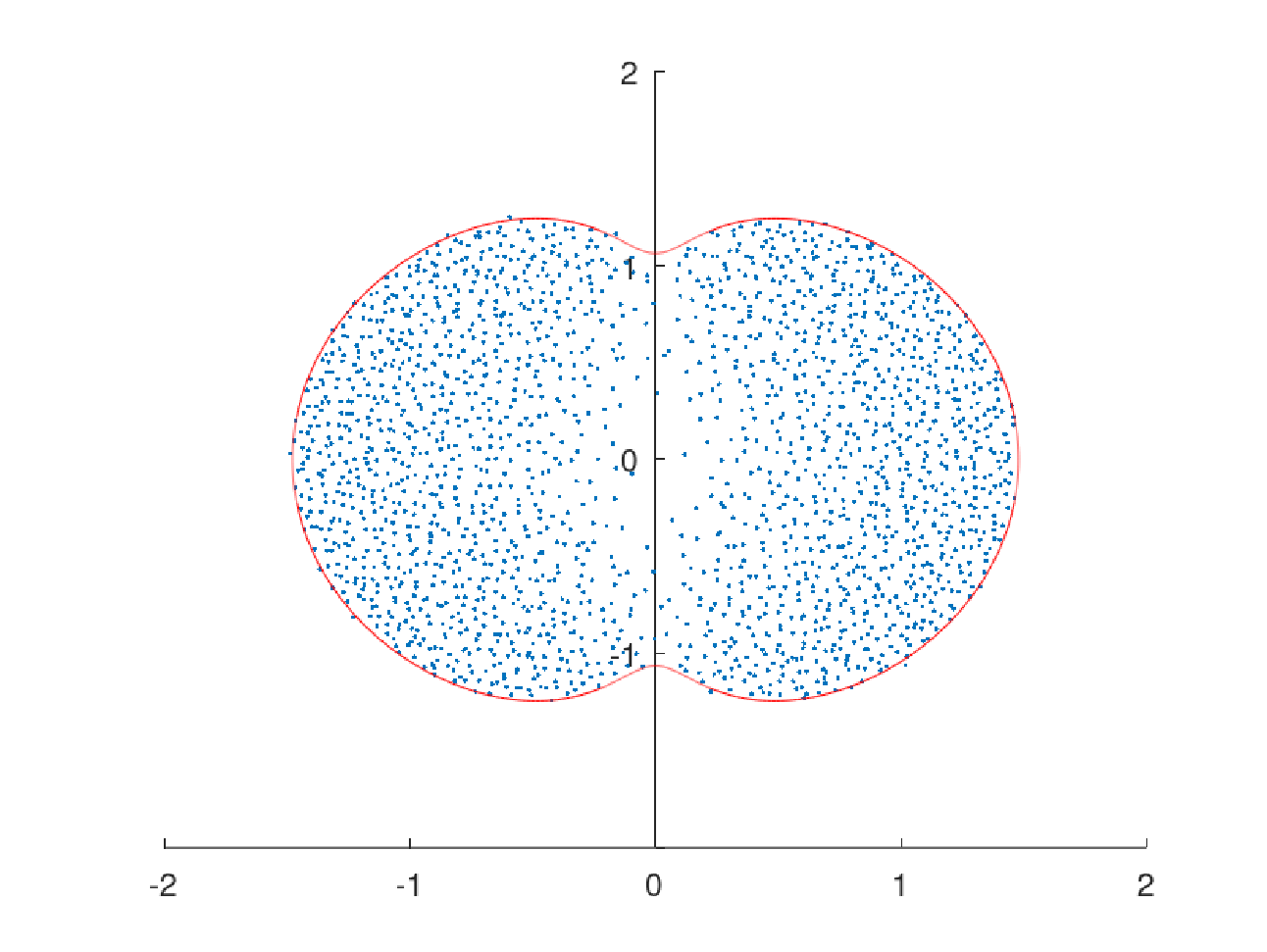}
	\caption{Critical edges at the origin for $\alpha<0$, $\alpha=0$, and $\alpha>0$ (from left to right), all with $\beta=0$; blue dots are $\spec(A+X)$, red curve is the edge, and blue solid lines are the tangent lines at the origin. }\label{fig:alpha}
\end{figure}

While we consider universality only under normal deformations $A$ in this paper, we see that new statistics 
at critical points can possibly emerge for non--normal deformations. In fact, as a byproduct of our analysis, in Lemma~\ref{lem:alpha} below, we show that the parameter $\alpha$, determining the local spectral statistics at the
 critical edges, has a wider range $\alpha\in (-1,1]$ for general deformations $A$ compared to the range $\alpha \in [-1/3,1]$ achievable by normal deformations. Proving universality for non--normal deformations goes beyond the scope and the techniques presented in this work. Furthermore, throughout the paper we only consider i.i.d. matrices $X$ with 
genuinely \emph{complex} entries. One could ask similar questions for real $X$, however we refrain from doing so 
 here since the local statistics close to critical edges in the real case are not 
 yet explicitly known (\cite{Liu-Zhang2023} considers only the complex case).
 Focusing on the complex case also  keeps the presentation shorter, but our approach 
 with relatively easy modifications would also work for the real case once the analogue of \cite{Liu-Zhang2023} 
 is established.

\subsection{Method and novelties} As it is customary in the study of spectral properties of non--Hermitian matrices we rely on Girko's Hermitization formula \cite{Girko1984} to express the spectral statistics of $A +X$ in terms of those of the Hermitization of $A+X-z$, for $z\in\C$,  which is defined by
\begin{equation}
\label{eq:hermintro}
H^z:=\left(\begin{matrix} 0 & A+X-z \\
(A+X-z)^* & 0
\end{matrix}\right).
\end{equation}
In fact, denoting the normalized empirical eigenvalue measure of $A+X$ by $\rho_{A+X}$,  Girko's formula asserts
\begin{equation}
\label{eq:girko}
\int_\C F(z)\dd \rho_{A+X}(z)=-\frac{1}{4\pi N}\int_\C \Delta F(z) \log\big|\mathrm{det}(H^z)\big|\,\dd^2 z
\end{equation}
for any smooth and compactly supported test function $F$ on $\C$.  The main advantage of \eqref{eq:girko}
is that it converts a  non--Hermitian problem into a one--parameter family of Hermitian problems
which are technically easier to handle.  In particular, using the identity
\begin{equation}\label{int}
  -\log\big|\mathrm{det}(H^z)\big| = \int_0^\infty\left[\Im \Tr G^z(\ii\eta) -\frac{2N\eta}{1+\eta^{2}}\right]\dd\eta 
\end{equation}
where $G^z(\ii\eta)= (H^z-\ii\eta)^{-1}$ is the Green function of $H^z$,  \emph{local laws} are available
that assert that $G^z(\ii\eta)$ becomes deterministic 
in the large $N$ limit if $\eta$ is not too small.

Given \eqref{eq:girko}--\eqref{int}, the proof of universality for the local eigenvalue statistics of $A+X$ can be divided into three steps. First, we prove that the statistics of $A+X$ and $A+X^{\mathrm{Gin}(\C)}$ near critical edges
are asymptotically the same. Using the local law from \cite[Theorem 3.3]{Campbell-Cipolloni-Erdos-Ji2024} as an input, this 
will be relatively easily achieved by a two moment matching Green function comparison argument as in \cite{Cipolloni-Erdos-Schroder2021}. Second, we show that the statistics of $A+X^{\mathrm{Gin}(\C)}$ are close to those of $B+X^{\mathrm{Gin}(\C)}$, where $B$ is a carefully constructed matrix with spectrum supported only on finitely many points. Finally, the statistics of $B+X^{\mathrm{Gin}(\C)}$ can be computed relying on \cite{Liu-Zhang2023}. The main novel contribution of this work is thus to show the second step above that involves Ginibre matrices only. For this reason in the rest of the introduction we write $X=X^{\mathrm{Gin}(\C)}$ for the Ginibre ensemble for notational simplicity.

One fundamental difficulty in the study of the spectrum of non--Hermitian matrices is to establish 
a good estimate for the lower tail of the smallest singular value of  $A+X-z$,
that handles the small $\eta$ regime in the integral \eqref{int} which is not covered by the local law.
 The presence of deformation $A$ is a major complication in proving such an estimate directly, 
 without relying on a comparison with a
 simpler reference ensemble. Typically such a direct method involves explicit computations using the Gaussianity of $X$ (e.g. supersymmetry), 
which becomes rather involved with a non--trivial
 deformation; see e.g. \cite{Shcherbina-Shcherbina2022} where the distribution of the smallest singular values of 
$A+X-z$ was computed for $z$ close to the sharp edge using a supersymmetric approach.
 Instead, in \cite{Campbell-Cipolloni-Erdos-Ji2024} this was achieved by comparing the law of small singular values of $A+X-z$, 
 when $z$  is close to the sharp edge of the spectrum of $A+X$, with those of 
$X-z$ for $z$ close to the sharp edge $1$ of $X$. 
For the Ginibre ensemble, the desired smallest singular value bound for $X-1$ was already known in \cite{Cipolloni-Erdos-Schroder2020pmp, Cipolloni-Erdos-Schroder2022s}.
  
 Likewise, in the present work, we also use the comparison method to bypass the difficulty in studying the smallest singular value of $A+X-z$ directly. In fact, thanks to a very useful observation, we can add 
 another layer of simplification: as far as the small singular values are concerned, we may compare
 the critical edge with a sharp edge of a reference ensemble. The point is that while the non--Hermitian 
 spectrum of $A+X$ is sensitive to critical edges (criticality), the spectral density $\mu^z$ of the Hermitization, hence the singular values
 of $A+X-z$, do not feel criticality; they are governed by the very same cubic singularity in both cases. 
 On the Hermitized level, criticality is manifested only in the way how $\mu^z$ changes as $z$ varies. 

 To be precise, for our proof we show that 
  the distribution of the singular values of $A+X-z$, for fixed $z$ close to the edge,
   does not depend on whether the edge is sharp or critical. 
 In particular, in the proof of Lemma~\ref{lem:smalletaerr} below we show that the singular values of 
$A+X-z$ 
have the same joint law as
those of $|A-z|+X$ and that, when $z$ is close to the critical edge of $A+X$, 
the latter ensemble $\absv{A-z}+X$ has a sharp edge near the origin. Then, to prove the desired bound
on the tail of the lowest singular value,
we simply use the result from \cite[Eq. (4.20)]{Campbell-Cipolloni-Erdos-Ji2024} on the singular value estimate around the sharp edge, which itself relies on the explicit computations for $X-1$ 
performed in \cite{Cipolloni-Erdos-Schroder2020pmp}.

Once the possibility to have a small singular value is excluded, the next step is to compare the statistics of 
$A+X$ close to a critical edge with those of $B+X$, using Girko's formula \eqref{eq:girko}--\eqref{int}. 
To do so, we evolve $\mathcal{A}_0:=A-z$ continuously along a carefully chosen flow $(\mathcal{A}_t)_{t\in [0,1]}$ such that $\mathcal{A}_t+X$ has a critical edge at the same point (origin) with the same parameter $\alpha$ for any $t\in [0,1]$, and so that 
$\mathcal{A}_1=B-z$ is 
supported on finitely many points. Along this flow, we show that the time derivatives of the
Hermitized resolvents are sufficiently small,  hence the
 eigenvalue and singular value statistics 
  are approximately preserved. 
   The construction of the flow $\caA_{t}$, which deforms a general $N\times N$ normal matrix into another normal matrix whose spectrum is supported on finitely many points, is our main technical novelty.

While the idea of tracking the eigenvalue statistics close to spectral edges along flows of deformations has been successfully implemented several times, for example
in the initial fundamental work \cite{Lee-Schnelli2015} 
on Hermitian regular edges and \cite{Campbell-Cipolloni-Erdos-Ji2024} on non-Hermitian sharp edges, the construction of the flow to study non--Hermitian critical edges is much more delicate. In fact, in all previous works 
flows of the form $\caA_{t}=c_{t}A+d_{t}I$ were used,
 where $A$ is the initial deformation and $c_{t},d_{t}$ are time--dependent scalars with $c_{0}=1$ and $c_{1}=0$. 
 Such  flows connect $A$ to a scalar multiple of the identity, which serves as a reference for the deformation.
On the other hand, 
in the current case 
the flow must go beyond such an affine family $\spann(A,B)$ even if $B$ is a reference deformation for which $B+X$ has a critical edge:
For fixed $A$ and $B$, the intersection of $\spann(A,B)$ and the set of $\caA$'s
for which $\caA+X$ has a critical edge is typically discrete, even if we allow for different values of the parameter $\alpha$.
Hence a continuous interpolation within this simple affine family $\caA_{t}=c_{t}A+d_{t}B$ of
matrices is not possible and this forces us to evolve $A$ truly as a matrix. Connecting $A$ with $B$ along a continuous 
flow while keeping both relevant parameters $\alpha, \beta$ fixed is a complicated  shooting problem
for a  matrix valued ODE-system.
For this reason,
even a slight extension of  this affine family (e.g. including more parameters but still in 
an explicit form) does not seem to  work, so we use a very different construction.

The point is that,  thanks to partial universality result in
\cite{Liu-Zhang2023}, we have more freedom on the final point of the flow $\caA_{1}$ compared to the sharp edge case, that is, $\caA_{1}$ need not be 
a fixed ($A$-independent) matrix, but
only needs to have finitely many distinct eigenvalues. The set of such matrices can be made suitably dense in the space of normal matrices 
by increasing the number (still independent of $N$)
of distinct eigenvalues. Taking advantage of this fact, our flow $\caA_{t}$ is constructed as follows: $\caA_{t}$ stays normal with the same eigenvectors as $\caA_{0}=A-z$, and the eigenvalue distribution $\rho_{\caA_{t}}$ `shrinks' into a sum of finite (but large) number of Dirac masses. The shrinkage happens locally so that only eigenvalues close to each other end up at the same Dirac mass. In particular, evolutions of eigenvalues of $A$ that are far away from each other are unrelated.
Thus the core of  our analysis is still essentially local and inverse function theorems can be used 
instead of a solving a genuine shooting problem.
 The actual construction is more involved due to the additional restriction
  on the parameter $\alpha$ and  to the discreteness of the Dirac masses. 
  We refer to the preamble of Section \ref{sec:path} for more details.

To show that the time derivatives of eigenvalues and singular values along the flow are small, we use the careful construction of $\caA_{t}$ to identify  non--trivial cancellations within the derivatives, while to estimate the error terms we rely on the single resolvent local law from \cite[Theorem 3.3]{Campbell-Cipolloni-Erdos-Ji2024}. We point out that, while we cannot treat non--normal matrices, an advantage of this approach is that we solely rely on the single resolvent local law \cite[Theorem 3.3]{Campbell-Cipolloni-Erdos-Ji2024}.

	\subsection{Notations}
	In this section we introduce some common notations used throughout the paper.
	For integers $k,n\in\Z$ we write $\bbrktt{k,n}:=[k,n]\cap\Z$, $\rrbra k,n\rrbra=(k,n]\cap\Z$, et cetera, and abbreviate $\bbrktt{k}\equiv\bbrktt{1,k}$. We denote the Lebesgue measure on $\C\cong\R^2$ by $\rmd^2 z$, and that on $\C^{k}\cong\R^{2k}$ by $\dd^{2k}\bsz$. For a square matrix $A\in\C^{d\times d}$, we write its normalized trace as $\brkt{A}\deq d^{-1}\Tr A$ and denote $\absv{A}\deq\sqrt{A\adj A}$. We use $\bsu\adj\bsy$ to denote the standard Euclidean inner product of two vectors $\bsu,\bsy\in\C^d$. Furthermore, we write $c$ and $C$ for generic positive constants (independent of $N$) whose precise values may vary by lines. For positive quantities $f,g$ we write $f\lesssim g$ and $f\sim g$ if $f\le Cg$ and $cg\le f \le Cg$, respectively. We write $\C_{+}$ for the upper half plane $\C_{+}:=\{z\in\C:\im z>0\}$.
	
	\section{Main result}
	
	We start introducing the ensemble of random matrices we consider throughout this paper, as well as some preliminary definitions needed to state our main result.
	
	\begin{defn}\label{defn:X}
		An $N$-dependent sequence of complex 
		random matrices $X\equiv X^{(N)}\in\C^{N\times N}$ 
		is called an \emph{IID random matrix} if it satisfies the following:
			\begin{itemize}
				\item The entries of $X$ are independent.
				\item $\E[X]=0$, $\E[\absv{X_{ij}}^{2}]=1/N$,
				and $\E[X_{ij}^{2}]=0$. 
				\item The entries of $\sqrt{N}X$ have finite moments, i.e. for each $p\in\N$ there exists a constant $c_{p}>0$ with
				\beq
				\label{eq:momentass}
				\sup_{N\in\N}\max_{i,j}\E\absv{\sqrt{N}X_{ij}}^{p}\leq c_{p}.
				\eeq
			\end{itemize}
			A complex IID matrix with Gaussian entries is referred to as \emph{complex} 
			\emph{Ginibre matrix} and denoted by $X^{\mathrm{Gin}(\bbC)}$.
	\end{defn}

	\begin{defn}[Brown measure, \cite{Brown1986}]\label{defn:BM}
		Let $(\caM,\brkt{\cdot}_{\caM})$ be a $W\adj$-probability space, i.e., let $\caM$ be a von-Neumann algebra with a faithful, normal, tracial state $\brkt{\cdot}_{\caM}$. For an element $\bsa\in\caM$, we define its \emph{Brown measure} to be the distributional Laplacian
		\beq
		\rho_{\bsa}\deq \frac{1}{2\pi}\Delta\brkt{\log\absv{\bsa-\cdot}}_{\caM}.
		\eeq
		For a (random) matrix $B\in\bbC^{N\times N}$, we consistently write $\rho_{B}$ for its empirical eigenvalue distribution.
	\end{defn}
	
	\begin{defn}\label{defn:vN}
		Let $(\bsx_{ij})_{i,j\in\bbrktt{N}}\in\caM$ be a collection of $*$-free circular elements\footnote{More precisely, $\bsx_{ij}=(\bss_{ij,1}+\ii\bss_{ij,2})/\sqrt{2}$ where $\{\bss_{ij,1},\bss_{ij,2}:1\leq i,j\leq N\}$ is a collection of free semi-circular elements in $\caM$.}, and define $\bsx=(\bsx_{ij})\in\caM^{N\times N}$. We canonically embed $\C^{N\times N}$ in $\caM^{N\times N}$, so that $A+\bsx$ becomes a sum of two elements of $\caM^{N\times N}$. Finally, we consistently extend $\brkt{\cdot}_{\caM}:\caM^{N\times N}\to\C^{N\times N}$ to be the partial trace, i.e.
		\beq\label{eq:vN_pt}
		\brkt{\bsy}_{\caM}=(\brkt{\bsy_{ij}}_{\caM})_{1\leq i,j\leq N}\in\C^{N\times N},\qquad \bsy=(\bsy_{ij})_{1\leq i,j\leq N}\in\caM^{N\times N}.
		\eeq
		Notice that $(\caM^{N\times N}, \brkt{\brkt{\cdot}_{\caM}})$ is a $W\adj$-probability space, where the outer bracket $\brkt{\cdot}$ is the normalized trace on $\C^{N\times N}$.
	\end{defn}

	\begin{defn}\label{defn:A}
		For a matrix $A\in\C^{N\times N}$ and $\frC>1$, we say that \emph{$A+\bsx$} has a  \emph{critical edge} or \emph{criticality at the origin} (with parameter $\frC$) if it satisfies the following:
		\begin{itemize}
			\item[(i)] $\norm{A}\leq \frC$, $\norm{A^{-1}}\leq \frC$:
			\item[(ii)] $\brkt{\absv{A}^{-2}}=1$
			\item[(iii)] $\brkt{A^{-2}(A\adj)^{-1}}=0$.
		\end{itemize}
		We write $\mathrm{Crit}_{N,\frC}$ for the set of $A$'s satisfying (i)--(iii). For $A\in\bigcup_{\frC>0}\mathrm{Crit}_{N,\frC}$, we define the {\bf shape parameter} of $A$ as
		\begin{equation}
		\label{eq:defangle}
			\alpha(A)\deq\lambda_{2}/\lambda_{1},
		\end{equation}
		where $\lambda_{1}\geq\lambda_{2}$ are the eigenvalues of the Hessian matrix $\caH\in\R^{2\times 2}$ of 
		the real-valued function $\R^{2}\ni(x,y)
		\mapsto \brkt{\absv{A-x-\ii y}^{-2}}$ evaluated at the origin. 
	\end{defn}

	\begin{rem}\label{rem:assump} We explain the role of each assumption in Definition \ref{defn:A}. The second and third assumptions guarantee that the origin is an edge point of $\rho_{A+X}$ and that it is a critical edge, respectively, so that they are part of the setup. The first part of (i) on $\norm{A}$ is technical, it comes mainly from the stability of the associated Dyson equation (see \eqref{eq:MDE} below) whose theory has traditionally been developed for bounded $A$ even though this condition could be relaxed. 
		
	The second part of (i) on $\norm{A^{-1}}$ is more important, in fact the clean dichotomy between sharp and critical edge may fail without it, see \cite[Theorem 14]{Erdos-Ji2023circ} for an example of an intermediate situation. Some spectrum of $A$ close to the origin may be allowed, but the precise condition is delicate. For example, one potential way to relax $\norm{A^{-1}}\leq\frC$ is to allow for a tiny fraction of the spectrum of $A$ (assuming that $A$ is normal) around the origin. In such a case we expect the same result to hold. However, when the ``outliers'' of $A$ are too close to the origin, they may affect the local statistics: See \cite[Theorem~1.2]{Liu-Zhang2023} for an extreme case where $A$ is allowed to have $O(1)$ number of eigenvalues \emph{exactly} at the origin.
	\end{rem}
	
    Without loss of generality, thanks to a trivial shift,  we may assume throughout the paper that the criticality occurs at the origin. We will refer to the origin as  critical edge  or criticality interchangeably. The shape parameter indeed determines the shape of the limiting spectrum $\rho_{A+\bsx}$ around the origin; see Figure \ref{fig:alpha} above. More precisely, when $\alpha\leq0$ the support asymptotically lies within the double cone $\absv{\re z}\geq\sqrt{-\alpha}\absv{\im z}$. When $\alpha>0$, the support is simply connected around the origin, but still $\alpha$ determines the density (see \eqref{eq:dens_norm} below).

	We now show some restrictions on the range of the angle $\alpha(A)$ according to properties of $A$. The proof of this lemma is postponed to Appendix~\ref{sec:aux}.
	
	\begin{lem}[Ranges of $\alpha(A)$]\label{lem:alpha}\
		\begin{itemize}
			\item[(i)] The set of possible values of $\alpha(A)$ is the interval $(-1,1]$,  i.e.
			\begin{equation}
				(-1,1]=\bigcup_{N\in\N,\frC>1}\left\{\alpha(A):A\in \mathrm{Crit}_{N,\frC}\right\}.
			\end{equation}
			\item[(ii)] If we restrict to normal $A$'s, the set of possible values of $\alpha(A)$ is $[-1/3,1]$, i.e.
			\begin{equation}
				[-1/3,1]=\bigcup_{N\in\N,\frC>1}\left\{\alpha(A):A\in\mathrm{Crit}_{N,\frC}, AA\adj=A\adj A\right\}.
			\end{equation}
			
			\item[(iii)] If we further restrict to Hermitian $A$'s, the value of $\alpha(A)$ is always $-1/3$, i.e. 			\begin{equation}
				\{-1/3\}=\bigcup_{N\in\N,\frC>1}\left\{\alpha(A):A\in\mathrm{Crit}_{N,\frC}, A=A\adj\right\}.
			\end{equation}
			Moreover, in the family of normal matrices $\alpha(A)=-1/3$ iff $A=e^{\ii\varphi}A^*$ for some $\varphi\in [0,2\pi)$.
		\end{itemize}
	\end{lem}

For a matrix $A$ having a criticality at the origin, after proper scaling, the Brown measure of $A+\bsx$ (hence the limiting 
eigenvalue density of $A+X$) near the origin depends only on $\alpha(A)$ up to the leading order. More precisely, we introduce the {\bf scaling parameter}
\begin{equation}\label{eq:gamma}
	\gamma\equiv\gamma(A)\deq \frac{2(\brkt{A^{-2}(A\adj)^{-2}})^{1/2}}{\brkt{\absv{A}^{-4}}^{1/4}}\e{\ii\theta}\equiv \frac{(\Tr \caH)^{1/2}}{\brkt{\absv{A}^{-4}}^{1/4}}\e{\ii\theta}\in\C,
\end{equation}
where $\e{\ii\theta}\in\C\cong\R^{2}$ is the eigenvector of $\caH$ corresponding to the larger eigenvalue\footnote{\label{fn:alpha=0} When $\lambda_{1}=\lambda_{2}$, i.e. $\alpha(A)=1$, the Hessian is a constant multiple of the identity. In this case we choose $\theta=0$ for definiteness, yet this choice does not affect \eqref{eq:dens_norm} since $\rho_{A+\bsx}$ is radially symmetric up to leading order.} 
$\lambda_{1}$, i.e. 
\begin{equation}\label{eq:theta}
	\caH\begin{pmatrix}\cos\theta \\ \sin\theta\end{pmatrix}
=\lambda_{1}\begin{pmatrix}\cos\theta \\ \sin\theta \end{pmatrix}.
\end{equation}
Now with the change of variable
 $z\mapsto \gamma^{-1}z$, by  \cite[Eq. (2.14)]{Erdos-Ji2023circ}, we find that 
\begin{equation}\label{eq:dens_norm}
	\rho_{\gamma(A+\bsx)}(z)=\frac{1}{\absv{\gamma}^{2}}\rho_{A+\bsx}(\gamma^{-1} z)
	=\frac{\lone(\absv{A-\gamma^{-1}z}^{-2}\geq 1)}{8\pi}
	\left[\frac{x^{2}+\alpha(A) y^{2}}{1+\alpha(A)}
	+2\frac{x^{2}+\alpha(A)^{2}y^{2}}{(1+\alpha(A))^{2}}+O(\absv{z}^{3})\right].
\end{equation}
Note that the quadratic form $(x,y)\caH(x,y)\tp$ is nonnegative up to leading order when $x+\ii y\in\supp \rho_{A+\bsx}$. Indeed, the support is exactly the super-level set $\{z:\brkt{\absv{A-z}^{-2}}\geq 1\}$ by \cite[Theorem 4.6]{Zhong2021}, and 
at criticality the function $\brkt{\absv{A-z}^{-2}}$ has value $1$, vanishing first derivatives, and Hessian $\caH$ at the origin.

\subsection{Main result}

The main result of this paper is the universality of local statistics close to criticality, when the deformation $A$ is normal (see Theorem~\ref{theo:mainthm} below). However, there are further interesting examples and phenomena appearing close to criticality that we will discuss in Section~\ref{sec:extension}. We now state our main result:
\begin{thm}
\label{theo:mainthm}
	Let $X=X^{(N)}$ be an $N\times N$ IID matrix and let $A=A^{(N)}$ be an $N\times N$ normal matrix having a criticality at the origin with parameter $\frC>1$. Assume further that there are $N$-independent constants $\alpha_{\infty}\in (-1/3,1]$ and $\frc>0$ such that
	\begin{equation}\label{eq:alpha_assu}
			\absv{\alpha(A)-\alpha_{\infty}}\leq N^{-\frc}.
	\end{equation} 
	Fix $k\in\N$ and let $p_{k}^{(N)}$ be the $k$-point correlation function of the eigenvalues of $N^{1/4}\gamma(A)(A+X)$ near the origin. Then, for any $F\in C_{c}^{\infty}(\C^{k})$, we have
	\begin{equation}
		\lim_{N\to\infty}\int_{\C^{k}} F(\bsz)(p_{k}^{(N)}(\bsz)-p_{k,\alpha_{\infty}}(\bsz))\dd^{2k}\bsz=0,
	\end{equation}
	where $p_{k,\alpha}:\C^{k}\to\R$ is an explicit symmetric function depending only on the shape parameter $\alpha$ (see Remark \ref{rem:Liu}).
\end{thm}
\begin{rem}\label{rem:Liu}
We point out that Liu and Zhang in \cite{Liu-Zhang2023} computed an explicit formula for $p_{k,\alpha}$ as the limit of $p_{k}^{(N)}$ when $X$ is a 
 complex Ginibre matrix and $A$ is a normal matrix of a quite
 specific form (roughly speaking $A$ has only a few, independent of $N$, eigenvalues that are
 $N$-independent themselves and have fixed relative multiplicities; see \cite[Eqs. (1.3)-(1.5)]{Liu-Zhang2023}).  We refer the interested reader to \cite[Theorem 1.2]{Liu-Zhang2023} for the explicit formula for $p_{k,\alpha}$, we do not repeat it here since it is rather lengthy and involved containing a $k$--dependent number contour integrals. In particular, $p_{k,\alpha}$ does not seem to have a determinental structure so there is a separate formula for each $k$. Our Theorem~\ref{theo:mainthm} substantially extends this universality result in two directions: i) we consider general IID $X$, ii) the deformation $A$ is required to satisfy only Definition~\ref{defn:A}, normality, and \eqref{eq:alpha_assu}.
\end{rem}

\begin{rem}[On the conditions \eqref{eq:alpha_assu}]\label{rem:Liu_gen}
We believe that this condition can be relaxed merely to the existence of the limit $\alpha(A)\to \alpha_\infty$
without any effective control. The reason for the current \eqref{eq:alpha_assu} lies in the non-optimal
conditions in~\cite{Liu-Zhang2023}. Note that 
	our proof follows a flow\footnote{Throughout the paper, we consistently use calligraphic letters (e.g. $\caA, \caB, \caV$) to denote the flows that connect  objects denoted by the corresponding Roman letters (e.g. $A, B, V$).}
	$\caA_{t}+X$, where the deformation $\caA_{t}$ connects a given normal matrix 
	$A=\caA_{0}$ to another normal matrix $\caA_{1}$ with known correlation functions, while keeping the value of $\alpha(A_{t})$ (almost) constant. Since the limiting correlation functions at a critical edge are only available from \cite{Liu-Zhang2023}, the final point $\caA_{1}$ must be covered by this work. In particular,  the
	conditions of~\cite{Liu-Zhang2023}  demand that 
	the spectral distribution of  $\caA_{1}$  be essentially $N$-independent (up to $O(N^{-1})$).
	This implies that $\alpha(\caA_{1}) = \alpha_\infty + O(N^{-1})$ with some $N$-independent $\alpha_\infty$.
	In fact, with an extra twist, our flow method allows the initial value $\alpha(A)$ for a much bigger $O(N^{-\frc})$ room
	in~\eqref{eq:alpha_assu}
	than the straightforward $O(N^{-1})$ approximation error from~\cite{Liu-Zhang2023} would predict.
	
	Nonetheless, a natural refinement of the  method of~\cite{Liu-Zhang2023} 
	should  allow $\caA_{1}$ to depend on $N$, as long as  $\caA_{1}\in\mathrm{Crit}_{N,\frC}$,
	   the number of distinct eigenvalues of $\caA_{1}$ remains bounded, and  $\alpha(\caA_{1})$ has a limit. 
	   The current  proof in~\cite{Liu-Zhang2023} uses a
	  completely $N$-independent saddle point analysis  even though the phase functions have mild $N$-dependence.
	   The restrictive conditions in~\cite{Liu-Zhang2023}  stem from crudely estimating this difference.  
	  We believe that, with some extra work, a more accurate analysis can be performed by expanding around
	  the exact $N$-dependent saddle points, 
	   for example, keeping the location of the saddle and the Hessian at the saddle  $N$-dependent. 
	   In private communication with D.-Z. Liu, he agreed with this assessment.
	Given such a generalization of \cite{Liu-Zhang2023}, we can easily replace
	 the assumption \eqref{eq:alpha_assu} with the minimal condition $\alpha(A)\to \alpha_\infty$.
	 \end{rem}

\begin{rem}[On the condition $\alpha_{\infty}>-1/3$]
	According to Lemma \ref{lem:alpha} (iii), by assuming $\alpha_{\infty}>-1/3$ we have excluded Hermitian $A$ in Theorem \ref{theo:mainthm}. Our approach for non-Hermitian $A$ does not directly 
		apply to Hermitian $A$ for two reasons: First, if the initial matrix $A$ is Hermitian (i.e. $\alpha(A)=-1/3$), then the whole path $\caA_{t}$ must stay Hermitian in order to keep $\caA_{t}$ normal and $\alpha(\caA_{t})=-1/3$ constant by Lemma~\ref{lem:alpha}~(iii). On the other hand, along the flow the eigenvalues of $\caA_{t}$ are not allowed to `cross' the origin due to the regularity condition $\norm{A^{-1}}=O(1)$. Thus the ($N$-dependent) numbers of positive and negative eigenvalues of $\caA_{t}$ remain unchanged throughout the flow. As discussed in Remark \ref{rem:Liu_gen}, this prohibits us from connecting to models covered in \cite{Liu-Zhang2023} when these numbers heavily depend on $N$, for example if $\absv{\rho_{A}(0,\infty)-\rho_{A}(-\infty,0)}=1/\log N$ where $\rho_{A}$ is the spectral measure. Nonetheless, in Appendix~\ref{sec:Herm_A} we prove that any two Hermitian matrices in $\mathrm{Crit}_{N,\frC}$ can be connected as long as their numbers of positive eigenvalues are the same. Hence we can cover fully general Hermitian $A$'s once we have an $N$-dependent generalization of \cite{Liu-Zhang2023} as explained in Remark \ref{rem:Liu_gen}. Alternatively, assuming that the condition on $\lVert A^{-1}\rVert$ in Definition \ref{defn:A} (i) can be relaxed as described in Remark \ref{rem:assump} so that some eigenvalues of $A$ can be $o(1)$, then one may first calibrate the numbers of positive and negative eigenvalues by moving some some of them to the other side following a path $o(1)$-away from the real axis.
	
	The second reason is that our construction of the path $\caA_{t}$ is tailored for non-Hermitian (still normal) deformations; see Lemma \ref{lem:Jac}. The current approach capitalizes on the fact that, given $\alpha>-1/3$, the set of matrices $\{A:\alpha(A)=\alpha\}$ is of (real) co-dimension $2$ in the space of normal matrices. In contrast, as shown in Lemma \ref{lem:alpha}, the equation $\alpha(A)=-1/3$ is critical in that the co-dimension is $N$. Indeed, the path for Hermitian $A$ in Section \ref{sec:Herm_A} is completely different from that for non-Hermitian $A$ in Section \ref{sec:path}.
\end{rem}

\subsection{Examples and extensions}
\label{sec:extension}

This section is divided into two parts: in Section~\ref{sec:almostcrit} we describe some possible extensions of Theorem~\ref{theo:mainthm}, while in Section~\ref{sec:exampl} we present some concrete examples illustrating Lemma~\ref{lem:alpha}.

\subsubsection{Almost criticality}
\label{sec:almostcrit}

\begin{figure}
	\centering
	\includegraphics[width=0.32\textwidth]{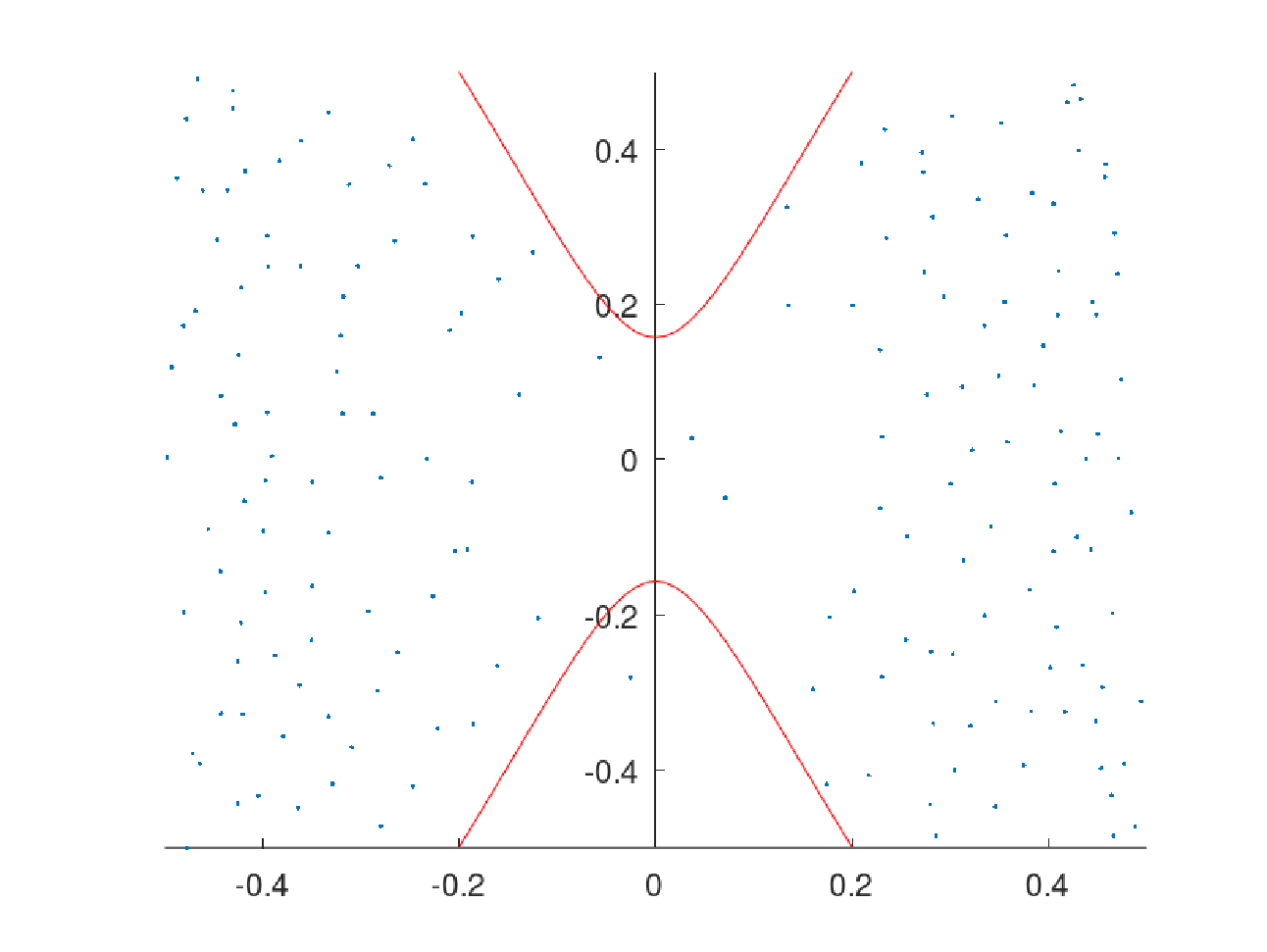}
	\includegraphics[width=0.32\textwidth]{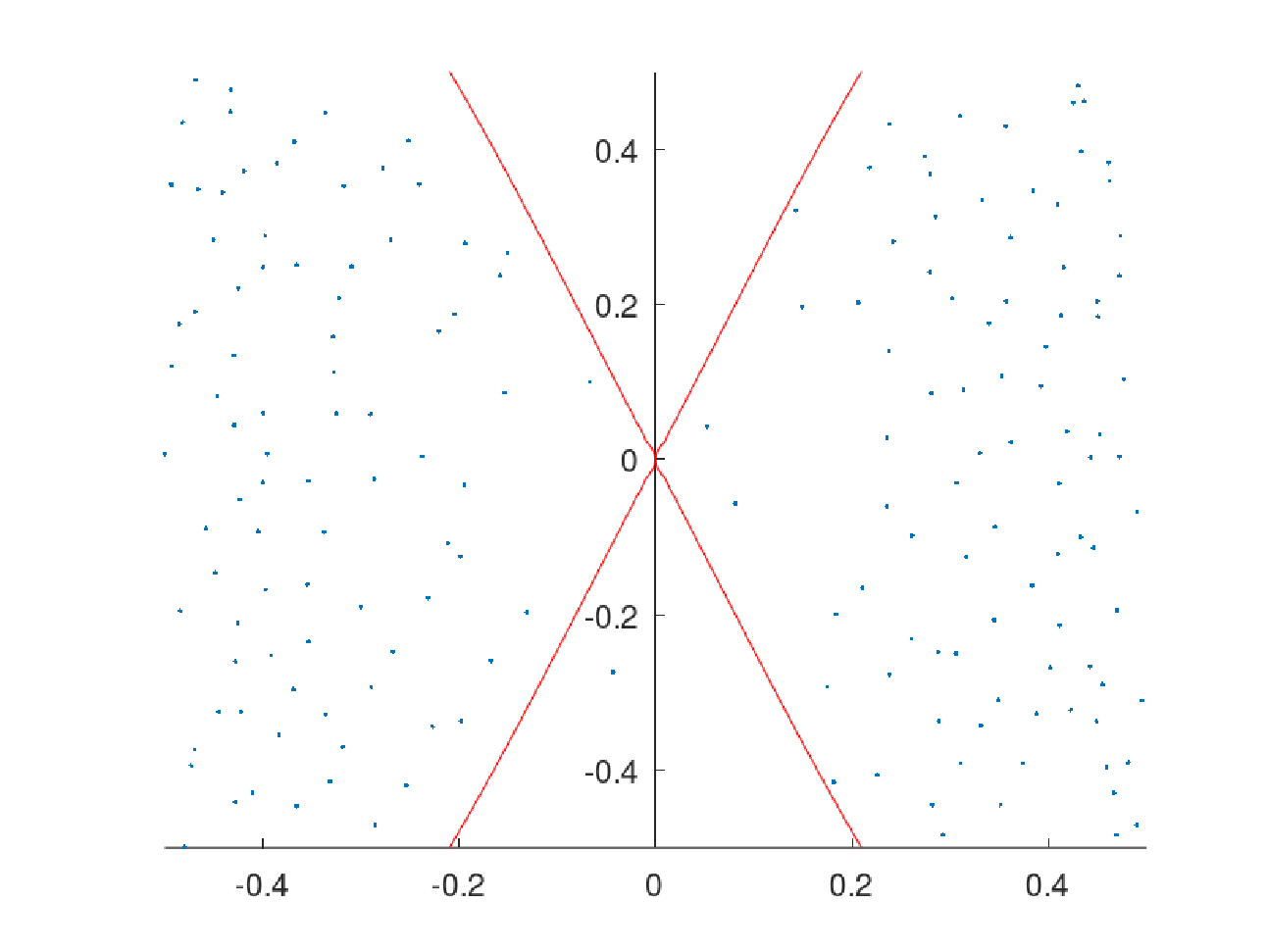}
	\includegraphics[width=0.32\textwidth]{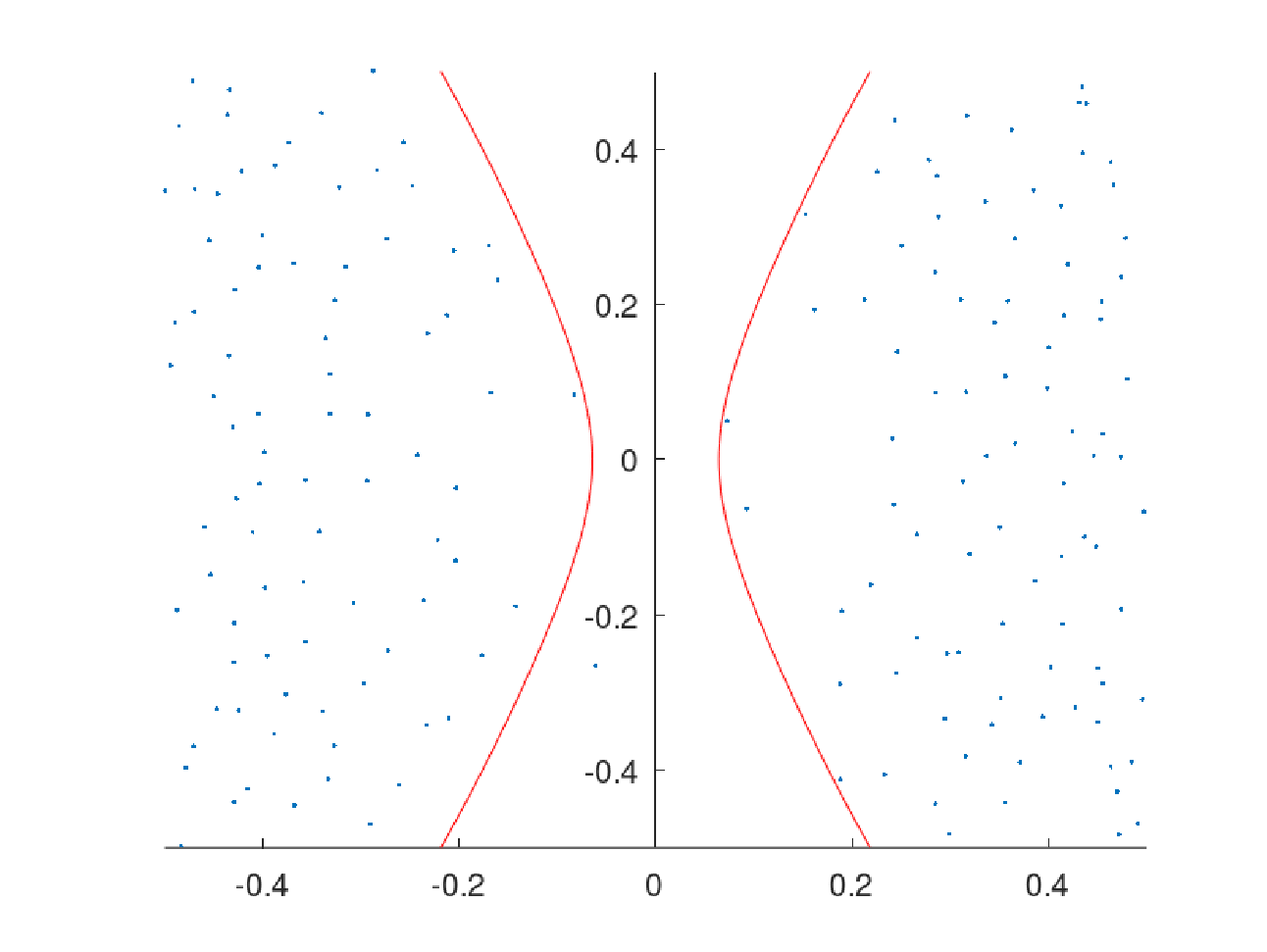}
	\caption{Almost criticalities at the origin for $\beta<0$, $\beta=0$, and $\beta>0$ (from left to right)}\label{fig:beta}
\end{figure}

In Theorem~\ref{theo:mainthm} we stated the universality result when $A+\bsx$ has an exact criticality at the origin, showing that the local fluctuation scale of the eigenvalues close to the critical edge\footnote{We point out that in the special case $\alpha(A)=0$ new statistics, with a new anisotropic scaling, may emerge in the special direction $\ii y$. However, the study of this atypical statistics is not in the scope of the current work.} is of order $N^{-1/4}$, in contrast with the $N^{-1/2}$--scale close to sharp edges (see e.g. \cite[Theorem~2.7]{Campbell-Cipolloni-Erdos-Ji2024}). It is thus natural to ask if it is possible to prove an analog of Theorem~\ref{theo:mainthm} in the regime of \emph{almost criticality}, i.e. when the density of $A+\bsx$ does not have an exact criticality at zero but  the support of the density locally around zero consists of two different component at a distance $N^{-1/4}$ or of one component whose length in one direction is of order $N^{-1/4}$. 
Our proof is strong enough to also prove universality in the regime of \emph{almost criticality} monitoring an additional parameter
\begin{equation}
\label{eq:defbneta}
\beta(A):=N^{1/2}\big(1-\langle |A-z|^{-2}\rangle\big),
\end{equation}
giving rise to a two parameter universal family of local statistics (see e.g. \cite[Theorem~1.2]{Liu-Zhang2023}, where the parameter $\beta$ is denoted as $\wh{\tau}$). In particular, the parameter $\beta(A)$ monitors the distance between the two components of the spectrum in this regime of \emph{almost criticality}. In fact, precisely at \emph{criticality} we would have $\beta(A)=0$, while for $\beta(A)>0$ the density has a gap of order $N^{-1/4}$, and the parameter $\beta(A)$ describes its size. We omit the details for brevity and only provide a pictorial illustration in Figure \ref{fig:beta}. Furthermore, we expect that if $\rho_{A+\bsx}$ consists of two different component at a distance $\gg N^{-1/4}$ around the origin then the local statistics close to edges of these two components are (asymptotically) the same as those close to a sharp edge of a complex Ginibre matrix. We believe that this fact is also provable using the methods developed in this paper, we omit this to keep the presentation short and simple. Lastly, one would expect that in the regime when $\rho_{A+\bsx}$, locally around zero, consists of one component whose length in one direction is of order $\gg N^{-1/4}$ then the local statistics around zero would asymptotically agree with those in the bulk of the complex Ginibre ensemble; this is out of the scope of this paper since we are interested in the edge behavior of $A+X$.

\subsubsection{Examples describing Lemma~\ref{lem:alpha}}
\label{sec:exampl}

\begin{figure}
	\includegraphics[width=0.5\textwidth]{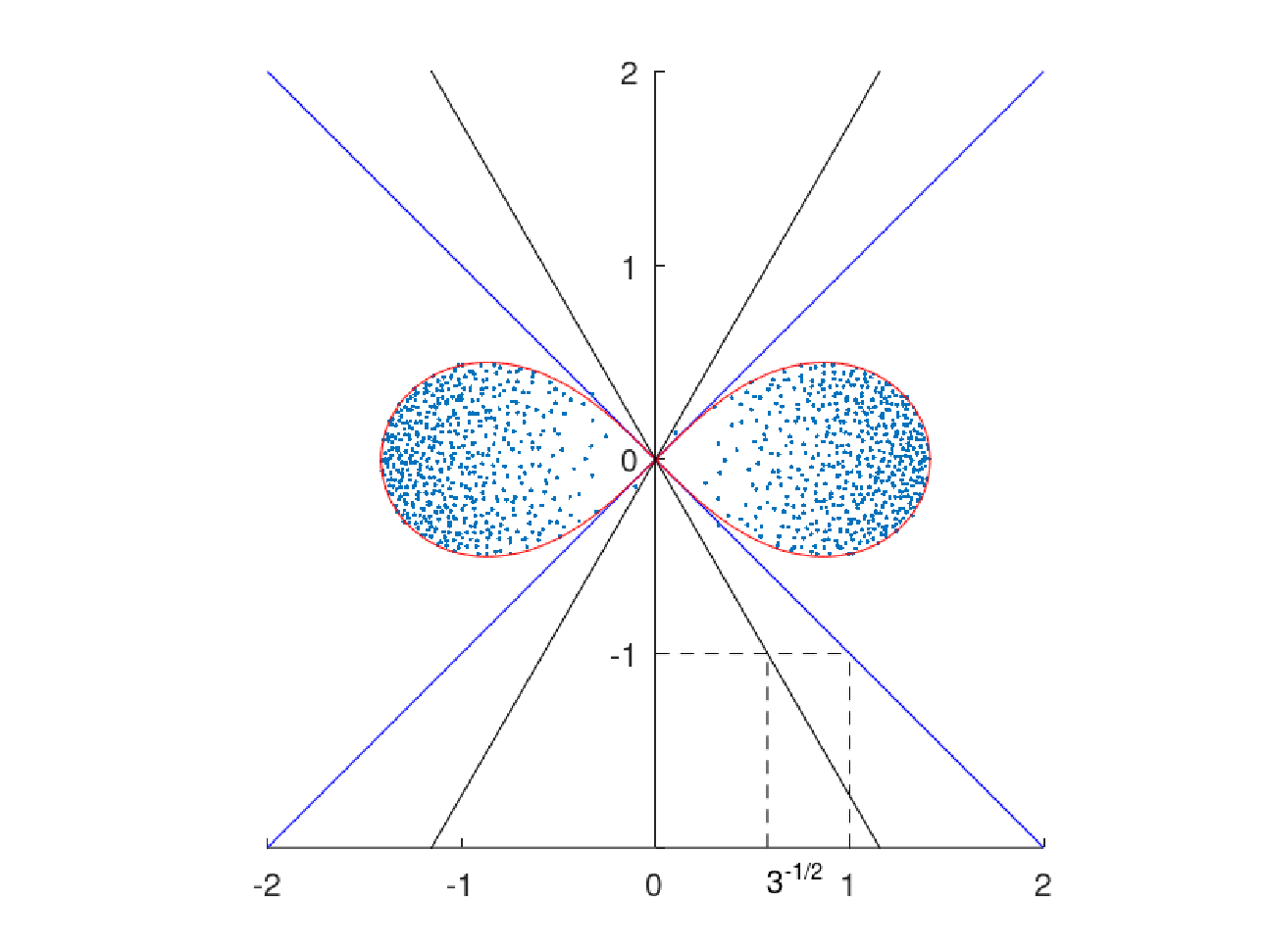}
	\caption{$\spec(A+X)$ for $A$ as in \eqref{eq:shapeparex1} with $c\gg 1$; black lines are of slope $\sqrt{3}$, corresponding to the lower bound $\alpha\geq-1/3$ for normal $A$.}\label{fig:nonnorm}
\end{figure}

First, we give an example of a family of non--normal matrices $A$ for which the shape parameter $\alpha(A)$ from \eqref{eq:defangle} goes below $-1/3$, in fact we will see that it can be arbitrarily close to $-1$. Consider the $N\times N$ deterministic matrix
\begin{equation}\label{eq:shapeparex1}
A=A(c):=(1+c^{2}/2)^{-1/2}\begin{pmatrix} -1 & c \\ 0 & 1 \end{pmatrix}^{\oplus N/2},\qquad\quad c>0.
\end{equation}
Then we can easily see that $A$ has a criticality at the origin, and by explicit computations we get
\begin{equation}
\label{eq:shapeparex}
\alpha(A)=\frac{\mathcal{H}_{22}}{\mathcal{H}_{11}}=-\frac{2+2c^2}{6+2c^2},
\end{equation}
where by $\mathcal{H}_{ij}$ we denoted the $(i,j)$--entry of the $2\times 2$ Hessian matrix $\mathcal{H}$ from Definition~\ref{defn:A}. Note that as $c\in \R_+$ varies $\alpha(A)$ covers the whole interval $(-1,-1/3]$. See Figure \ref{fig:nonnorm} for an illustration.

Next, we give an example to show that the non--normality of $A$ does not necessarily decrease the value of $\alpha(A)$, i.e. we show that there are also non--normal matrices for which $\alpha(A)> -1/3$. Consider
\begin{equation}
A:=(1+c^{2}/2)^{-1/2}\left(A_1\oplus A_2\right)^{\oplus N/4},
\quad\quad
A_1=A_1(c):=\begin{pmatrix} -1 & c \\ 0 & -1 \end{pmatrix},
\quad\quad
A_2=A_2(c):=\begin{pmatrix} 1 & c \\ 0 & 1 \end{pmatrix},
\quad\quad
c>0.
\end{equation}
By explicit computations we thus find that
\begin{equation}
\alpha(A)=-\frac{2+4c^2}{6+20c^2}.
\end{equation}
As $c\in \R_+$ varies $\alpha(A)$ covers the interval $[-1/3, -1/5)$, showing that for non--normal matrices it is possible to have $\alpha(A)>-1/3$.

\section{Universality along deformation path}
In this section, we introduce several notations used throughout the comparison of the eigenvalue statistics of $\caA_{t}+X$, where $\caA_{t}$ is a flow of deformations with a parameter $\alpha_t\equiv\alpha(\mathcal{A}_{t})$ which changes only very mildly as $t$ varies. Most notations are consistent with the ones previously introduced in \cite{Campbell-Cipolloni-Erdos-Ji2024} where the sharp edge statistics were considered. Practically, 
 the only difference is that the $w$-variable in this paper (see \eqref{eq:def_gamma} and \eqref{eq:def_Herm} below) carries a different scaling that the one used in \cite{Campbell-Cipolloni-Erdos-Ji2024} (see (4.37) and (4.38) therein). 
 The main technical result of this section is then a theorem showing that the local statistics close to a criticality of $\caA_{t}+X$ are preserved along the flow (see Theorem~\ref{theo:mainthmflow} below). 
\begin{assump}[Deformation path]\label{assump:path}
	Let $\frC_{1}>1$, let $\frc_{1}>0$ be a small constant and let $[0,1]\ni t\mapsto \caA_{t}\in\C^{N\times N}$ be a $C^{1}$ path that satisfies the following for all $t\in[0,1]$:
	\begin{itemize}
		\item[(i)] $\caA_{t}+\bsx$ has a criticality at the origin with parameter $\frC_{1}>0$:
		\item[(ii)] $\absv{\dd\alpha(\caA_{t})/\dd t}\leq N^{-\frc_{1}}$.
		\item[(iii)] $\norm{\frac{\dd \caA_{t}}{\dd t}}\leq \frC_{1}\log N$.
	\end{itemize}
\end{assump}

For such a choice of $\caA_{t}$, we will show that the local eigenvalue statistics near the criticality of $\caA_{t}+X$ is independent of $t$ in the large $N$ limit, after proper rescaling. 
Introducing the time-dependent auxiliary quantities $I_{4},\wt{I}_{4}$:
\begin{equation}\label{eq:def_I}
I_{4}\equiv I_{4}(t)\deq \Brkt{\frac{1}{\absv{\caA_{t}}^{4}}},\qquad 
\wt{I}_{4}\equiv \wt{I}_{4}(t)\deq \Brkt{\frac{1}{\caA_{t}^{2}(\caA_{t}\adj)^{2}}},
\end{equation}
the scaling parameters of our flow, $c_{t}$ and $\gamma_{t}$, are defined as
\begin{equation}\label{eq:def_gamma}
c_{t}\deq I_{4}(t)^{-1/4},\qquad  \gamma_{t}\deq I_{4}(t)^{-1/4}\wt{I}_{4}(t)^{1/2}\e{\ii\theta_{t}},
\end{equation}
where $e^{\ii\theta_t}$ is the eigenvector of the Hessian $\caH$ of $\brkt{\absv{\caA_{t}-z}^{-2}}$ at the origin as in \eqref{eq:theta}. The parameter $\gamma_{t}$ serves the same purpose as $\gamma$ in \eqref{eq:gamma}, i.e., rescales the non-Hermitian eigenvalue density, and $c_{t}$ rescales the singular value density of $A+X$, so that both densities depend only on $\alpha(A)$ (see \eqref{eq:cubic_nor} below).

We now state our main technical result: if two deformation matrices can be connected through a flow satisfying  Assumption~\ref{assump:path} then their local statistics close to criticality are universal.
\begin{thm}[Universality along the flow]
\label{theo:mainthmflow}
	Let $X$ be an IID matrix and, for $t\in [0,1]$, let $\mathcal{A}_t\in\C^{N\times N}$ as in Assumption~\ref{assump:path}. 
	Fix $k\in\N$ and let $p_{k,t}^{(N)}$ be the $k$-point correlation function of the eigenvalues of $N^{1/4}\gamma_t(\mathcal{A}_t+X)$ near the origin. Then, for any $F\in C_{c}^{\infty}(\C^{k})$, we have
	\begin{equation}
		\int_{\C^{k}} F(\bsz)\big(p_{k,t}^{(N)}(\bsz)-p_{k,0}^{(N)}(\bsz)\big)\dd^{2k}\bsz=O(N^{-c})
	\end{equation}
for some constant $c>0$, uniformly in $t\in [0,1]$.
\end{thm}

We now introduce some notations that will be used in the proof of Theorem~\ref{theo:mainthmflow}. For $w\in\C$ we introduce the Hermitizations and the regularizing parameter $\eta_{t}$
\begin{equation}\begin{aligned}\label{eq:def_Herm}
	W&\deq \begin{pmatrix}
		0 & X \\ 
		X\adj & 0
	\end{pmatrix}, \\
	H_{t}^{w}&\deq\begin{pmatrix}
		0 & \caA_{t}+X-\gamma_{t}^{-1}N^{-1/4}w \\
		(\caA_{t}+X-\gamma_{t}^{-1}N^{-1/4}w)^{\adj} & 0 
	\end{pmatrix},	\\
	\bsh_{t}^{w}&\deq\begin{pmatrix}
		0 & \caA_{t}+\bsx-\gamma_{t}^{-1}N^{-1/4}w \\
		(\caA_{t}+\bsx-\gamma_{t}^{-1}N^{-1/4}w)^{\adj} & 0 
	\end{pmatrix},\\
	J_{t}^{w}&\deq \begin{pmatrix}
		0 & \gamma_{t}^{-1}w \\
		\ol{\gamma}_{t}^{-1}\ol{w} & 0
	\end{pmatrix},	\\
	\eta_{t}&\deq c_{t}^{-1}\eta_{\infty},\qquad \eta_{\infty}\deq N^{-3/4-\delta},
\end{aligned}\end{equation}
for some small $\delta>0$. We denote the resolvent of $H_{t}^{w}$ (restricted to the imaginary axis) by
\begin{equation}\label{Gtw}
	G_{t}^{w}(\ii\eta)=(H_{t}^{w}-\ii\eta)^{-1},\qquad \eta>0.
\end{equation}
By the local law in Theorem~\ref{theo:llaw} below, it follows that the resolvent $G_{t}^{w}(\ii\eta)$ for $\eta>0$ concentrates around the deterministic matrix
\begin{equation}\label{eq:def_M}
M_{t}^{w}(\ii\eta):=\brkt{(\bsh_{t}^{w}-\ii\eta)^{-1}}_{\caM}\in\C^{2N\times 2N}.
\end{equation}
The matrix $M_{t}^{w}(\ii\eta)$ can be characterized in terms of a Matrix Dyson Equation (MDE) which we explain in detail in the next section.

\subsection{Matrix Dyson equation around criticality}
In this section we present some known properties of the deterministic approximation of the Hermitization of $A+X$.

\begin{lem}[{\cite[Theorem 2.1]{Helton-Rashidi-Speicher2007} and \cite[Theorem 3.8]{Zhong2021}}] \label{lem:Msoleq}
	Let $A\in\C^{N\times N}$ and $\wh{z}\in\C_{+}$. Then the $(2N\times 2N)$ matrix Dyson equation
	\begin{equation}\label{eq:MDE}
		\frac{1}{M}=\begin{pmatrix} 0 & A \\ A\adj & 0 \end{pmatrix} -\wh{z}-\brkt{M}
	\end{equation}
	has a unique solution $M\equiv M_{A}(\wh{z})\in\C^{N\times N}$ with $\im M>0$. Furthermore, the solution of \eqref{eq:MDE}
	can be expressed as
	\begin{equation}
		\label{eq:detM}
		M_{A}(\wh{z}) = \Brkt{\begin{pmatrix} -\wh{z} & A+\bsx \\ (A+\bsx)\adj & -\wh{z}
			\end{pmatrix}^{-1}}_{\caM}.
	\end{equation}
\end{lem}

In the following we will often use $M_{A-z}$ to denote the solution to \eqref{eq:MDE} with $A$ replaced with $A-z$. As an immediate consequence of Lemma \ref{lem:Msoleq}, we find that the deterministic matrix $M_{t}^{w}(\ii\eta)$ defined in \eqref{eq:def_M} can be written as $M_{t}^{w}(\ii\eta)=M_{\caA_{t}-\gamma_{t}^{-1}N^{-1/4}w}(\ii\eta)$, or equivalently, as a solution of
\begin{equation}
	\frac{1}{M_{t}^{w}(\ii\eta)}=\E H_{t}^{w}-\ii\eta-\brkt{M_{t}^{w}(\ii\eta)}.
\end{equation}
When $A$ has a criticality, the solution $M$ to the MDE has the following asymptotics on the imaginary axis
$\wh z = \ii\eta$:
\begin{lem}
	\label{lem:asm}
	Let $C>0$ be any fixed constant. The following holds uniformly over $A\in\C^{N\times N}$, $|z|\le C$, and $\eta>0$ with $\norm{A}\leq C$, $\norm{(A-z)^{-1}}\leq C$, and $\eta\in[0,C]$: 
	\begin{equation}\label{eq:M_asymp}
		\im\brkt{M_{A-z}(\ii\eta)}\sim\begin{cases}
			(\brkt{\absv{A-z}^{-2}}-1)+\eta^{1/3}, & \brkt{\absv{A-z}^{-2}}\geq 1, \\
			\dfrac{\eta}{1-\brkt{\absv{A-z}^{-2}}+\eta^{2/3}}, & \brkt{\absv{A-z}^{-2}}<1.
		\end{cases}
	\end{equation}
	Furthermore, if $A$ has a criticality at the origin, then
	$v \deq \im\brkt{M_{A-z}(\ii\eta)}+\eta$
	satisfies the following asymptotic cubic equation in the small $|z|, \eta$--regime.
	\begin{equation}\label{eq:cubic}
		\Brkt{\frac{1}{\absv{A}^{4}}}v^{3}-\frac{1}{2}\begin{pmatrix} \re z \\ \im z \end{pmatrix}\tp\caH\begin{pmatrix} \re z \\ \im z \end{pmatrix}v-\eta 
		=O\big(\absv{z}^{4}+\eta^{4/3}\big).
	\end{equation}
\end{lem}

We omit the proof since it is nearly identical to that of \cite[Lemma 3.4]{Campbell-Cipolloni-Erdos-Ji2024}. The only difference in the statement\footnote{Note that the ``base point" $z\in\C$ in \cite{Campbell-Cipolloni-Erdos-Ji2024} has now become the origin.} is that the error  in \eqref{eq:cubic} is slightly bigger, as a natural byproduct of the proof; here the error involves $\absv{z}^{4}+\eta^{4/3}$ whereas the corresponding term in \cite[Eq. (3.10)]{Campbell-Cipolloni-Erdos-Ji2024} is $\absv{z}^{5/2}+\eta^{5/3}$. Nonetheless, in the regime of our interest, the cubic equation \eqref{eq:cubic} is still effective in the sense that
the error is much smaller than $v^{3}$ hence it serves its purpose. Indeed, we almost always work in the regime $\absv{z}=O(N^{-1/4})$ and $\eta\sim N^{-3/4-o(1)}$, in which case \eqref{eq:M_asymp} implies $v\gtrsim N^{-1/4-o(1)}$ (note also that the natural scale for $z$ in \cite{Campbell-Cipolloni-Erdos-Ji2024} was $\absv{z}\sim N^{-1/2}$).

We now state the local law for the Hermitization of $A+X$. To keep the presentation simple and short, we state the following local law only in the average case and close to a criticality. We point out that in \cite[Theorem 3.3]{Campbell-Cipolloni-Erdos-Ji2024} we actually proved both averaged and isotropic local laws in an order one neighborhood of the edge of the spectrum of $A+X$.

\begin{thm}[Theorem 3.3 of \cite{Campbell-Cipolloni-Erdos-Ji2024}]
\label{theo:llaw}
Let $X$ be an IID matrix and $A$ a deterministic matrix as in Definition \ref{defn:A}. Let $G^z(\ii\eta):=(H^z-\ii\eta)^{-1}$, with $H^z$ being the Hermitization of $A+X-z$ defined as in \eqref{eq:hermintro},
and let $M_{A-z}$ be the solution of \eqref{eq:MDE} with $A$ replaced with $A-z$. There exists a (small enough) constant $c>0$ such that for any deterministic matrix $B\in\C^{2N\times 2N}$ we have
	\begin{equation}\label{eq:locallawgoal}
		\absv{ \brkt{ B(G^{z}(\ii\eta)-M_{A-z}(\ii\eta)) } }\prec \norm{B}\frac{1}{N\eta},
	\end{equation}
	uniformly in $A$ from Definition \ref{defn:A}, $|z|<1+c$, and $N^{-1}\leq\eta\leq N^{100}$.
\end{thm}

With the rescaling (for $w$ and $\eta$) incorporated in $M_{t}^{w}(\ii\eta_{t})$, the cubic equation \eqref{eq:cubic} is equivalent to
\begin{equation}\label{eq:cubic_nor}
	(I_{4}(t)^{1/4}v_{t}^{w})^{3}-N^{-1/2}\frac{1}{2}\frac{(\re w)^{2}+\alpha (\im w)^{2}}{1+\alpha}(I_{4}(t)^{1/4}v_{t}^{w})-\eta_{\infty}=O(N^{-1}),
\end{equation}
where $v_{t}^{w}:=\im\brkt{M_{t}^{w}(\ii\eta_{t})}+\eta_{t}$ and $\alpha=\alpha(A)$. 
 The implicit constant in the $O(N^{-1})$ error is uniform in $t\in [0,1]$. Notice that the rescaled solution $I_{4}(t)^{1/4}v_{t}^{w}$ is independent of $t$ except for the error term, and that the equation \eqref{eq:cubic_nor}, hence its rescaled 
 solution, depends on $A$ only via $\alpha(A)$.
For simplicity, we introduce the solution $\mr{M}_{t}^{w}$ to the degenerate MDE
\begin{equation}
	\mr{M}_{t}^{w}:=(\E H_{t}^{w})^{-1}.
\end{equation}

We finally introduce
\begin{equation}\label{eq:def_Lt}\begin{aligned}
	L_{t}(w)
	&\deq \Tr \log(\absv{H_{t}^{w}-\ii\eta_{t}}) -\Tr\brkt{\log(\absv{\bsh_{t}^{w}-\ii\eta_{t}})}_{\caM},
\end{aligned}\end{equation}
where $\brkt{\cdot}_{\caM}$ is applied entrywise as described in Definition \ref{defn:vN}. Notice that $L_{t}(w)$ can be expressed as an integral of $\brkt{G_{t}^{w}(\ii\eta)-M_{t}^{w}(\ii\eta)}$ in $\eta$, i.e.
\begin{equation}\label{eq:dlog=G}
	L_{t}(w)=2N\int_{\eta_{t}}^{\infty}\im\brkt{G_{t}^{w}(\ii\eta)-M_{t}^{w}(\ii\eta)}\dd\eta.
\end{equation}
Lastly, we introduce the following abbreviation for $\bsw\in\C^{k}$:
\begin{equation}
\label{eq:defL}
	\bsL_{t}(\bsw)\deq \prod_{j=1}^{k}L_{t}(w_{j}).
\end{equation}

We first study the evolution of the deterministic part in the definition of $L_t(w)$ in \eqref{eq:def_Lt} (see Lemma~\ref{lem:determ} below), showing that (at leading order) it is given by a harmonic function. Then, in Proposition~\ref{prop:univ_Gin}, we show that $t\mapsto \bsL_{t}(\bsw)$ is asymptotically constant (modulo a harmonic function). The proof of both these facts is presented in Section~\ref{sec:techest} below.

\begin{lem}\label{lem:determ}
For $\absv{w}\leq cN^{1/4}$, we have
	\begin{equation}
	\frac{\dd}{\dd t}\brkt{\brkt{\log\absv{\bsh_{t}^{w}-\ii\eta_{t}}}_{\caM}}=-\frac{\dd}{\dd t}\brkt{\log\absv{\mr{M}_{t}^{w}}}+O(N^{-1-\mathfrak{c}_1})
	\end{equation}
	uniformly in $t\in [0,1]$. Moreover, the function $w\mapsto\brkt{\log\absv{\mr{M}_{t}^{w}}}$ and its time derivatives are harmonic in $\absv{w}\leq cN^{1/4}$. Here $c>0$ is a constant   depending only on $\frC_{1}$ from Assumption \ref{assump:path}.
\end{lem}

\begin{prop}
\label{prop:univ_Gin}
	Let $A_{t}$ satisfy Assumption \ref{assump:path} and fix $k\in\N,C>0$. Then, there exists a constant $c>0$ such that following holds uniformly over $t\in[0,1]$, $w\in\C$, and $\bsw\in\C^{k}$ with $\absv{w},\absv{\bsw}\leq C$:
\begin{align}
\frac{\dd}{\dd t}\E \bsL_{t}(\bsw)=h_t(\bsw)+O(N^{-c}), \label{eq:Lbound}
\end{align}
with $\Delta_{w_1}\dots\Delta_{w_k} h_t=0$.
\end{prop}

We are now ready to prove the main result of this section:
\begin{proof}[Proof of Theorem~\ref{theo:mainthmflow}]

We now show that if two matrices $\mathcal{A}_s$, $\mathcal{A}_t$ are connected by a path as in Assumption~\ref{assump:path}, then they asymptotically have the same local statistics around a criticality (recall that in Assumption~\ref{assump:path} we assume that along the entire path there is always a criticality at zero). Several steps in the proof of Theorem~\ref{theo:mainthmflow} are very similar to the proof of \cite[Theorem 2.7]{Campbell-Cipolloni-Erdos-Ji2024}, for these reason we omit some details and explain the main differences. Furthermore, without loss of generality throughout the proof we assume that $X$ is a complex Ginibre matrix. In fact, if this is not the case, i.e. $X$ is an IID matrix, then  by a simple Green's function comparison argument (GFT) with two moment matching, we can show that the local statistics of $A+X$ and $A+X^{\mathrm{Gin}}(\C)$ are asymptotically the same. This kind of GFT was performed
 in \cite[Section 7]{Campbell-Cipolloni-Erdos-Ji2024} (see also \cite{Cipolloni-Erdos-Schroder2021}) and it applies to the current case verbatim. 
 In fact, this ``two moment matching GFT" only relies on the single resolvent local law \eqref{eq:locallawgoal}, which is insensitive to whether $z$ is close to a sharp edge or to a criticality.

Let $F_1,\dots, F_k\in C_c^\infty(\C)$. For $j\in [k]$ and $l=s,t$, define
\begin{equation}
\widetilde{F}_j^{(l)}(\cdot ):=NF_j\big(N^{1/4}\gamma_l\cdot\big).
\end{equation}
Proceeding analogously to \cite[Eqs. (4.2)--(4.6)]{Campbell-Cipolloni-Erdos-Ji2024}, to conclude Theorem~\ref{theo:mainthmflow} it is enough to prove\footnote{The implicit constant in $O(\cdot)$ depends on $k$ and the $F_j$ through $\lVert F_j \rVert_{L^\infty}$, $\lVert \Delta F_j\rVert_{L^1}$, and $\mathrm{diam}(\supp F_j)$.}
\begin{equation}
\label{eq:main}
\E\prod_{j=1}^{k}\left(\int_{\C}\wt{F}_{j}^{(s)}(w)\dd (\rho_{\mathcal{A}_s+X}-\rho_{\mathcal{A}_s+\bsx})(w)\right)-\E\prod_{j=1}^{k}\left(\int_{\C}\wt{F}_{j}^{(t)}(w)\dd (\rho_{\mathcal{A}_t+X}-\rho_{\mathcal{A}_t+\bsx})(w)\right) =O(N^{-c}),
\end{equation}
for some small constant $c>0$. Next, using Girko's formula, for both $l=s,t$, we write
\begin{equation}
\begin{split}
\int_{\C}\wt{F}_{j}^{(l)}(w)\dd\rho_{A_l+X}(w)=&N\int_{\C}F_{j}(w)\dd\rho_{N^{1/4}\gamma_l(A_l+X)}(w)	\\
=&-\frac{N}{2\pi}\int_{\C}\Delta_{w}F_{j}(w)\brkt{\log\absv{N^{1/4}\gamma_l(A_l+X)-w}}\dd^{2}w \\
=&-\frac{N}{2\pi}\int_{\C}\Delta_{w}F_{j}(w)\brkt{\log \absv{A_l+X-\gamma_l^{-1}N^{-1/4}w}}\dd^{2}w.
\end{split}
\end{equation}
We point out that in the third line we subtracted $\log(N^{1/4}|\gamma_l|)$ from the integrand using that $\int \Delta F=0$ for compactly supported $F$. We now split
\begin{equation}
\int_{\C}\wt{F}_{j}^{(l)}(w)\dd(\rho_{A_l+X}-\rho_{A_l+\bsx})(w)=I^{(j,l)}_{1}+I^{(j,l)}_{2},
\end{equation}
where for $l=s,t$ we defined
\begin{equation}
\begin{split}
I^{(j,l)}_{1}\deq&-\frac{1}{2\pi}\int_{\C}\Delta_{w}F_{j}(w) L_l(w)\dd^{2}w,\\
I^{(j,l)}_{2}\deq&-\frac{1}{2\pi} \int_{\C}\Delta_{w}F_{j}(w) \int_{0}^{\eta_l}N\Brkt{\frac{\eta}{\absv{A_l+X-\gamma_l^{-1}w}^{2}+\eta^{2}}-\im M_{A_l-\gamma_l^{-1} N^{-1/4}  w}(\ii\eta)}\dd\eta\dd^{2}w,
\end{split}
\end{equation}
where for any $z\in\C$ the matrix $M_{A-z}$ is defined as the solution of \eqref{eq:MDE} for $\widehat{z}=z$.

We will now see that the small $\eta$--regime $I^{(j,l)}_{2}$ contributes with a negligible error and so that the main contribution comes from $I^{(j,l)}_{1}$. The proof of this lemma is postponed to the end of this section.
\begin{lem}
\label{lem:smalletaerr}
Let $X$ be a complex Ginibre matrix. Then, there exists a small $c>0$ such that, for any $j\in [k]$ and $l=s,t$, we have
\begin{equation}
\label{eq:smalletaerr}
\E\big|I^{(j,l)}_{2}\big|\le N^{-c}.
\end{equation}
\end{lem}
Using \eqref{eq:smalletaerr} and that $L_l(w)\prec 1$ uniformly in $w$ by the local law \eqref{eq:locallawgoal} and $\lVert \Delta F_j\rVert\lesssim 1$, we readily obtain
\begin{equation}
\label{eq:almthere}
\E\prod_{j=1}^{k}\left(\int_{\C}\wt{F}_{j}^{(l)}(w)\dd (\rho_{A_l+X}-\rho_{A_l+\bsx})(w)\right)=\E\prod_{j=1}^{k} I_{1}^{(j,l)}+O(N^{-c}).
\end{equation}
Finally, by \eqref{eq:almthere} together with \eqref{eq:Lbound} (see \cite[Proof of Lemma 4.3]{Campbell-Cipolloni-Erdos-Ji2024} for the detailed proof), we immediately conclude \eqref{eq:main}, where again used $\lVert \Delta F_j\rVert\lesssim 1$.

\end{proof}

\begin{proof}[Proof of Lemma~\ref{lem:smalletaerr}]

For $A=A_l$, $\eta=\eta_l$, and $X$ being a complex Ginibre matrix, uniformly over $\absv{w}=O(1)$, we claim
\begin{equation}
\label{eq:smallestsvb}
\mathbb{P}\big[\lambda_1(A+X-N^{-1/4}w)\le\eta\big]\lesssim N^{-c},
\end{equation}
where $\lambda_1(\cdot)$ denotes the smallest singular value\footnote{The smallest singular values $\lambda_1(\cdot)$ should not be confused with the $\lambda_1$ introduced in \eqref{eq:defangle}, denoting the largest eigenvalues if the Hessian matrix $\mathcal{H}$.}. Given \eqref{eq:smallestsvb} and $\absv{\gamma_{l}}\sim 1$, then \eqref{eq:smalletaerr} follows immediately as in \cite[Lemma 4]{Cipolloni-Erdos-Schroder2021}.

We now turn to the proof of \eqref{eq:smallestsvb}. First, we notice that by Definition~\ref{defn:A} (i)--(iii) we have
\begin{equation}
\label{eq:implfneed}
\Brkt{\frac{1}{|A-N^{-1/4}w|^2}}=1+O\big(N^{-1/2}|w|^2\big).
\end{equation}
Define
\begin{equation}
f(E):=\Brkt{\frac{1}{\big||A-N^{-1/4}w|+N^{-1/2}E\big|^2}}.
\end{equation}
By \eqref{eq:implfneed}, we can thus find an $E\in\R, E=O(1)$ such that
\begin{equation}
\label{eq:usrelE}
f(E)=1.
\end{equation}
Here we used that
\[
\partial_E f|_{E=0}=-2N^{-1/2}\brkt{\absv{A-N^{-1/4}w}^{-3}}\sim N^{-1/2},
\]
and so that we can apply the implicit function theorem to show the existence of $E$ so that \eqref{eq:usrelE} is satisfied. In fact, let $\widetilde{A}=\widetilde{A}(E,w):=|A-N^{-1/4}w|+N^{-1/2}E$, then from $\wt{A}> 0$ and $\lVert \wt{A}^{-1}\rVert=O(1)$, we can easily check that $(\wt{A},0)$ satisfies \cite[Assumption 2.6]{Campbell-Cipolloni-Erdos-Ji2024} (written for $(A,z_0)$ therein) and in particular that $\wt{A}$ has a sharp edge at the origin. 
Next, by the unitary invariance of $X$, i.e. that $X$ and $XU$ have the same distribution for deterministic unitary $U\in\C^{N\times N}$, we have
\begin{equation}
\big|A+X-N^{-1/4}w\big|\stackrel{\mathrm{d}}{=}\big|X+|A-N^{-1/4}w|\big|=\big|\widetilde{A}+X-N^{-1/2}E\big|.
\end{equation}
For $\wt{A}+X-N^{-1/2}E$ we can thus follow the proof of the smallest singular value bound from \cite[Eq. (4.20)]{Campbell-Cipolloni-Erdos-Ji2024} to get the desired result (recall that $X$ is a complex Ginibre matrix)
\begin{equation}
\begin{split}
\mathbb{P}\big[\lambda_1(A+X-N^{-1/4}w)\le\eta\big]&=\mathbb{P}\big[\lambda_1(\widetilde{A}+X-N^{-1/2}E)\le\eta\big] \\
&\le 2N\eta\E\Brkt{\frac{\eta}{|\widetilde{A}+X-N^{-1/2}E|^2+\eta^2}} \\
&= 2N\eta_\infty\E\Brkt{\frac{\eta_\infty}{|X-1-N^{-1/2}\tilde{\gamma}E|^2+\eta_\infty^2}}+O(N^{-c}) \\
&\le N^{-c},
\end{split}
\end{equation}
where we used the short--hand notation $\tilde{\gamma}:=\big(\brkt{|A|^{-4}A^*}\brkt{|A|^{-4}}^{-1/2}\big)$. We point out that to go from the second to the third line we used \cite[Eq. (4.12)]{Campbell-Cipolloni-Erdos-Ji2024} for $A=\widetilde{A}$, $z_0=0$, $w=E$, $z_1=1$, and $\eta_1=\eta_\infty$, using the notation therein. Finally, in the last inequality we used \cite[Eq. (28)]{Cipolloni-Erdos-Schroder2021}.
\end{proof}

\section{Evolution of the $\log$--determinant}
\label{sec:techest}

In this section we study the evolution of products of $\log$--determinants as well as of their expectation, i.e. we present the proofs of Lemma~\ref{lem:determ} and Proposition~\ref{prop:univ_Gin}.

\begin{proof}[Proof of Lemma~\ref{lem:determ}]
	We closely follow that of \cite[Lemma 4.8]{Campbell-Cipolloni-Erdos-Ji2024}, hence we only point out the main steps.
	First, the proof of \cite[Lemma 4.8]{Campbell-Cipolloni-Erdos-Ji2024} applies verbatim until Eq. (4.57) therein;
	\begin{equation}
		\begin{split}\label{eq:determ_0}
			&\frac{\dd}{\dd t}\brkt{\brkt{\log\absv{\bsh_{t}^{w}-\ii\eta_{t}}}_{\caM}}
			=\Brkt{\frac{\dd (\E H_{t}^{w}-\ii\eta_{t})}{\dd t}M_{t}^{w}}	\\
			=&-\frac{\dd \brkt{\log\absv{\mr{M}_{t}^{w}}}}{\dd t}
			+\frac{(v_{t}^{w})^{2}}{2}\frac{\dd}{\dd t}\brkt{(\mr{M}_{t}^{w})^{2}}	
			-\frac{(v_{t}^{w})^{4}}{4}\frac{\dd}{\dd t}\brkt{(\mr{M}_{t}^{w})^{4}}
			+v_{t}^{w}\frac{\dd\eta_{t}}{\dd t}\Brkt{(\mr{M}_{t}^{w})^{2}}+O(N^{-3/2}\log N),
		\end{split}
	\end{equation}
	where $v_{t}^{w}=\im\brkt{ M_{t}^{w}(\ii\eta_{t})}+\eta_{t}$.
	Then we expand each term of \eqref{eq:determ_0} except the first around $w=0$. More precisely, for $w=x+\ii y$, we have
	\begin{equation}\label{eq:determ_1}\begin{aligned}
		\brkt{(\mr{M}_{t}^{w})^{2}}=&1+N^{-1/2}\frac{1}{2}\begin{pmatrix}
			\re [\gamma_{t}^{-1}w] \\ \im[\gamma_{t}^{-1}w]
		\end{pmatrix}\tp\caH \begin{pmatrix}
		\re [\gamma_{t}^{-1}w] \\ \im[\gamma_{t}^{-1}w]
		\end{pmatrix}+O(N^{-3/4})	\\
		=&1+N^{-1/2}\frac{I_{4}(t)^{1/2}}{2}\frac{x^{2}+\alpha y^{2}}{1+\alpha}+O(N^{-3/4}),	\\
		\brkt{(\mr{M}_{t}^{w})^{4}}=&I_{4}(t)+O(N^{-1/4}).
	\end{aligned}\end{equation}
	Likewise we can prove the same estimates 
	 for the time derivatives of both sides, except that the errors have an extra factor of $\log N$ from Assumption \ref{assump:path} (iii). Indeed, in the Taylor expansions for the time derivatives, the coefficient of the error is bounded by $\norm{\frac{\dd \caA_{t}}{\dd t}}\norm{(\E H_{t}^{w})^{-1}}^{k}$ for a suitable $k$, say $k=10$.
	
	We recall $\eta_\infty$ from \eqref{eq:def_Herm} and point out that in the following we use the short--hand notation $\alpha_t:=\alpha(\mathcal{A}_t)$. Then, plugging \eqref{eq:determ_1} and its derivative into \eqref{eq:determ_0}, we obtain 
	\begin{equation}\begin{aligned}\label{eq:determ_2}
		&\frac{\dd}{\dd t}\brkt{\brkt{\log\absv{\bsh_{t}^{w}-\ii\eta_{t}}}_{\caM}} \\
		&\qquad=-\frac{\dd}{\dd t}\brkt{\log\absv{\mr{M}_{t}^{w}}}+\frac{(v_t^w)^2I_4(t)^{1/2}}{4 N^{1/2}}\cdot\frac{\dd \alpha_t}{\dd t}\cdot\frac{y^2(1+\alpha_t)-(x^2+\alpha_t y^2)}{(1+\alpha_t)^2}	\\
		&\qquad\quad-\frac{v_{t}^{w}}{4}I_{4}(t)^{-3/4}\frac{\dd I_{4}(t)}{\dd t}\left[(I_{4}(t)^{1/4}v_{t}^{w})^{3}-
		(I_{4}(t)^{1/4}v_{t}^{w})N^{-1/2}\frac{x^{2}+\alpha_t y^{2}}{1+\alpha_t}
		-\eta_{\infty}\right]+O(N^{-5/4}\log N) \\
		&\qquad=-\frac{\dd}{\dd t}\brkt{\log\absv{\mr{M}_{t}^{w}}}-\frac{v_{t}^{w}}{4}I_{4}(t)^{-3/4}\frac{\dd I_{4}(t)}{\dd t}\left[(I_{4}(t)^{1/4}v_{t}^{w})^{3}-
		(I_{4}(t)^{1/4}v_{t}^{w})N^{-1/2}\frac{x^{2}+\alpha_t y^{2}}{1+\alpha_t}
		-\eta_{\infty}\right]+O(N^{-1-\mathfrak{c}_1}).
	\end{aligned}\end{equation}
	Since the big square bracket 
	is exactly the cubic equation for $v_{t}^{w}$, by \eqref{eq:cubic_nor}, we obtain
	\begin{equation}
		\frac{\dd}{\dd t}\brkt{\brkt{\log\absv{\bsh_{t}^{w}-\ii\eta_{t}}}_{\caM}}
		=-\frac{\dd}{\dd t}\brkt{\log\absv{\mr{M}_{t}^{w}}}+O(N^{-1-\mathfrak{c}_1}),
	\end{equation}
	where we used Assumption~\ref{assump:path} (ii) and that $v_t^w=O(N^{-1/4}$).

	The harmonicity of $ w\mapsto  \brkt{\log \absv{\mr{M}_{t}^{w}}}$ follows from the fact that $\norm{(\caA_{t}^{0})^{-1}}\leq \frC_{1}$, due to Assumption \ref{assump:path} (i). We refer to \cite[Eq. (4.46)]{Campbell-Cipolloni-Erdos-Ji2024} for more details. This concludes the proof of Lemma \ref{lem:determ}.
\end{proof}

We now compute the evolution of $\bsL_{t}(\bsw)$. In the following, for a subset $S\in \{1,\dots, k\}$ we denote (cf. \eqref{eq:defL})
\begin{equation}
\bsL_t^{(S)}:=\prod_{j=1, \atop j\notin S}^k L_t(w_j),
\end{equation}
with $L_t$ from \eqref{eq:def_Lt}; a similar notation is used for $\bsw^{(S)}$ as well. We start with (from now on we often use the short-hand notation $G_t^w=G_t^w(\ii\eta_t)$)
\beq
\begin{split}
\frac{\dd L_{t}(w)}{\dd t}
=&\Tr\frac{\dd (H_{t}^{w}-\ii\eta_{t})}{\dd t} G_{t}^{w}-2N\frac{\dd\brkt{\brkt{\log\absv{\bsh_{t}^{w}-\ii\eta_{t}}}_{\caM}}}{\dd t}	\\
=&\Tr\frac{\dd (\E H_{t}^{w_{j}}-\ii\eta_{t})}{\dd t} G_{t}^{w}-2N\frac{\dd\brkt{\log\absv{\mr{M}_{t}^{w}}}}{\dd t}+O(N^{-\mathfrak{c}_1}),
\end{split}
\eeq
where in the last equality we used Lemma~\ref{lem:determ}. We thus obtain
\begin{align}
\label{eq:dL=dHG}
\E\frac{\dd\bsL_{t}(\bsw)}{\dd t}=\sum_{j=1}^{k}\E\Tr\left[\frac{\dd(\E H_{t}^{w_{j}}-\ii\eta_{t})}{\dd t}G^{w_{j}}_{t}(\ii\eta_{t})\right]\bsL_{t}^{(j)}(\bsw^{(j)})+h_{t}(\bsw)+O(N^{-\mathfrak{c}_1+\xi}),
\end{align}
for any arbitrary small $\xi>0$ an explicit $h_{t}(\bsw)$ such that $\Delta_{w_1}\dots\Delta_{w_k} h_{t}(\bsw)=0$. The harmonicity of $h_{t}(\bsw)$ in fact follows from the harmonicity of $\brkt{\log\absv{\mr{M}_{t}^{w}}}$. We point out that in \eqref{eq:dL=dHG} we used that $\absv{L_t(w_j)}\lesssim N^\xi$ with very high probability for any small $\xi>0$, by the local law \eqref{eq:locallawgoal}.

Then Proposition \ref{prop:univ_Gin} will readily follow by the following proposition:
\begin{prop}
\label{prop:dt}
Let $\bbF=\bbC$, $k\in\N$, $C>0$, and $\epsilon>0$ be fixed. Then the following holds uniformly over $\caA_{t}$ satisfying Assumption \ref{assump:path}, $t\in[0,1]$, and $\absv{w},\norm{\bsw}\leq C$: 
\begin{align}
\sum_{j=1}^{k} \E\left[\Brkt{\frac{\dd(\E H_{t}^{w_{j}}-\ii\eta_{t})}{\dd t}G^{w_{j}}_{t}(\ii\eta_{t})}\bsL_{t}^{(j)}(\bsw^{(j)})\right]&=h_{t}(\bsw)+O(N^{\epsilon}\Psi^{5}+N^{-1/2-\mathfrak{c}_1+\epsilon}\Psi^2),	\label{eq:comp_logd_d}\\
\frac{\dd}{\dd t}\E \frac{\brkt{G_{t}^{w}}}{c_{t}}&=O(N^{1+\epsilon}\Psi^{6}+\Psi^3N^{1/2-\mathfrak{c}_1+\epsilon}), \label{eq:comp_sing_d}
\end{align}
where $\Delta_{w_{1}}\cdots\Delta_{w_{k}}h_{t}\equiv 0$ and $\Psi:=1/(N\eta_\infty)=N^{-1/4+\delta}$. 
\end{prop}

We now first conclude the proof of Proposition \ref{prop:univ_Gin} and then present the proof of Proposition~\ref{prop:dt}.

\begin{proof}[Proof of Proposition \ref{prop:univ_Gin}]

Combining \eqref{eq:dL=dHG} with \eqref{eq:comp_logd_d}, and choosing $\delta,\epsilon$ small in term of $\mathfrak{c}_1$, we immediately obtain \eqref{eq:Lbound}.

\end{proof}

\begin{proof}[Proof of Proposition~\ref{prop:dt}]

We first focus on the proof of \eqref{eq:comp_logd_d}, and then explain the very minor changes to obtain \eqref{eq:comp_sing_d}, given \eqref{eq:comp_logd_d}.

To compute the expectation in the lhs. of \eqref{eq:comp_logd_d} we write the resolvent $G_t^{w_j}$ as (recall $\mathring{M}_t^{w_j}:=(\E H_t^{w_j})^{-1}$)
\begin{equation}
\label{eq:eqG}
G_t^{w_j}=\mathring{M}_t^{w_j}+\mathring{M}_t^{w_j}(\ii\eta_t-W)G_t^{w_j},
\end{equation}
and then use Stein's Lemma (Guassian integration by parts) to obtain an expression involving only traces $\Brkt{G_t^{w_j}}$. More precisely, we repeatedly use the equation \eqref{eq:eqG} and Stein's Lemma to factorize any trace containing products of $G$'s and deterministic matrices. These quite technical and non--trivial computations are presented in full detail in \cite{Campbell-Cipolloni-Erdos-Ji2024}, since in the current case the analogous computations would be completely identical they are omitted. In fact, defining
\[
B_{d,t}^{w_j}:=\frac{\dd \E H_t^{w_j}}{\dd t} \mathring{M}^{w_j},
\]
and proceeding exactly as in \cite[Eqs. (5.13)--(5.31)]{Campbell-Cipolloni-Erdos-Ji2024}, we obtain
\beq\begin{split}
\label{eq:cancel_0}
&\E\Brkt{\frac{\dd (\E H_{t}^{w_{j}}-\ii\eta_{t})}{\dd t}G_t^{w_{j}}}\bsL_t^{(j)}	\\
=&\E\Brkt{2\frac{\dd (\E H_{t}^{w_{j}}-\ii\eta_{t})}{\dd t}G_t^{w_{j}}E_{1}}\bsL_t^{(j)}	\\
=&\Brkt{2B_{d,t}^{w_{j}}E_{1}}\E\bsL_t^{(j)} \\
&+\left(\brkt{2B_{d,t}^{w_{j}}(\mr{M}_t^{w_{j}})^{2}E_{1}}-\frac{\brkt{2B_{d,t}^{w_{j}}(\mr{M}_t^{0})^{4}E_{1}}}{I_{4}}\left(\brkt{(\mr{M}_t^{w_{j}})^{2}}-1\right)\right)\E\brkt{G_t^{w_{j}}}^{2}\bsL_t^{(j)}	\\
&+\left(\brkt{4B_{d,t}^{w_{j}}(\mr{M}_t^{0})^{2}E_{1}}-\frac{\brkt{2B_{d,t}^{w_{j}}(\mr{M}_t^{0})^{4}E_{1}}}{I_{4}(t)}+\frac{\dd c_{t}}{\dd t}\frac{1}{c_{t}}\right)\ii\eta\E\brkt{G_t^{w_{j}}}\bsL_t^{(j)}	\\
&+\sum_{\ell}^{(j)}\left(\brkt{2B_{d,t}^{w_{j}}\mr{M}_t^{w_{\ell}}\mr{M}_t^{w_{j}}E_{1}}-\frac{\brkt{2B_{d,t}^{w_{j}}(\mr{M}_t^{0})^{4}E_{1}}}{I_{4}(t)}(\brkt{2\mr{M}_t^{w_{\ell}}\mr{M}_t^{w_{j}}E_{1}}-1)\right)\E\frac{\brkt{G_t^{w_{\ell}}E_{2}G_t^{w_{j}}E_{1}}}{N}\bsL_t^{(j,\ell)}	\\
&+\sum_{\ell}^{(j)}\frac{1}{2N}\left(\brkt{2B_{d,t}^{w_{j}}(\mr{M}_t^{w_{j}})^{2}E_{1}}-\frac{\brkt{2B_{d,t}^{w_{j}}(\mr{M}_t^{0})^{4}E_{1}}}{I_{4}(t)}\right)\E\bsL_t^{(j,\ell)}
+O(N^{\epsilon}\Psi^{5}),
\end{split}
\eeq
where $\sum^{(j)}$ denotes that the index $j$ is missing from the summation, and the matrices $E_1,E_2\in\C^{2N\times 2N}$ are defined as (here $I_N$ denotes the $N\times N$ identity matrix)
\begin{equation}
E_1:=\left(\begin{matrix}
I_N & 0 \\
0 & 0
\end{matrix}\right), \qquad\quad E_2:=\left(\begin{matrix}
0 & 0 \\
0 & I_N
\end{matrix}\right).
\end{equation}
We now show that the real part of the rhs. of \eqref{eq:cancel_0} consists of a harmonic function plus a negligible (of size $O(\Psi^5)$) error. Several computations are similar to \cite[Eqs. (5.32)--(5.42)]{Campbell-Cipolloni-Erdos-Ji2024} and so omitted; here we focus on highlighting the main differences. For the term in the first line we write
\begin{equation}
\re\brkt{2B_{d,t}^{w_{j}}E_{1}}=\Brkt{\frac{\dd \E H_{t}^{w_{j}}}{\dd t}\mr{M}_{t}^{w_{j}}}=-\frac{\dd}{\dd t}\brkt{\log\absv{\mr{M}_{t}^{w_{j}}}},
\end{equation}
which is harmonic in $w_j$ by Lemma~\ref{lem:determ}. Notice that also the term in the last line of \eqref{eq:cancel_0} is harmonic since it does not depend on $w_{\ell}$; we thus ignore this term in the following.  For the term in the fourth line, we use that (recall that $\alpha_t=\alpha(\mathcal{A}_t)$)
\begin{equation}
\label{eq:usefexp}
\begin{split}
\brkt{2B_{d,t}^{w_{j}}(\mr{M}_t^{w_{j}})^{2}E_{1}}&=-\frac{1}{2}\frac{\dd}{\dd t} \Brkt{(\mathring{M}_t^{w_j})^2} \\
&=-\frac{1}{8}N^{-1/2}I_4(t)^{-1/2}\frac{\dd I_4(t)}{\dd t}\cdot\frac{x_j^2+\alpha_t y_j^2}{1+\alpha_t} \\
&\quad-\frac{I_4(t)^{1/2}}{4N^{1/2}}\cdot \frac{\dd \alpha_t}{\dd t}\cdot\frac{y_j^2(1+\alpha_t)-(x_j^2+\alpha_t y_j^2)}{(1+\alpha_t)^2}+O(N^{-3/4}\log N) \\
&=-\frac{1}{8}N^{-1/2}I_4(t)^{-1/2}\frac{\dd I_4(t)}{\dd t}\cdot\frac{x_j^2+\alpha_t y_j^2}{1+\alpha_t}+O(N^{-1/2-\mathfrak{c}_1}) \\
\brkt{2B_{d,t}^{w_{j}}(\mr{M}_t^{0})^{4}E_{1}}&=\Brkt{\frac{\dd \E H_t^{w_j}}{\dd t} \mathring{M}_t^{w_j} (\mathring{M}^0)_t^4 E_1}=-\frac{1}{4}\frac{\dd I_4(t)}{\dd t}+O(N^{-1/4}\log N),
\end{split}
\end{equation}
where we used \eqref{eq:determ_1}. In fact, this gives
\begin{equation}
\brkt{2B_{d,t}^{w_{j}}(\mr{M}_t^{w_{j}})^{2}E_{1}}-\frac{\brkt{2B_{d,t}^{w_{j}}(\mr{M}_t^{0})^{4}E_{1}}}{I_{4}}\left(\brkt{(\mr{M}_t^{w_{j}})^{2}}-1\right)=O(\Psi^3+N^{-1/2-\mathfrak{c}_1}).
\end{equation}
This concludes the estimate for the fourth line. For the term in the fifth line, recalling $c_t=I_4(t)^{-1/4}$ and using \eqref{eq:usefexp}, we obtain
\begin{equation}
\brkt{4B_{d,t}^{w_{j}}(\mr{M}^{0})^{2}E_{1}}-\frac{\brkt{2B_{d,t}^{w_{j}}(\mr{M}^{0})^{4}E_{1}}}{I_{4}}+\frac{\dd c_{t}}{\dd t}\frac{1}{c_{t}}=\frac{1}{4}I_4(t)^{-1}\frac{\dd I_4(t)}{\dd t}-\frac{1}{4}I_4(t)^{-1}\frac{\dd I_4(t)}{\dd t}+O(N^{-1/4}\log N)=O(\Psi).
\end{equation}

We are now only left with the term in the penultimate line of \eqref{eq:cancel_0}. First, we notice that, by \eqref{eq:determ_1} and simple Taylor expansion, we have
\begin{equation}
\label{eq:esterr}
\lVert B_{d,t}^{w_j}-B_{d,t}^0\lVert=O(N^{-1/4}), \qquad\quad \brkt{2\mr{M}_t^{w_{j}}\mr{M}_t^{w_{\ell}}E_{1}}-1=O(N^{-1/2}).
\end{equation}
Then, using (here we use that $\mr{M}_t^{w_{j}}E_{1}=E_2\mr{M}_t^{w_{j}}$) 
\begin{equation}
\brkt{2B_{d,t}^{w_{j}}\mr{M}_t^{w_{\ell}}\mr{M}_t^{w_{j}}E_{1}}=\Brkt{2\mr{M}_t^{w_{j}}\frac{\dd \E H_{t}^{w_{j}}}{\dd t}\mr{M}_t^{w_{j}}\mr{M}_t^{w_{\ell}}E_{2}}=-\Brkt{2\frac{\dd \mr{M}_t^{w_{j}}}{\dd t}\mr{M}_t^{w_{\ell}}E_{2}},
\end{equation}
together with \eqref{eq:esterr}, we can write the penultimate line of \eqref{eq:cancel_0} as
\begin{equation}
\label{eq:cancel_2pt_0}
2\re\sum_{j>\ell}\left(-\Brkt{\frac{\dd (\mr{M}_t^{w_{j}}\mr{M}_t^{w_{\ell}})}{\dd t}E_{2}}-\frac{\brkt{B_{d,t}^{0}(\mr{M}_t^{0})^{4}}}{I_{4}}\left(\brkt{2\mr{M}_t^{w_{j}}\mr{M}_t^{w_{\ell}}E_{2}}-1\right)\right)\frac{\E\brkt{G_t^{w_{\ell}}E_{2}G_t^{w_{j}}E_{1}}}{N}\bsL^{(j,\ell)}+O(\Psi^5),
\end{equation}
where we used that $N^{-1}\brkt{G_t^{w_{\ell}}E_{2}G_t^{w_{j}}E_{1}}=O(\Psi^2)$. We point out that the bound $N^{-1}\brkt{G^{w_{\ell}}E_{2}G^{w_{j}}E_{1}}=O(\Psi^2)$ follows from a Schwarz inequality $|\brkt{G_t^{w_{\ell}}E_{2}G_t^{w_{j}}E_{1}}|\le \sqrt{\brkt{\Im G_t^{w_{\ell}}}\brkt{\Im G_t^{w_{\ell}}}}/\eta$ together with $\brkt{\Im G_t^w}=O(N^{-1/4})$ by Lemma~\ref{lem:asm} and Theorem~\ref{theo:llaw}. Next, we notice
\begin{equation}
\begin{split}
&-\Brkt{\frac{\dd (\mr{M}_t^{w_{j}}\mr{M}_t^{w_{\ell}})}{\dd t}E_{2}}-\frac{\brkt{B_{d,t}^{0}(\mr{M}_t^{0})^{4}}}{I_{4}}\left(\brkt{2\mr{M}_t^{w_{j}}\mr{M}_t^{w_{\ell}}E_{2}}-1\right) \\
&\qquad\quad= -\frac{1}{2}I_4(t)^{1/2}\frac{\dd}{\dd t}\left[I_4(t)^{-1/2}\left(\brkt{2\mr{M}_t^{w_{j}}\mr{M}_t^{w_{\ell}}E_{2}}-1\right)\right] \\
&\qquad\quad=\frac{1}{2}I_4(t)^{1/2}\frac{\dd}{\dd t}\left[I_4(t)^{-1/2}\left(I_4(t)^{1/2}N^{-1/2}\cdot\frac{(x_\ell+x_j)^2+\alpha_t (y_\ell+y_j)^2}{2(1+\alpha_t)}+O(N^{-3/4})\right)\right] \\
&\qquad\quad=O(N^{-3/4}\log N)=O(\Psi^3+N^{-1/2-\mathfrak{c}_1}).
\end{split}
\end{equation}
This, together with $N^{-1}\brkt{G_t^{w_{\ell}}E_{2}G_t^{w_{j}}E_{1}}=O(\Psi^2)$, concludes the proof of \eqref{eq:comp_logd_d}.

Given \eqref{eq:comp_logd_d}, we readily obtain \eqref{eq:comp_sing_d} as well. In fact, \eqref{eq:comp_sing_d} follows by simply differentiating \eqref{eq:comp_logd_d}, for $k=1$, with respect to $\ii \eta_t$. Note that by differentiating \eqref{eq:comp_logd_d} the error term deteriorates by a factor $N^{3/4+\delta}$. We point out that by taking the derivative of the error term in \eqref{eq:comp_logd_d} we mean that inspecting the proof of \eqref{eq:comp_logd_d} one can see that if each term is differentiated then the error term deteriorates by a factor $1/\eta_t=N^{3/4+\delta}$; this gives the desired error term in \eqref{eq:comp_sing_d}. We omit the details since they are already given in \cite[Proof of Proposition 5.2, (5.5)]{Campbell-Cipolloni-Erdos-Ji2024}; nothing changes in the current case. 
\end{proof}

\section{Path construction} \label{sec:path}

In this section, we construct a  $C^1$-path  $\caA_{t}$, $t\in [0,1]$, 
 satisfying Assumption~\ref{assump:path}, starting from a given normal matrix
$\caA_0:=A$ and ending at another matrix $\caA_{1}$ with additional properties. More precisely, the final point $\caA_{1}$ must have (i) $O(1)$ number of distinct eigenvalues\footnote{Hereafter we refer to this property as \emph{having a finite spectrum}, where ``finite" really means $N$-independent.} and (ii) an $N$-independent spectral distribution up to $O(N^{-1})$ errors. 
In fact, our final path will be the concatenation of two paths that establish properties (i) and (ii), respectively. The first path is constructed in a general setting in Proposition \ref{prop:path_finite} below. For the second path, we first prove in Proposition \ref{prop:path_fix} below that two matrices with finite spectra can be connected  by  a path if their spectral distributions are close, i.e., $o(1)$ away in total variation distance.
Then along the proof of Theorem \ref{theo:mainthm}, in Section \ref{sec:pf_main}, we use Proposition \ref{prop:path_fix} together with a compactness argument to construct the second path.

Notice that all conditions involved in Assumption \ref{assump:path} are invariant under unitary conjugations $\caA_{t}\mapsto U\caA_{t}U\adj$ and rotations $\caA_{t}\mapsto \e{\ii\varphi}\caA_{t}$, $\varphi\in \R$.
Note also that the normalization $\brkt{\absv{\caA_{t}}^{-2}}=1$ can be easily achieved by taking $\brkt{\absv{\caA_{t}}^{-2}}^{1/2}\caA_{t}$ if $\caA_{t}$ satisfies all the rest of Assumption \ref{assump:path}, at the cost of increasing $\frC_{1}$ to $\frC_{1}^{3}$. 
With these in mind, we make a number of preliminary simplification in the following lemma; in words, we diagonalize $\caA_{t}$, rotate $\caA_{t}$ so that $\brkt{\caA_{t}^{-3}(\caA_{t}\adj)^{-1}}\geq0$, drop the normalization $\brkt{\absv{\caA_{t}}^{-2}}=1$, and
we work with the inverse matrix 
    $\caA_{t}^{-1}$ instead of $\caA_{t}$ to make formulas shorter. The conditions will be expressed in terms of a new flow $\caB_t$ that is a simply transformed version 
 of $\caA_t$,  see \eqref{eq:B_to_A} below.
\begin{lem}\label{lem:path_inv}
	Let $\frC>1,\frc>0$ be constants and let $\caB_{\cdot}:[0,1]\to\C^{N\times N}$ be a matrix-valued path that satisfies the following\footnote{\label{fn:path_cond}We use the same numbering for items (i)--(v) consistently through out this section; (i) for diagonality, (ii) for the (inverse) norm bounds, (iii) for the derivative-norm bound, and (iv) for the criticality. The last condition (v) varies slightly by context (in which case we  write (v$'$), (v$''$), or (v$'''$)), but always concerns the parameter $\chi$.} for all $t\in[0,1]$.
	\begin{itemize}
		\item[(i)] $\caB_{t}$ is a diagonal matrix
		\item[(ii)] $\norm{\caB_{t}}\leq \frC, \norm{\caB_{t}^{-1}}\leq\frC$.
		\item[(iii)]  $\norm{\dd\caB_{t}/\dd t}\leq \frC.$
		\item[(iv)] $\brkt{\caB_{t}^{2}\caB_{t}\adj}=0$.
		\item[(v)] $\chi(\caB_{t})\geq 0$ and $\absv{\dd\chi(\caB_{t})/\dd t}\leq N^{-\frc}$, where for any $B\in \C^{N\times N}$
		 we define
		\begin{equation}
			\chi(B)\deq \frac{\brkt{B^{3}B\adj}}{\brkt{\absv{B^{2}}^{2}}}.
		\end{equation}
	\end{itemize}
	Then, for any unitary matrix $U\in\caU(N)$ and $\varphi\in(0,2\pi]$, the path
	\begin{equation}\label{eq:B_to_A}
		\caA_{t}\deq \e{\ii\varphi}\brkt{\absv{\caB_{t}}^{-2}}^{-1/2}U\caB_{t}^{-1}U\adj
	\end{equation}
	satisfies Assumption \ref{assump:path} with parameters $(\frC_{1},\frc_{1})$ depending only on $(\frC,\frc)$ and $\frc_1>0$ can be made arbitrary small if $\frc>0$ is small.
	Furthermore, the number $\alpha(\caA_{t})$ in Assumption \ref{assump:path} (ii) is given by
	\begin{equation}\label{eq:alpha_chi}
		\alpha(\caA_{t})=\frac{1-2\chi(\caB_{t})}{1+2\chi(\caB_{t})}.
	\end{equation}
\end{lem}
\begin{proof}
	The fact that $\caA_{t}$ has a criticality at the origin is an immediate consequence of the observations above. The estimate for $\norm{\dd\caA_{t}/\dd t}$
	follows from that of $\caB_{t}$ and the assumption that $\norm{\caB_{t}}$ and $\norm{\caB_{t}^{-1}}$ are uniformly bounded over $t\in[0,1]$.
	
	Finally, to see \eqref{eq:alpha_chi} and also Assumption \ref{assump:path} (ii), we only need to notice that if $A$ is a general matrix having a criticality at the origin, the two eigenvalues of the 
	Hessian $\caH$ are given by (see \eqref{eq:Hess_comput} for a proof)
	\begin{equation}
		\lambda_{1}=2\Brkt{\frac{1}{\absv{A^{2}}^{2}}}+4\Absv{\Brkt{\frac{1}{A^{3}A\adj}}},\qquad
		\lambda_{2}=2\Brkt{\frac{1}{\absv{A^{2}}^{2}}}-4\Absv{\Brkt{\frac{1}{A^{3}A\adj}}}.
	\end{equation}
\end{proof}

When Lemma \ref{lem:path_inv} is used in practice (see Section \ref{sec:pf_main}), for a given normal matrix $A$ having a criticality at the origin we consider the diagonal matrix $\caB_{0}=B$ given by\footnote{As in Footnote \ref{fn:alpha=0}, we make the cosmetic choice $\varphi=0$ when $\brkt{A^{-3}(A\adj)^{-1}}=0$.}
\begin{equation}\label{eq:A_to_B}
	B=\exp\left(-\ii\varphi\right)U\adj A^{-1}U,\qquad \varphi=\frac{\arg\brkt{A^{-3}(A\adj)^{-1}}}{2}\in[0,\pi),\qquad U\in\caU(N).
\end{equation}
Then, given a path $\caB_{t}$ starting from $\caB_{0}$ (notice that $\chi(B)=\chi(\caB_{0})\geq0$ due to our choice of $\varphi$), our path of deformations $\caA_{t}$ will be as in \eqref{eq:B_to_A} where $U,\varphi$ are chosen in \eqref{eq:A_to_B} so that $\caA_{0}=A$. Notice that only the eigenvalues of such an $\caA_{t}$ are evolving over $t\in[0,1]$,
 its eigenvectors remain intact. 

In what follows, we construct two flows $\caB_{t}$ that shall be converted into the final 
$\caA_{t}$-flows
by Lemma~\ref{lem:path_inv}.
In Section \ref{sec:path_finite}, Proposition \ref{prop:path_finite}, we construct a flow $\caB_{t}$ that shall give us the first path, that is, such that the final matrix $\caB_{1}$ has a finite spectrum. 
Then in Section \ref{sec:path_fix}, precisely in Proposition \ref{prop:path_fix}, we construct a flow between two given matrices $\caB_{0},\caB_{1}$ with finite spectra and with spectral distributions close to each other, such that the value of $\chi(\caB_{t})$ has a bounded derivative along the flow.

We remark that we assume that $B=\caB_{0}$ has a genuinely complex spectrum, i.e. that $1-\chi(B)\gtrsim 1$. This in turn corresponds to genuinely complex $A$ and $\alpha(A)+1/3\gtrsim 1$ via \eqref{eq:alpha_chi}, and is guaranteed by our assumption \eqref{eq:alpha_assu} on $A$ and $\alpha_\infty\in (-1/3, 1]$.
 This restriction will help in stabilizing the map 
\begin{equation}
	B\mapsto (\brkt{B^{2}B\adj},\chi(B)),
\end{equation}
that is, in showing that any value in a small vicinity of $(\brkt{B^{2}B\adj},\chi(B))$ has a corresponding matrix achieving that value; 
see Lemma \ref{lem:Jac} below.

\subsection{Step 1: From general $A$ to finitely supported $A$}\label{sec:path_finite}
\begin{prop}\label{prop:path_finite} 
	Let $\frC>1$ be a constant and let $B\in\C^{N\times N}$ be a diagonal matrix with 
	\begin{itemize}
		\item[(a)] $\norm{B},\norm{B^{-1}}\leq \frC$,
		\item[(b)] $\brkt{B^{2}B\adj}=0$,
		\item[(c)] $0\leq \chi(B)\leq 1-1/\frC$.
	\end{itemize}
	Then there exist constants $\frC_{1},M>0$, depending only on $\frC$, and a $C^{1}$ map $[0,1]\ni t\mapsto\caB_{t}\in\C^{N\times N}$ such that the following\footnote{Continued from Footnote \ref{fn:path_cond}, we consistently number a condition as (0), (0$'$), or (0$''$) when it concerns the initial and final matrices of a path.} hold:
	\begin{itemize}
		\item[(0)] $\caB_{0}=B$ and $\rho_{\caB_{1}}$, the spectral measure of $\caB_1$,
		 is supported on at most $M$ points.
		\item[(i)--(iv)] Items (i)--(iv) of Lemma \ref{lem:path_inv} hold for all $t\in[0,1]$, with the constant $\frC_{1}$ in place of $\frC$.
		\item[(v$'$)] $\chi(\caB_{t})=\chi(B)$, i.e. $\brkt{\caB_{t}^{3}\caB_{t}\adj}-\chi(B)\brkt{\absv{\caB_{t}}^{4}}=0$ for all $t\in[0,1]$. 
	\end{itemize}
\end{prop}

In the next lemma we show that the two key relations
\begin{equation}\label{eq:key}
	\brkt{B^{2}B\adj}=0,\qquad \brkt{B^{3}B\adj}-\chi(B) \brkt{\absv{B}^{4}}=0
\end{equation}
are stable under a small simultaneous perturbation of two eigenvalues $z_1, z_2$ of $B$  (possibly with different relative multiplicites, expressed by 
the parameter $p$ below) as long as they are
on  different sides of the imaginary axis, by bounding the corresponding Jacobian from below. 
The imaginary axis plays this special role since we  rotated $B$ so that $\chi(B)\geq0$. 

\begin{lem}\label{lem:Jac}
	For $\chi\in[0,1]$ and $p\in[0,1]$, define the map $F_{\chi,p}\equiv (F_{1},F_{2}):\C^{2}\to\C^{2}$ by
	\begin{equation}\begin{aligned}\label{eq:def_F}
			F_{1}(z_{1},z_{2})\deq& p{z_{1}^{2}\ol{z}_{1}}+(1-p){z_{2}^{2}\ol{z}_{2}}, \\
			F_{2}(z_{1},z_{2})\deq& p{z_{1}^{3}\ol{z}_{1}}+(1-p){z_{2}^{3}\ol{z}_{2}}-\chi\left[p{\absv{z_{1}}^{4}}+(1-p){\absv{z_{2}}^{4}}\right].
	\end{aligned}\end{equation}
	For any $c\in(0,1/2)$, there exists $C_{1}\equiv C_{1}(c)>0$ such that if $z_{1},z_{2}$, $\chi$, and $p$ satisfy
	\begin{equation}\label{eq:z1z2}
		c\leq\absv{z_{i}}\leq 1/c, \quad 
		(\re z_{1})(\re z_{2})\leq 0,\quad 
		\absv{\re z_{1}}+\absv{\re z_{2}}\geq c,\quad 
		\chi\leq 1-c,\quad c\leq p\leq 1-c,
	\end{equation}
	then 
	\begin{equation}
		\norm{DF_{\chi,p}(z_{1},z_{2})^{-1}}\leq C_{1}\equiv C_{1}(c)<\infty,
	\end{equation}
	where $DF_{\chi,p}\in\R^{4\times 4}$ is the differential of $F_{\chi,p}:\R^{4}\to\R^{4}$ 
	after canonically identifying $\C$ with $\R^{2}$.
\end{lem}
\begin{proof}
	We first write $F_{i}(z_{1},z_{2})=pH_{i}(z_{1})+(1-p)H_{i}(z_{2})$, where
	\begin{equation}
		H_{1}(z)\deq z^{2}\ol{z},\qquad H_{2}(z)\deq z^{3}\ol{z}-\chi\absv{z}^{4}.
	\end{equation} 
	We aim at showing that if $(z_{1},z_{2},\chi,p)$
	satisfies \eqref{eq:z1z2} then
	\begin{equation}\label{eq:z1z2_Jac}
		\Absv{\det\begin{pmatrix}
				DH_{1}(z_{1}) & DH_{1}(z_{2})	\\
				DH_{2}(z_{1}) & DH_{2}(z_{2})
		\end{pmatrix}}\geq\delta_{1}(c)>0,
	\end{equation}
	for some constant $\delta_{1}(c)>0$. Indeed, then we easily find from $c<p<1-c$  that
	\begin{equation}\label{eq:Jac_p}
		\absv{\det DF_{\chi,p}(z_{1},z_{2})}=p(1-p)\Absv{\det\begin{pmatrix}
				DH_{1}(z_{1}) & DH_{1}(z_{2})	\\
				DH_{2}(z_{1}) & DH_{2}(z_{2})
		\end{pmatrix}}\geq c\delta_{1}(c),
	\end{equation}
	and from the trivial bound $\norm{DF(z_{1},z_{2})}\leq 100 (\absv{z_{1}}\wedge\absv{z_{2}})^{3}\leq 100 c^{-3}$ we conclude
	\beq
	\norm{DF(z_{1},z_{2})^{-1}}\leq \absv{\det DF(z_{1},z_{2})}^{-1}\norm{DF(z_{1},z_{2})}^{3}\leq c_{1}c^{9}\delta_{1}(c)^{-1},
	\eeq
	with a numeric constant $c_{1}>0$.
	
	After an easy (yet lengthy) explicit algebraic computation, we find that the determinant in \eqref{eq:z1z2_Jac} can be computed as
	\begin{equation}\begin{aligned}
			\Absv{\det\begin{pmatrix}
				DH_{1}(z_{1}) & DH_{1}(z_{2})	\\
				DH_{2}(z_{1}) & DH_{2}(z_{2})
		\end{pmatrix}}\sim&\frac{1}{8\absv{z_{1}}^{2}\absv{z_{2}}^{2}}\det\begin{pmatrix}
				DH_{1}(z_{1}) & DH_{1}(z_{2}) \\
				DH_{2}(z_{1}) & DH_{2}(z_{2})
			\end{pmatrix}
			\\
			=&\left[
			-x_{1}x_{2}((y_{1}^{2}(x_{2}^{2}(7-4\chi) +2y_{2}^{2}(2-\chi)))
			+x_{1}^{2}(y_{2}^{2}(7-4\chi) +6x_{2}^{2}(1-\chi)))\right]	\\
			&+\left[3x_{1}^{4}(x_{2}^2+y_{2}^2) (1-\chi)+3x_{2}^{4}(x_{1}^{2}+y_{1}^{2})(1-\chi)\right] \\
			&+\left[3y_{1}^{2}y_{2}^{2}(1+\chi)(y_{1}-y_{2})^{2}
			+2x_{1}^{2}x_{2}^{2}(3y_{1}^{2}-(2+\chi)y_{1}y_{2}+3y_{2}^{2})\right]	\\
			&+y_{1}^{2}x_{2}^{2}\left[
			3y_{1}^{2}(1+\chi) 
			-y_{1}y_{2}(7+4\chi)+6y_{2}^{2}\right]	\\
			&+x_{1}^{2}y_{2}^{2}\left[6y_{1}^{2}
			-y_{1}y_{2}(7+4\chi)+3y_{2}^{2}(1+\chi)\right]	\\
			=:&(I)+ (II) +(III) +y_{1}^{2}x_{2}^{2}(IV) +x_{1}^{2}y_{2}^{2}(V),
	\end{aligned}\end{equation}
	where we write $z_l = x_l +\ii y_{l}$ for $l=1,2$.
	Recalling $\chi\in[0,1-c]$, 
	 we easily find that when $(z_{1},z_{2})$ satisfy \eqref{eq:z1z2} then
	\begin{equation}
		(I)+(III)\geq0,\qquad (II)>\delta_{2}(c)>0,
	\end{equation}
	for a constant $\delta_{2}(c)>0$ depending only on $c$. Thus it only remains to show that $(IV)$ and $(V)$ are nonnegative, in which case \eqref{eq:z1z2_Jac} immediately follows with the choice $\delta_{1}(c)=8c^{4}\delta_{2}(c)$.
	
	Notice that $(IV)\geq 0$ holds for all $y_{1},y_{2}\in\R$ if and only if 
	\begin{equation}
		\frac{7+4\chi}{\sqrt{18(1+\chi)}}\leq 2,
	\end{equation}
	which trivially holds since $\chi\in[0,1]$. The same applies to $(V)$ by interchanging $y_{1}$ and $y_{2}$.
\end{proof}

\begin{rem}
	The condition $\chi\leq 1-c$ is necessary for Lemma \ref{lem:Jac}. Indeed, $\chi=1$ corresponds to Hermitian $A$, and the Jacobian of $F_{\chi=1,p}$ vanishes for $z_{1},z_{2}\in\R$ since $F_{2}$ is identically zero.
\end{rem}

The next lemma forms the basis of our proof of Proposition \ref{prop:path_finite}, asserting that two small (in diameter)  ``clusters" of eigenvalues of $B$ on the left and right half-planes can be shrunken into two point masses while keeping the two key quantities (on the left-hand sides of \eqref{eq:key}) constant. We explain more details 
 how we use the lemma after its proof.
\begin{lem}\label{lem:box}
	Fix $c\in(0,1/2)$ and let $z_{1},z_{2},\chi$ satisfy \eqref{eq:z1z2}. Then there exist positive 
	constants $h=h(c)$ and $\frC_{1}=\frC_{1}(c)$ such that the following holds. 
	For any $m_{1},m_{2}\in\N$ with $m_{1}/m_{2}\in(2c,1/(2c))$ and $V_{i}\in\C^{m_{i}}$ with $\norm{V_{i}-z_{i}\lone_{m_{i}}}_{\infty}\leq h$ for $i=1,2$, there exist a $C^{1}$ map $\caV_{i}= (\caV_{i1}, \caV_{i2}, \ldots )
	:[0,1]\to \C^{m_{i}}$ and a point 
	$\wt{z}_{i}\in \C$ that satisfies\footnote{Notice that item (i) is absent; the diagonality is always guaranteed since the corresponding matrix is $\diag(\caV_{1}(t),\caV_{2}(t))$.}
	\begin{itemize}
		\item[(0$'$)] $\caV_{i}(0)=V_{i}$ and $\caV_{i}(1)=\wt{z}_{i}\lone_{m_{i}}$.
		\item[(ii)] $\caV_{i}(t)\in (D(0,\frC_{1})\setminus D(0,1/\frC_{1}))^{m_{i}}$ for all $t\in[0,1]$ (where $D(0, r): = \{ z\in \C\; : \; |z|\le r\}$).
		\item[(iii)] $\norm{\caV_{i}'(t)}_{\infty}\leq \frC_{1}$ for all $t\in[0,1]$.
		\item[(iv)] For all $t\in[0,1]$, \begin{equation}\label{eq:box_cond_crit}
				\frac{\dd}{\dd t}\sum_{i=1,2}\sum_{j=1}^{m_{i}}{\caV_{ij}(t)^{2}\ol{\caV_{ij}(t)}}=0.
			\end{equation}
		\item[(v$''$)] For all $t\in[0,1]$, \begin{equation}\label{eq:box_cond_chi}
				\frac{\dd}{\dd t}\sum_{i=1,2}\sum_{j=1}^{m_{i}}\left({\caV_{ij}(t)^{3}\ol{\caV}_{ij}(t)}-\chi{\absv{\caV_{ij}(t)}^{4}}\right)=0.
			\end{equation}
	\end{itemize}
\end{lem}
\begin{proof}
	In this proof, we equip the space $\C^{m_{1}+m_{2}}$ with the max norm $\norm{\cdot}_{\infty}$ 
	and $\C^{2}$ with the usual Euclidean norm. A matrix in $\C^{(m_{1}+m_{2})\times 2}\cong \R^{2(m_{1}+m_{2})\times 4}$ (and its dual $\R^{4\times 2(m_{1}+m_{2})}$) is endowed with the operator norm inherited from $\norm{\cdot}_{\infty}$, which we denote by $\norm{\cdot}_{\infty,2}$.
	
	Consider the function $\boldsymbol{\caF}\equiv (\caF_{1},\caF_{2}):\C^{m_{1}}\times\C^{m_{2}}\to \C^{2}$ given by 
	\begin{equation}
		\caF_{1}(U_{1},U_{2})\deq \frac{1}{m_{1}+m_{2}}\sum_{i=1,2}\sum_{j=1}^{m_{i}}{U_{ij}^{2}\ol{U_{ij}}},\qquad 
		\caF_{2}(U_{1},U_{2})\deq \frac{1}{m_{1}+m_{2}}\sum_{i=1,2}\sum_{j=1}^{m_{i}}\left({U_{ij}^{3}\ol{U}_{ij}}-\chi{\absv{U_{ij}}^{4}}\right),
	\end{equation}
	where $U_{i}=(U_{i1},\ldots,U_{im_{i}})$ for $i=1,2$.
	One can easily see that $\norm{D\boldsymbol{\caF}}_{\infty,2}\leq 100\norm{(U_{1},U_{2})}_{\infty}^{3}$. Similarly we have
	\begin{equation}\label{eq:DF_Lip}
		\norm{D\boldsymbol{\caF}(U_{1},U_{2})-D\boldsymbol{\caF}(U_{1}',U_{2}')}_{\infty,2}\leq 100C^{4}\norm{(U_{1},U_{2})-(U_{1}',U_{2}')}_{\infty},
	\end{equation}
	for $\norm{U_{i}}_{\infty},\norm{U_{i}'}_{\infty}\leq C$.
	
	We aim at applying the implicit function theorem to the map
	\begin{equation}
		\wt{F}(U_{1},U_{2},w_{1},w_{2})\deq \boldsymbol{\caF}(U_{1}+w_{1}\lone_{m_{1}},U_{2}+ w_{2}\lone_{m_{2}})\in\C^{2},\quad (U_{1},U_{2},w_{1},w_{2})\in\C^{m_{1}}\times\C^{m_{2}}\times \C\times\C
	\end{equation}
	around the point $(V_{1},V_{2},0,0)$. Computing the partial derivative with respect to the four real dimensional variable $\bsw=(w_{1},w_{2})$, we get
	\begin{equation}\label{eq:DtF}
		D_{\bsw}\wt{F}=D\boldsymbol{\caF} J_{m_{1},m_{2}}\in \R^{4\times 4},\quad 
		J_{m_{1},m_{2}}\deq
		\begin{pmatrix}
			\lone_{m_{1}} &  &&0 \\
			& \lone_{m_{1}} &&\\
			&&\lone_{m_{2}} & \\
			0&&&\lone_{m_{2}}
		\end{pmatrix}\in\R^{2(m_{1}+m_{2})\times 4}
	\end{equation}
	where $D\boldsymbol{\caF}$ is realized in $\R^{4\times 2(m_{1}+m_{2})}$. The other partial derivative $D_{\bsU}\wt{F}$ with respect to $\bsU=(U_{1},U_{2})$ is simply $D\boldsymbol{\caF}$.
	
	We use the following quantitative version of  the implicit function theorem. Its proof is rather elementary and standard yet provided in Appendix \ref{sec:Imp} for completeness.
	\begin{lem}\label{lem:Imp}
		Let $F:\R^{n}\times \R^{m}\to\R^{m}$ be a $C^{1}$ function, let $\norm{\cdot}_{1}$ and $\norm{\cdot}_{2}$ be arbitrary norms on $\R^{n}$ and $\R^{m}$, respectively.
		Assume that $F(0,0)=0$ and for some $h_{x},h_{y}>0$ that
		\begin{equation}\label{eq:Imp_a1}
			C_{1}\deq 1\vee \norm{(D_{y}F(0,0))^{-1}}_{2,2}<\infty ,\qquad 
			\sup_{(x,y)\in V_{h_{x}}\times U_{h_{y}}}\norm{I-(D_{y}F(0,0))^{-1}D_{y}F(x,y)}_{2,2}\leq 1/2,
		\end{equation}
		where $V_{h}\deq \{x\in\R^{n}:\norm{x}_{1}\leq h\}$, $U_{h}\deq \{y\in\R^{m}:\norm{y}_{2}\leq h\}$, and $\norm{\cdot}_{2,2}$ is the operator norm on $\R^{m\times m}$ induced by the norm $\norm{\cdot}_{2}$ on $\R^{m}$. Set
		\begin{equation}\label{eq:Imp_a2}
			C_{2}\deq \sup_{(x,y)\in V_{h_{x}}\times U_{h_{y}}}\norm{D_{x}F(x,y)}_{2,1},
		\end{equation}
		where $\norm{\cdot}_{2,1}$ is the induced operator norm on $\caL((\R^{n},\norm{\cdot}_{1}),(\R^{m},\norm{\cdot}_{2}))$.  Then there is a $C^{1}$ map $g:V_{\wt{h}_{x}}\to U_{h_{y}}$, where 
		$\wt{h}_{x}\deq h_{x}\wedge (h_{y}/(2C_{1}C_{2}))\le h_{x}$, such that $y=g(x)$ is the unique solution of the implicit equation
		$F(x,y)=0$ in $U_{h_{y}}$ for each $x\in V_{\wt{h}_{x}}$. The derivative of $g$ satisfies
		\begin{equation}\label{eq:Imp_d}
			Dg(x)=(D_{y}F(x,y))^{-1}D_{x}F(x,y),\qquad \norm{Dg(x)}_{2,1}\leq 2C_{1}C_{2}, \qquad \forall x\in V_{\wt{h}_{x}}.
		\end{equation}
	\end{lem}

	We now check the assumptions of Lemma \ref{lem:Imp} for $\wt{F}$ around $(V_{1},V_{2},0,0)$. First of all, \eqref{eq:DF_Lip} and \eqref{eq:DtF} (together with $\norm{J_{m_{1},m_{2}}}_{\infty,2}\leq 10$) imply
	\begin{equation}
		\norm{D_{\bsw}\wt{F}(U_{1},U_{2},w_{1},w_{2})-D_{\bsw}\wt{F}(z_{1}\lone_{m_{1}},z_{2}\lone_{m_{2}},0,0)}\leq 1000(1/c+h)^{4}h
	\end{equation}
	whenever $\absv{w_{i}}\leq h$ and $\norm{U_{i}-z_{i}\lone_{m_{i}}}_{\infty}\leq h$. On the other hand we also have
	\begin{equation}
		D_{\bsw}\wt{F}(z_{1}\lone_{m_{1}},z_{2}\lone_{m_{2}},0,0)=DF_{\chi,p}(z_{1},z_{2}),
	\end{equation}
	where $F_{\chi,p}$
	 was defined in \eqref{eq:def_F} and $p\deq m_{1}/(m_{1}+m_{2})$. Therefore by Lemma \ref{lem:Jac}, we can choose small enough $h_{0}\equiv h_{0}(c)>0$ so that for $\norm{V_{i}-z_{i}\lone_{m_{i}}}_{\infty}\leq h_{0}$,
	\begin{equation}
		\norm{D_{\bsw}\wt{F}(V_{1},V_{2},0,0)^{-1}}\leq 2\norm{DF_{\chi,p}(z_{1},z_{2})^{-1}}\leq 2C_{1}(c).
	\end{equation}
	Likewise, reducing $h_{0}$ further if necessary, whenever $\norm{V_{i}-z_{i}\lone_{m_{i}}}_{\infty}\leq h_{0}$
	we have
	\begin{equation}
		\norm{I-D_{\bsw}\wt{F}(V_{1},V_{2},0,0)^{-1}D_{\bsw}\wt{F}(U_{1},U_{2},w_{1},w_{2})}\leq\frac{1}{2}
	\end{equation}
	in the domain 
	\begin{equation}
		\{(U_{1},U_{2},w_{1},w_{2}):\norm{U_{i}-V_{i}}_{\infty}\leq h_{0},\absv{w_{i}}\leq h_{0}\}.
	\end{equation}
	Also \eqref{eq:DF_Lip} together with $D_{\bsU}\wt{F}=D\boldsymbol{\caF}$ immediately proves \eqref{eq:Imp_a2} with $C_{2}\leq200c^{-4}$ in the same domain. Finally, we further reduce $h_{0}$ so that $h_{0}\leq c/2$.
	
	Now we are ready to apply Lemma \ref{lem:Imp}; taking $h_{x}=h_{y}\deq h_{0}/(4C_{1}(c)C_{2})$, 
	which depends only on $c$, whenever $\norm{V_{i}-z_{i}\lone_{m_{i}}}_{\infty}\leq h$, there is a $C^{1}$ map 
	\begin{equation}
		g\equiv (g_{1},g_{2}) :\{(U_{1},U_{2}):\norm{U_{i}-V_{i}}_{\infty}\leq h\}\to \{(w_{1},w_{2}):\absv{w_{i}}\leq h_{0}\}
	\end{equation}
	that satisfies $\wt{F}(U_{1},U_{2},g(U_{1},U_{2}))=\wt{F}(V_{1},V_{2},0,0)$. Now we take 
	\begin{equation}
		\wt{\caV}_{i}(t)\deq (1-t)V_{i}+tz_{i}\lone_{m_{i}}, 
	\end{equation}
	and finally define 
	\begin{equation}
		\caV_{i}(t)\deq \wt{\caV}_{i}(t)+g_{i}(\wt{\caV}_{1}(t),\wt{\caV}_{2}(t))\lone_{m_{i}},
	\end{equation}
	which is well-defined since $\norm{\wt{\caV}_{i}(t)-z_{i}\lone_{m_{i}}}_{\infty}\leq h$ due to the assumption that $\norm{V_{i}-z_{i}\lone_{m_{i}}}_{\infty}\leq h$.
	
	We next check the properties for $\caV_{i}(t)$ claimed in Lemma~\ref{lem:box}. Since $g(V_{1},V_{2})=(0,0)$ is the unique solution of $F(V_{1},V_{2},\cdot,\cdot)=0$ by Lemma \ref{lem:Imp}, we have $\caV_{i}(0)=V_{i}$.
	Additionally, we have 
	\begin{equation}
		\caV_{i}(1)=z_{i}\lone_{m_{i}}+g_{i}(z_{1}\lone_{m_{1}},z_{2}\lone_{m_{2}})\lone_{m_{i}},
	\end{equation}
	proving (0$'$) with the choice
	\begin{equation}
		(\wt{z}_{1},\wt{z}_{2})=(z_{1},z_{2})+g(z_{1}\lone_{m_{1}},z_{2}\lone_{m_{2}}).
	\end{equation} 
	The second property (ii) is an immediate consequence of the fact that $\absv{g_{i}}\leq h_{0}$. The third property (iii) follows from the estimate (uniform in $t$)
	\begin{equation}
		\norm{\caV_{i}'}_{\infty}
		\leq\norm{(V_{1}-z_{1}\lone_{m_{1}},V_{2}-z_{2}\lone_{m_{2}})}_{\infty}(1+\norm{Dg(\wt{\caV}_{1},\wt{\caV}_{2})}_{\infty,2})\leq h(1+4C_{1}(c)C_{2}),
	\end{equation}
	where the last inequality follows from \eqref{eq:Imp_d}. Finally, from the definitions of $\boldsymbol{\caF}$ and $g$, we find
	\begin{equation}
		\boldsymbol{\caF}(\caV_{1}(t),\caV_{2}(t))=\wt{F}(V_{1},V_{2},0,0)=\boldsymbol{\caF}(V_{1},V_{2})
	\end{equation}
	for all $t\in[0,1]$, which implies (iv) and (v$''$). This completes the proof of Lemma \ref{lem:box}.
\end{proof}
In practice, Lemma \ref{lem:box} is used to shrink two ``clusters" $V_1, V_2$ of eigenvalues of $B$ around $z_1, z_2$
 into two point masses at $\wt{z}_1, \wt{z}_2$. Notice that we require the total masses of the two clusters $m_{1},m_{2}$ to be comparable to each other. Lemma \ref{lem:box} requires the two clusters to have 
  diameter at most $h$. Therefore, 
   below we will consider a mesh on $\spec(B)$ with mesh size $h$, and apply Lemma \ref{lem:box} to $V_{1},V_{2}$ that contain the eigenvalues in a pair of boxes in the mesh.

In the next lemma we first show that our initial matrix $B$ has a positive fraction of eigenvalues away from the imaginary axis on both sides. Later on, it will be used to guarantee that the eigenvalues of $B$ around any given point on the left half-plane can be paired with another (multi-)set of eigenvalues on the right half-plane, with the latter strictly away from the imaginary axis. Then the centers of boxes of such a pairing shall serve as the input $(z_{1},z_{2})$ of Lemma \ref{lem:Jac}, since the assumption \eqref{eq:z1z2} requires only one of $z_{1},z_{2}$ to be far away from the imaginary axis.

\begin{lem}\label{lem:box_ref}
	Let $\frC>1$ be a constant and let $B\in\C^{N\times N}$ be a diagonal matrix with $\norm{B},\norm{B^{-1}}\leq\frC$, $\brkt{B^{2}B\adj}=0$, and $\brkt{B^{3}B\adj}\geq 0$. Then, there exists a constant $c\equiv c(\frC)>0$ such that
	\begin{equation}\label{eq:mass_HP}
		\absv{\{i:\re B_{ii}<-c\}}>cN,\qquad \absv{\{i:\re B_{ii}>c\}}>cN.
	\end{equation}
\end{lem}
\begin{proof}
	We write $\bbH_{c,\pm}\deq \{\pm\re z>c\}$ and $\bbH_{c}\deq \bbH_{c,+}\cup\bbH_{c,-}$ for the complement of the strip of width $2c$ around the imaginary axis. 
	First we prove that the condition  $\brkt{B^{3}B\adj}\geq 0$ implies that not all $B_{ii}$ 
	may be in this strip. Notice that for all $0<c\leq 1/(3\frC)$ it holds that $\absv{z}\geq \frC^{-1}$ and $z\notin \bbH_{c}$ implies 
	\begin{equation}
		\cos (2\arg z)\leq -\frac{1}{2},
	\end{equation}
	so that
	\begin{equation}
		\re\frac{1}{N}\sum_{B_{ii}\notin \bbH_{c}}B_{ii}^{3}\ol{B}_{ii}\adj=\frac{1}{N}\sum_{B_{ii}\notin \bbH_{c}}\absv{B_{ii}}^{4}\cos(2\arg B_{ii})\leq-\frac{\#(B_{ii}\notin\bbH_{c})}{N}\frac{\frC^{-4}}{2}<0.
	\end{equation}
	Therefore we get
	\begin{equation}
		0\leq \brkt{B^{3}B\adj}=\re\brkt{B^{3}B\adj}
		\leq -\frac{\#(B_{ii}\notin \bbH_{c})}{N}\frac{\frC^{-4}}{2}
		+\frac{\#(B_{ii}\in\bbH_{c})}{N}\frC^{4},
	\end{equation}
	which in turn implies for any $c\leq 1/(3\frC)$ that
	\begin{equation}\label{eq:mass_2HP}
		\#(B_{ii}\in\bbH_{c})\geq \frac{1}{1+2\frC^{8}}N.
	\end{equation}
	
	Now we make use of the condition $\brkt{B^{2}B\adj}=0$ to prove that the mass is not concentrated on either of the half-planes $\bbH_{c,\pm}$. By \eqref{eq:mass_2HP}, we may assume without loss of generality that 
	\begin{equation}
		\#(B_{ii}\in\bbH_{1/(3\frC),+})\geq \frac{1}{2(1+2\frC^{8})}N.
	\end{equation} 
	Thus it only remains to prove that $\#(B_{ii}\in\bbH_{c,-})>cN$ for some $c>0$. Notice that for all $c\leq 1/(3\frC)$
	\begin{equation}
		\re\frac{1}{N}\sum_{B_{ii}\in\bbH_{c,+}}B_{ii}^{2}\ol{B}_{ii}\geq \frac{\#(B_{ii}\in\bbH_{1/(3\frC),+})}{N}\frac{1}{\frC^{2}}\frac{1}{3\frC}\geq \frac{1}{20\frC^{11}},
	\end{equation} 
	where in the first inequality we used that the partial sum on the leftmost side increases as $c>0$ decreases. 
	Since $\brkt{B^{2}B\adj}=0$, we thus have that
	\begin{equation}\label{eq:mass_HP+}
		-\re\frac{1}{N}\sum_{B_{ii}\notin\bbH_{c,+}}B_{ii}^{2}\ol{B}_{ii}\geq \frac{1}{20\frC^{11}}
	\end{equation}
	for all $c\leq 1/(3\frC)$. Then we plug the trivial bounds
	\begin{equation}
		-\re\frac{1}{N}\sum_{B_{ii}\in\bbH_{c,-}}B_{ii}^{2}\ol{B}_{ii}
		\leq \frac{\#(B_{ii}\in\bbH_{c,-})}{N}\frC^{3},\qquad 
		\Absv{\re \frac{1}{N}\sum_{B_{ii}\notin\bbH_{c}}B_{ii}^{2}\ol{B}_{ii}}\leq \frC^{2}c
	\end{equation}
	into \eqref{eq:mass_HP+} so that
	\begin{equation}\label{eq:mass_1HP}
		\frac{\#(B_{ii}\in\bbH_{c,-})}{N}\frC^{3}\geq \frac{1}{20\frC^{11}}-\frC^{2}c.
	\end{equation}
	Choosing small enough $c\equiv c(\frC)$ in \eqref{eq:mass_1HP} proves $\#(B_{ii}\in\bbH_{c,-})> cN$ as desired. This concludes the proof of Lemma \ref{lem:box_ref}.
\end{proof}

When $B$ is a perturbation (in the operator norm) of the two-point-mass case, Lemmas \ref{lem:box} and \ref{lem:box_ref} are already sufficient to shrink the eigenvalues of $B$ into two point masses. Going beyond this perturbative regime, say if $\spec(B)$ has only a point mass on the left half-plane but a diffuse density on the right, one immediately observes that the point mass has to be split into many pieces in order to match with the geometric partition (that is, a mesh) of spectrum on the right. Likewise, in the general case when $\spec(B)$ is scattered on both half-planes, one has to partition the meshes on both sides further to make pairings to apply Lemma~\ref{lem:box}. The fact that we use only 
 part of the mass at a point (or in a mesh) in each pair does no affect our proof, since our key relation \eqref{eq:key} is a plain sum over $\spec(B)$ so that we only need to keep the key quantities constant for each pair and then sum them up for all pairs.

The last delicate problem in applying Lemma \ref{lem:box} is that the mesh size, $h$, depends on the masses $m_{1},m_{2}$ of the two pair via their ratio $m_{1}/m_{2}$. While this is a natural consequence of our proof (due to \eqref{eq:Jac_p}), it forces us to fix the ratio \emph{a priori}, in particular without knowing how the mesh looks like. 

In the last preparatory step for the proof of Proposition \ref{prop:path_finite}, we resolve the above two issues with an elementary result concerning set partitions in Lemma \ref{lem:set_part}, which we apply to a mesh on $\spec(B)$. Recall that a \emph{set partition} of a set $\caI$ is a collection\footnote{Here, $\caP(\caI)$ means the power set of $\caI$.} $\caS\subset\caP(\caI)$ of subsets of $\caI$ such that $\bigsqcup_{\caJ\in\caS}\caJ=\caI$, i.e. elements of $\caS$ are disjoint, and their union is $\caI$. For two partitions $\caS,\caS'$ of $\caI$, we say that $\caS'$ is a \emph{refinement} of $\caS$ if and only if any element $\caJ'\in\caS'$ is a subset of some $\caJ\in\caS$.
\begin{lem}\label{lem:set_part}
	Fix a constant $c\in(0,1]$, and let $N_{1},N_{2},m_{1},m_{2}\in\N$ satisfy 
	\begin{equation}\label{eq:int_part_assump}
		c\leq \frac{N_{1}}{N_{2}}\leq \frac{1}{c}, \qquad N_{1},N_{2}\geq \frac{8m_{1}m_{2}}{c}.
	\end{equation}
	If $\caI^{(1)},\caI^{(2)}$ are two sets with set partitions $\caS^{(1)}\subset \caP(\caI^{(1)}),\caS^{(2)}\subset \caP(\caI^{(2)})$ such that $\absv{\caI^{(j)}}=N_{j}$ and $\absv{\caS^{(j)}}=m_{j}$ for $j=1,2$, then there exist refinements $\wt{\caS}^{(j)}$ of $\caS^{(j)}$ with $\absv{\wt{\caS}^{(1)}}=\absv{\wt{\caS}^{(2)}}\leq m_{1}+m_{2}$ and a bijection $f:\wt{\caS}^{(1)}\to\wt{\caS}^{(2)}$ such that
	\begin{equation}\label{eq:set_part}
		\frac{c}{4}\leq \frac{\absv{f(\caJ)}}{\absv{\caJ}}\leq \frac{4}{c},\qquad \forall \caJ\in\wt{\caS}^{(1)}.
	\end{equation}
\end{lem}
\begin{rem}
	$N_{j}$'s must not be too small for Lemma \ref{lem:set_part}
	to be true, in particular when $c$ is small. The assumption \eqref{eq:int_part_assump} is designed to ensure this in a quantitative manner but not optimized.
\end{rem}
Before the detailed proof, we sketch the idea.
Our proof of Lemma \ref{lem:set_part} builds upon the model case $N_{1}=N_{2}$ that we discuss first. 
In this case, whenever we are given a bijection $\frg:\caI^{(1)}\to\caI^{(2)}$, the refinements
\begin{equation}
	\wt{\caS}^{(2)}:=\{\caJ^{(2)}\cap \frg(\caJ^{(1)}):\caJ^{(j)}\in\caS^{(j)}\},\qquad \wt{\caS}^{(1)}:=\{\frg^{-1}(\caJ):\caJ\in\wt{\caS}^{(2)}\}
\end{equation}
satisfy all the properties listed in Lemma \ref{lem:set_part} except $\absv{\wt{\caS}^{(1)}}\leq m_{1}+m_{2}$, where $f$ is simply defined by $f(\caJ)=\frg(\caJ)$. To control $\absv{\wt{\caS}^{(1)}}$, we take a special bijection $\frg=\iota_{2}^{-1}\circ\iota_{1}$ where the indexing bijection $\iota_{j}:\caI^{(j)}\to \bbrktt{N_{j}}$ is such that $\iota_{j}(\caJ)$ consists of consecutive integers (i.e. it has the form $\bbrktt{n_{1},n_{2}}$) for all $\caJ\in\caS^{(j)}$.
Note that the indexing bijection just realizes the partition of an abstract set $\caI^{(j)}$ of cardinality $N_j$ on the
set of integers $\bbrktt{N_{j}}$.

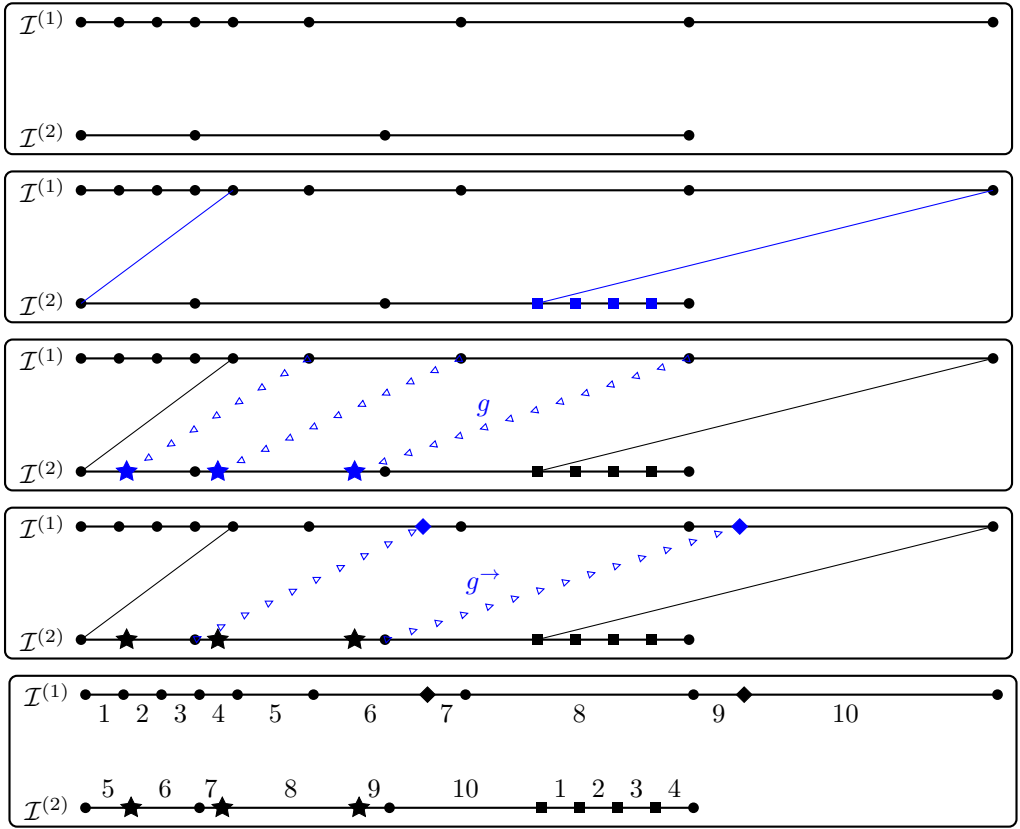
\begin{figure}
	\begin{tikzpicture}\color{black}
		\draw[black,thick, rounded corners=0.1cm](-1,-1.75) rectangle (12.25,0.25);
		\node at (-0.5,0) {$\caI^{(1)}$};
		\node at (-0.5,-1.5) {$\caI^{(2)}$}; 
		\draw[thick] (0, 0) -- (12, 0);
		\fill (0, 0) circle (2pt);    
		\fill (0.5, 0) circle (2pt);    
		\fill (1, 0) circle (2pt);    
		\fill (1.5, 0) circle (2pt);    
		\fill (2, 0) circle (2pt);    
		\fill (3, 0) circle (2pt);    
		\fill (5, 0) circle (2pt);    
		\fill (8, 0) circle (2pt);    
		\fill (12, 0) circle (2pt);    
		\node at (0.35, -.4) {${\caJ}_{1}^{(1)}$};
		\node at (1.1, -.4) {$\cdots$};
		\node at (1.85, -.4) {${\caJ}_{4}^{(1)}$};
		
		\draw[thick] (0, -1.5) -- (8, -1.5);
		\fill (0, -1.5) circle (2pt);   
		\fill (1.5, -1.5) circle (2pt); 
		\fill (4, -1.5) circle (2pt);   
		\fill (8, -1.5) circle (2pt);   
		
		\node at (6, -1.2) {$\caJ_{m_{2}}^{(2)}$};
	\end{tikzpicture}
	\medskip
	
	\begin{tikzpicture}{\color{black}
		\draw[black,thick, rounded corners=0.1cm](-1,-2.4) rectangle (12.25,0.25);
	\node at (-0.5,0) {$\caI^{(1)}$};
	\node at (-0.5,-1.5) {$\caI^{(2)}$}; 
		\draw[thick] (0, 0) -- (12, 0);
		\fill (0, 0) circle (2pt);    
		\fill (0.5, 0) circle (2pt);    
		\fill (1, 0) circle (2pt);    
		\fill (1.5, 0) circle (2pt);    
		\fill (2, 0) circle (2pt);    
		\fill (3, 0) circle (2pt);    
		\fill (5, 0) circle (2pt);    
		\fill (8, 0) circle (2pt);    
		\fill (12, 0) circle (2pt);    
		
		\draw[thick] (0, -1.5) -- (8, -1.5);
		\fill (0, -1.5) circle (2pt);   
		\fill (1.5, -1.5) circle (2pt); 
		\fill (4, -1.5) circle (2pt);   
		\fill (8, -1.5) circle (2pt);   
		\node at (5, -2) {$\wh{\caJ}_{m_{2}}^{(2)}$};
		
		\node at (6.45, -2) {$\wh{\caJ}_{m_{2}+1}^{(2)}$};
		\node at (7.2, -2) {$\cdots$};
		\node at (7.95, -2) {$\wh{\caJ}_{m_{2}+4}^{(2)}$};
		
		}
		{\color{blue}
		\fill[rectangle] (7.5, -1.5) ++(-2pt, -2pt) rectangle ++(4pt, 4pt);
		\fill[rectangle] (7, -1.5) ++(-2pt, -2pt) rectangle ++(4pt, 4pt);
		\fill[rectangle] (6.5, -1.5) ++(-2pt, -2pt) rectangle ++(4pt,4pt);
		\fill[rectangle] (6, -1.5) ++(-2pt, -2pt) rectangle ++(4pt,4pt);
		\draw (2,0) -- (0,-1.5);
		\draw (12,0) -- (6,-1.5);	
	}
		
\end{tikzpicture}
\medskip

	\begin{tikzpicture}{\color{black}
		\draw[black,thick, rounded corners=0.1cm](-1,-1.75) rectangle (12.25,0.25);
		\node at (-0.5,0) {$\caI^{(1)}$};
		\node at (-0.5,-1.5) {$\caI^{(2)}$}; 
		\draw[thick] (0, 0) -- (12, 0);
		\fill (0, 0) circle (2pt);    
		\fill (0.5, 0) circle (2pt);    
		\fill (1, 0) circle (2pt);    
		\fill (1.5, 0) circle (2pt);    
		\fill (2, 0) circle (2pt);    
		\fill (3, 0) circle (2pt);    
		\fill (5, 0) circle (2pt);    
		\fill (8, 0) circle (2pt);    
		\fill (12, 0) circle (2pt);    
		
		\draw[thick] (0, -1.5) -- (8, -1.5);
		\fill (0, -1.5) circle (2pt);   
		\fill (1.5, -1.5) circle (2pt); 
		\fill (4, -1.5) circle (2pt);   
		\fill (8, -1.5) circle (2pt);   
		\fill[rectangle] (7.5, -1.5) ++(-2pt, -2pt) rectangle ++(4pt, 4pt);
		\fill[rectangle] (7, -1.5) ++(-2pt, -2pt) rectangle ++(4pt, 4pt);
		\fill[rectangle] (6.5, -1.5) ++(-2pt, -2pt) rectangle ++(4pt,4pt);
		\fill[rectangle] (6, -1.5) ++(-2pt, -2pt) rectangle ++(4pt,4pt);
		\draw (2,0) -- (0,-1.5);
		\draw (12,0) -- (6,-1.5);	
	}
	{\color{blue}
		\node [star, star points=5, star point ratio=2, minimum size=8pt, inner sep=0pt, draw,fill] at (3/5,-1.5) {};
		\node [star, star points=5, star point ratio=2, minimum size=8pt, inner sep=0pt, draw,fill] at (9/5,-1.5) {};
		\node [star, star points=5, star point ratio=2, minimum size=8pt, inner sep=0pt, draw,fill] at (18/5,-1.5) {};
		\draw[decorate,decoration=triangles] (3,0) -- (3/5,-1.5);
		\draw[decorate,decoration=triangles] (5,0) -- (9/5,-1.5);
		\draw[decorate,decoration=triangles] (8,0) -- (18/5,-1.5);
		\node at (58/10-0.5,-0.65) {$g$};
	}
	\end{tikzpicture}
\medskip

\begin{tikzpicture}{\color{black}
	\draw[black,thick, rounded corners=0.1cm](-1,-1.75) rectangle (12.25,0.25);
	\node at (-0.5,0) {$\caI^{(1)}$};
	\node at (-0.5,-1.5) {$\caI^{(2)}$}; 
	\draw[thick] (0, 0) -- (12, 0);
	\fill (0, 0) circle (2pt);    
	\fill (0.5, 0) circle (2pt);    
	\fill (1, 0) circle (2pt);    
	\fill (1.5, 0) circle (2pt);    
	\fill (2, 0) circle (2pt);    
	\fill (3, 0) circle (2pt);    
	\fill (5, 0) circle (2pt);    
	\fill (8, 0) circle (2pt);    
	\fill (12, 0) circle (2pt);    
	
	\draw[thick] (0, -1.5) -- (8, -1.5);
	\fill (0, -1.5) circle (2pt);   
	\fill (1.5, -1.5) circle (2pt); 
	\fill (4, -1.5) circle (2pt);   
	\fill (8, -1.5) circle (2pt);   
	\fill[rectangle] (7.5, -1.5) ++(-2pt, -2pt) rectangle ++(4pt, 4pt);
	\fill[rectangle] (7, -1.5) ++(-2pt, -2pt) rectangle ++(4pt, 4pt);
	\fill[rectangle] (6.5, -1.5) ++(-2pt, -2pt) rectangle ++(4pt,4pt);
	\fill[rectangle] (6, -1.5) ++(-2pt, -2pt) rectangle ++(4pt,4pt);
	\node [star, star points=5, star point ratio=2, minimum size=8pt, inner sep=0pt, draw,fill] at (3/5,-1.5) {};
	\node [star, star points=5, star point ratio=2, minimum size=8pt, inner sep=0pt, draw,fill] at (9/5,-1.5) {};
	\node [star, star points=5, star point ratio=2, minimum size=8pt, inner sep=0pt, draw,fill] at (18/5,-1.5) {};
	\draw (2,0) -- (0,-1.5);
	\draw (12,0) -- (6,-1.5);
}
{\color{blue}
	\node[diamond, draw, inner sep=1.5pt, fill ] at (2+2.5,0) {};
	\node[diamond, draw, inner sep=1.5pt, fill ] at (2+4*5/3,0) {};
	\draw[decorate,decoration=triangles] (1.5,-1.5) -- (2+2.5,0);
	\draw[decorate,decoration=triangles] (4,-1.5) -- (2+4*5/3,0);		
	\node at (58/10-0.5,-0.75) {$g^{\rightarrow}$};
}
\end{tikzpicture}
\medskip

\begin{tikzpicture}\color{black}
	\draw[black,thick, rounded corners=0.1cm](-1,-1.75) rectangle (12.25,0.25);
	\node at (-0.5,0) {$\caI^{(1)}$};
	\node at (-0.5,-1.5) {$\caI^{(2)}$}; 
	\draw[thick] (0, 0) -- (12, 0);
	\fill (0, 0) circle (2pt);    
	\fill (0.5, 0) circle (2pt);    
	\fill (1, 0) circle (2pt);    
	\fill (1.5, 0) circle (2pt);    
	\fill (2, 0) circle (2pt);    
	\fill (3, 0) circle (2pt);    
	\node[diamond, draw, inner sep=1.5pt, fill ] at (2+2.5,0) {};
	\fill (5, 0) circle (2pt);    
	\node[diamond, draw, inner sep=1.5pt, fill ] at (2+4*5/3,0) {};
	\fill (8, 0) circle (2pt);    
	\fill (12, 0) circle (2pt);    
	
	\node at (0.25,-0.25) {1};
	\node at (0.75,-0.25) {2};
	\node at (1.25,-0.25) {3};
	\node at (1.75,-0.25) {4};
	\node at (2.5,-0.25) {5};
	\node at (3.75,-0.25) {6};
	\node at (4.75,-0.25) {7};
	\node at (6.5,-0.25) {8};
	\node at (50/6,-0.25) {9};
	\node at (10,-0.25) {10};
	
	\draw[thick] (0, -1.5) -- (8, -1.5);
	\fill (0, -1.5) circle (2pt);   
	\node [star, star points=5, star point ratio=2, minimum size=8pt, inner sep=0pt, draw,fill] at (3/5,-1.5) {};
	\fill (1.5, -1.5) circle (2pt); 
	\node [star, star points=5, star point ratio=2, minimum size=8pt, inner sep=0pt, draw,fill] at (9/5,-1.5) {};
	\node [star, star points=5, star point ratio=2, minimum size=8pt, inner sep=0pt, draw,fill] at (18/5,-1.5) {};
	\fill (4, -1.5) circle (2pt);   
	\fill[rectangle] (6, -1.5) ++(-2pt, -2pt) rectangle ++(4pt,4pt);
	\fill[rectangle] (6.5, -1.5) ++(-2pt, -2pt) rectangle ++(4pt,4pt);
	\fill[rectangle] (7, -1.5) ++(-2pt, -2pt) rectangle ++(4pt, 4pt);
	\fill[rectangle] (7.5, -1.5) ++(-2pt, -2pt) rectangle ++(4pt, 4pt);
	\fill (8, -1.5) circle (2pt);   
	\node at (3/10,-1.25) {5};
	\node at (21/20,-1.25) {6};
	\node at (33/20,-1.25) {7};
	\node at (27/10,-1.25) {8};
	\node at (38/10,-1.25) {9};
	\node at (5,-1.25) {10};
	\node at (6.25,-1.25) {1};
	\node at (6.75,-1.25) {2};
	\node at (7.25,-1.25) {3};
	\node at (7.75,-1.25) {4};
\end{tikzpicture}
\caption{Illustration of Lemma \ref{lem:set_part}\label{fig:part}: In the first graph, the two horizontal lines are realizations of $\caI^{(j)}$ on $\bbrktt{N_{j}}$, and the segments between marked points denote the partitions $\caS^{(j)}$. In the second graph, the blue lines mark the beginning and the end of almost-linear matching; the remaining four small segments in each line are exactly matched in size.
	The last graph depicts our construction of the refinements $\wt{\caS}^{(j)}$ and the bijection $f$ between them; see the proof of Lemma \ref{lem:set_part} for more details.}
\end{figure}

In the general case when $N_{1}>N_{2}$, we take $\frg=\iota_{2}^{-1}\circ g \circ\iota_{1}$ where $g:\bbrktt{N_{1}}\to\bbrktt{N_{2}}$ is  a `linear' map.
In general, strict linearity is not compatible with the discreteness of $\bbrktt{N_{1}}$ and $\bbrktt{N_{2}}$;
our map will be only approximately linear.
In practice one has to first deal with too small elements of $\caS^{(1)}$, for otherwise some image of the 'linear' map may be empty.
 Thus in the actual proof 
 $\wt{\caS}^{(1)}$ has to be chosen more carefully to make it a refinement of $\caS^{(1)}$. See Figure \ref{fig:part} for an illustration. 

\begin{proof}
	Without loss of generality assume $N_{2}\leq N_{1}\leq N_{2}/c$. We order the elements of the partitions
	increasingly in their size as
	\begin{equation}
		\caS^{(j)}=\{\caJ^{(j)}_{1},\ldots,\caJ^{(j)}_{m_{j}}\},\qquad  \absv{\caJ^{(j)}_{1}}\leq \absv{\caJ^{(j)}_{2}}\leq\cdots\leq \absv{\caJ^{(j)}_{m_{j}}},\qquad j=1,2.
	\end{equation}
	
	We first remove the ``small" elements of $\caS^{(1)}$; this step is shown in Figure \ref{fig:part} as the passage from the top graph to the second. Define $k\in\bbrktt{0,m_{1}}$ such that 
	\begin{equation}\label{eq:int_part_cutoff}
		\absv{\caJ^{(1)}_{k}}\leq\frac{4}{c}<\absv{\caJ^{(1)}_{k+1}}.
	\end{equation} 
	Then $k\leq m_{1}-1$ by \eqref{eq:int_part_assump} since
	\begin{equation}
		\absv{\caJ^{(1)}_{m_{1}}}\geq\frac{N_{1}}{m_{1}}>\frac{4}{c}.
	\end{equation}
	It also follows from \eqref{eq:int_part_assump} that the largest element of $\caS^{(2)}$ satisfies
	\begin{equation}
		\absv{\caJ^{(2)}_{m_{2}}}-\sum_{i\in\bbrktt{k}}\absv{\caJ^{(1)}_{i}}\geq \frac{N_{2}}{m_{2}}-\frac{4m_{1}}{c}\geq \frac{cN_{2}-4m_{1}m_{2}}{cm_{2}}>0,
	\end{equation}
	so that we may take a partition of $\caJ^{(2)}_{m_{2}}$ into $k+1$ elements as 
	\begin{equation}
		\caJ^{(2)}_{m_{2}}=\bigsqcup_{i=0}^{k}\wh{\caJ}^{(2)}_{m_{2}+i},\qquad \absv{\wh{\caJ}^{(2)}_{m_{2}+i}}=\absv{\caJ^{(1)}_{i}},\,\,\,i=1,\ldots,k,
	\end{equation}
where $\wh{\caJ}^{(2)}_{m_{2}}$ contains the rest.  Thus the largest element of $\caS^{(2)}$ takes
care of all ``small" elements of $\caS^{(1)}$ and now we consider the remaining elements of both partitions. 
These are
	\begin{equation}
		\wh{\caS}^{(1)}\deq\{\caJ^{(1)}_{k+1},\cdots,\caJ^{(1)}_{m_{1}}\},\qquad \wh{\caS}^{(2)}\deq\{\caJ^{(2)}_{1},\ldots,\caJ^{(2)}_{m_{2}-1},\wh{\caJ}^{(2)}_{m_{2}}\},
	\end{equation}
	on the sets $\wh{\caI}^{(j)}\deq\bigsqcup_{\caJ\in\wh{\caS}^{(j)}}\caJ$, 
	i.e. $\wh{\caS}^{(j)}$ is a set partition of $\wh{\caI}^{(j)}$. 
	
	Now it suffices to find refinements $\check{\caS}^{(j)}$ of $\wh{\caS}^{(j)}$ with 
	\begin{equation}
		\absv{\check{\caS}^{(1)}}=\absv{\check{\caS}^{(2)}}\leq m_{1}+m_{2}-k
	\end{equation}
	and a bijection $\check{f}:\check{\caS}^{(1)}\to\check{\caS}^{(2)}$ satisfying \eqref{eq:set_part}. Indeed, given $\check{\caS}^{(j)}$ and such an $\check{f}$, it is easy to see that
	\begin{equation}\begin{aligned}
		\wt{\caS}^{(1)}&\deq\{\caJ^{(1)}_{1},\ldots,\caJ^{(1)}_{k}\}\sqcup \check{\caS}^{(1)},\qquad
		\wt{\caS}^{(2)}\deq\{\wh{\caJ}^{(2)}_{1},\ldots,\wh{\caJ}^{(2)}_{k}\}\sqcup\check{\caS}^{(2)},\\
		f(\caJ)&\deq\begin{cases}
			\wh{\caJ}^{(2)}_{m_{2}+i} & \text{if }\caJ\in \wt{\caS}^{(1)}\setminus\check{\caS}^{(1)}, \,\,\,\text{where $i$ is such that }\caJ=\caJ^{(1)}_{i}, \\
			\check{f}(\caJ) & \text{if }\caJ\in\check{\caS}^{(1)}
			\end{cases}
	\end{aligned}\end{equation}
	satisfy all required properties in Lemma~\ref{lem:set_part}. 
	 We thus focus on constructing $\check{\caS}^{(j)}$ in the rest of the proof.
	
	We first check that the sizes $\wh{N}_{j}\deq\absv{\wh{\caI}^{(j)}}$ are still comparable. Note that
	\begin{equation}\label{eq:hatNcomp}
		1\geq \frac{\wh{N}_{2}}{\wh{N}_{1}}=\frac{N_{2}-\sum_{i\in\bbrktt{k}}\absv{\caJ_{i}^{(1)}}}{N_{1}-\sum_{i\in\bbrktt{k}}\absv{\caJ_{i}^{(1)}}}> \frac{c}{2},
	\end{equation}
	where the first inequality is a trivial consequence of $N_{1}\geq N_{2}$, and the second can be checked from \eqref{eq:int_part_assump}, \eqref{eq:int_part_cutoff}, and $k\leq m_{1}$ via
	\begin{equation}
		N_{2}-\frac{c}{2}N_{1}\geq \frac{N_{2}}{2}\geq \frac{4m_{1}m_{2}}{c}\geq \frac{4m_{1}}{c}\geq \sum_{i\in\bbrktt{k}}\absv{\caJ_{i}^{(1)}}> \left(1-\frac{c}{2}\right)\sum_{i\in\bbrktt{k}}\absv{\caJ_{i}^{(1)}}.
	\end{equation}
	
	Next, for $j=1,2$ we introduce  indexing bijections $\iota_{j }:\wh{\caI}^{(j)}\to\bbrktt{\wh{N}_{j}}$ so that each image $\iota_{j}(\caJ)$ of $\caJ\in\wh{\caS}^{(j)}$ is a consecutive interval. That is, for some points
	 $x_{1},\ldots, x_{m_{1}-k}\in\bbrktt{\wh{N}_{1}}$ and $y_{1},\ldots,y_{m_{2}}\in\bbrktt{\wh{N}_{2}}$ we have 
	\begin{equation}\begin{aligned}
			\iota_{1}(\caJ_{k+i}^{(1)})&=\rrbra x_{i-1},x_{i}\rrbra,\,\,\,i\in\bbrktt{m_{1}-k},	\\
		\iota_{2}(\caJ_{\ell}^{(2)})&=\rrbra y_{\ell-1},y_{\ell}\rrbra,\,\,\,\ell\in\bbrktt{m_{2}-1},\qquad 
		\iota_{2}(\wh{\caJ}_{m_{2}}^{(2)})=\rrbra y_{m_{2}-1},y_{m_{2}}\rrbra,
	\end{aligned}\end{equation}
	with the convention $x_{0}=0=y_{0}$. 	
	Note that $x_{m_{1}-k}=\wh{N}_{1}$ and $y_{m_{2}}=\wh{N}_{2}$ by definition, and also that
\begin{equation}\label{eq:n1_not_small}
	x_{\ell}-x_{\ell-1}=\absv{\caJ^{(1)}_{k+\ell}}>\frac{4}{c}>2\wh{N}_{1}/\wh{N}_{2},\qquad \forall\ell\in\bbrktt{m_{1}-k}
\end{equation}
by \eqref{eq:int_part_cutoff} and \eqref{eq:hatNcomp}.

	We consider the `linear' map\footnote{For $x\in\R$ we adopt the usual notations $\lceil x\rceil\deq\min\{n\in\Z:n\geq x\}$ and $\lfloor x\rfloor\deq\max\{n\in\Z:n\leq x\}$.} $g:\bbrktt{\wh{N}_{1}}\to\bbrktt{\wh{N}_{2}}$ given by $g(x)\deq \lceil{x\cdot\wh{N}_{2}/\wh{N}_{1}}\rceil$. We also need a specific `linear' stretching $g^{\rightarrow}:\bbrktt{\wh{N}_{2}}\to\bbrktt{\wh{N}_{1}}$ defined by
	\begin{equation}
		g^{\rightarrow}(y)\deq\begin{cases}
			x_{i}& \text{if }(y-1)\cdot \wh{N}_{1}/\wh{N}_{2} <x_{i} \leq y\cdot \wh{N}_{1}/\wh{N}_{2}, \\
			\lfloor y\cdot \wh{N}_{1}/\wh{N}_{2}\rfloor& \text{otherwise}.
		\end{cases}
	\end{equation}
	Then $g^{\rightarrow}$ is injective by \eqref{eq:n1_not_small} (recall also $\wh{N}_{1}\geq\wh{N}_{2}$),
	and $g(g^{\rightarrow}(y))=y$ for any $y\in\bbrktt{1,\wh{N}_{2}}$. Furthermore it also follows that $g^{\rightarrow}(g(x_{i}))=x_{i}$ for each $i$.

	Finally we define the desired refinements and show \eqref{eq:set_part}.
	The set 
	\begin{equation}
		\caY\deq \{g(x_{i}):1\leq i\leq m_{1}-k\}\cup \{y_{\ell}:1\leq \ell\leq m_{2}\}\subset \bbrktt{1,\wh{N}_{2}}
	\end{equation}
	consists of at most $m_{1}+m_{2}-k$ elements and contains $\wh{N}_{2}=y_{m_{2}}$. In Figure \ref{fig:part}, we constructed $\caY$ in the shorter horizontal line in the third graph, and then took its image under $g^{\rightarrow}$ in the fourth graph.
	Let $\check{y}_{1}<\cdots<\check{y}_{\check{m}}$ be the elements of $\caY$, from which we define the partitions $\check{\caS}^{(j)}$ as (with the convention $\check{y}_{0}=0=g^{\rightarrow}(\check{y}_{0})$)
	\begin{equation}\begin{aligned}
		\check{\caS}^{(2)}&\deq\{\iota_{2}^{-1}(\rrbra\check{y}_{i-1},\check{y}_i\rrbra):i=1,\ldots,\check{m}\},	\\
		\check{\caS}^{(1)}&\deq\{\iota_{1}^{-1}(\rrbra g^{\rightarrow}(\check{y}_{i-1}),g^{\rightarrow}(\check{y}_{i})\rrbra):i=1,\ldots,\check{m}\}. 
	\end{aligned}\end{equation}

	Since $\{y_{\ell}\}\subset\caY$, it follows that $\check{\caS}^{(2)}$ is a refinement of $\wh{\caS}^{(2)}$. Likewise $\check{\caS}^{(1)}$ refines $\wh{\caS}^{(1)}$ as $\{x_{i}\}\subset g^{\rightarrow}(\caY)$. The two partitions are of the same length $\check{m}$ since $g^{\rightarrow}$ is injective. Finally, \eqref{eq:set_part} follows from the definition of $g$ since 
	\begin{equation}
		1\leq g^{\rightarrow}(y+1)-g^{\rightarrow}(y)\leq 2\frac{\wh{N}_{1}}{\wh{N}_{2}}\leq \frac{4}{c},\qquad \forall y\in\caY.
	\end{equation}
	See the last graph of Figure \ref{fig:part} for an illustration of how we match the two final partitions.
	This concludes the proof of Lemma \ref{lem:set_part}. 
\end{proof}

\begin{proof}[Proof of Proposition \ref{prop:path_finite}]
	We first apply Lemma \ref{lem:box_ref} to find a $\frC$-dependent parameter $0<c_{0}<1/(2\frC)$ satisfying $\rho_{B}(\{\re z<-c_{0}\})>c_{0}$ and $\rho_{B}(\{\re z>c_{0}\})>c_{0}$. We partition the total index set $\bbrktt{N}$ into four 
	sets as $\bbrktt{N}=\caI_{\mathrm{i}}^{+}\sqcup\caI_{\mathrm{i}}^{-}\sqcup\caI_{\mathrm{o}}^{+}\sqcup\caI_{\mathrm{o}}^{-}$, where
	\begin{equation}\begin{aligned}
			\caI_{\mathrm{i}}^{+}\deq \{i:0<\re B_{ii}\leq c_{0}\}, \qquad
			\caI_{\mathrm{i}}^{-}\deq \{i:-c_{0}<\re B_{ii}\leq 0\},	\\
			\caI_{\mathrm{o}}^{+}\deq \{i:\re B_{ii}>c_{0}\}, \qquad
			\caI_{\mathrm{o}}^{-}\deq \{i:\re B_{ii}\leq -c_{0}\},
	\end{aligned}\end{equation}
	and define $N_{\mathrm{i}}^{\pm}\deq \absv{\caI_{\mathrm{i}}^{\pm}}$ and similarly $N_{\mathrm{o}}^{\pm}$. 
	Here the Roman indices $\mathrm{i}$ and $\mathrm{o}$ stand for ``inside'' and ``outside'', respectively, while the plus and minus signs indicate $B_{ii}$'s in the half plane right and left to the imaginary axis.
	Now we divide $\caI_{\mathrm{o}}^{\pm}$ further into
	\begin{equation}\label{eq:oipm}
		\caI_{\mathrm{o}}^{\pm}=\caI_{\mathrm{oo}}^{\pm}\sqcup\caI_{\mathrm{oi}}^{\pm},	\qquad
		N_{\mathrm{oi}}^{\pm}\deq \absv{\caI_{\mathrm{oi}}^{\pm}}=\left\lceil \frac{c_{0}}{2}N_{\mathrm{i}}^{\mp}\right\rceil,
	\end{equation}
	where we used $N_{\mathrm{o}}^{\pm}>Nc_{0}\geq N_{\mathrm{i}}^{\mp}c_{0}$ from Lemma \ref{lem:box_ref}. The specific choice of $\caI_{\mathrm{oi}}^{\pm}$ does not affect the proof, only their cardinality matters as 
	we will match the index set $\caI_{\mathrm{oi}}^{\pm}$ with $\caI_{\mathrm{i}}^{\mp}$, 
	and $\caI_{\mathrm{oo}}^{+}$ with $\caI_{\mathrm{oo}}^{-}$. 
	Notice that $N_{\mathrm{oi}}^{\pm}$ is comparable to $N_{\mathrm{i}}^{\mp}$ with the flipped sign.	Then by  \eqref{eq:oipm} one can see that there exist $0<c<c_{0}$ depending only on $c_{0}$ so that 
	\begin{equation}
		4c<\frac{N_{\mathrm{oi}}^{\pm}}{N_{\mathrm{i}}^{\mp}},	\frac{N_{\mathrm{oo}}^{\pm}}{N_{\mathrm{oo}}^{\mp}}<\frac{1}{4c}.
	\end{equation}
	We apply Lemma \ref{lem:box} to find the parameters $h>0$ and $\frC_{1}$, which still depend only on $\frC>0$ via $c$.
	
	Next, we consider the grid in $\C\cong\R^{2}$ of mesh size $h/2$, i.e. the collection of boxes 
	\begin{equation}
		Q_{k\ell}\deq (kh/2,(k+1)h/2]\times (\ell h/2,(\ell+1)h/2],\qquad k,\ell\in\Z.
	\end{equation}
	The grid $(Q_{k\ell})$ naturally defines a partition of any given subset $\caI$ of $\bbrktt{N}$ as 
	\begin{equation}
		\caI=\bigsqcup_{k,\ell\in\Z^{2}}\left(\caI\cap\{i:B_{ii}\in Q_{k\ell}\}\right),
	\end{equation}
	and the union contains at most $4\norm{B}^{2}/h^{2}\leq 4\frC^{2}/h^{2}$ non-empty sets. Now we apply Lemma \ref{lem:set_part} to the sets $\caI_{\mathrm{oi}}^{+}$ and $\caI_{\mathrm{i}}^{-}$ and their partitions
	\begin{equation}\begin{aligned}
		\caS_{\mathrm{oi}}^{+}\deq \{\caI_{\mathrm{oi}}^{+}\cap\{i:B_{ii}\in Q_{k\ell}\}:k,\ell\in\Z,\,\,\spec(B)\cap Q_{k\ell}\neq \emptyset\},\\
		\caS_{\mathrm{i}}^{-}\deq \{\caI_{\mathrm{i}}^{-}\cap\{i:B_{ii}\in Q_{k\ell}\}:k,\ell\in\Z,\,\,\,\spec(B)\cap Q_{k\ell}\neq \emptyset\},
	\end{aligned}\end{equation}
	to find refinements $\wt{\caS}_{\mathrm{oi}}^{+}$ and $\wt{\caS}_{\mathrm{i}}^{-}$ with a common size $\absv{\wt{\caS}_{\mathrm{oi}}^{+}}=m=\absv{\wt{\caS}_{\mathrm{i}}^{-}}$. Properly indexing $\wt{\caS}_{\mathrm{oi}}^{+}=\{\caJ_{\mathrm{oi}}^{+}(j):1\leq j\leq m\}$ and similarly for $\wt{\caS}_{\mathrm{i}}^{-}$, we eventually have the disjoint unions
	\begin{equation}
		\caI_{\mathrm{oi}}^{+}=\bigsqcup_{j=1}^{m}\caJ_{\mathrm{oi}}^{+}(j),\qquad \caI_{\mathrm{i}}^{-}=\bigsqcup_{j=1}^{m}\caJ_{\mathrm{i}}^{-}(j),
	\end{equation}
	where, due to
	Lemma \ref{lem:set_part}, we have $m\leq \absv{\caS_{\mathrm{oi}}^{+}}+\absv{\caS_{\mathrm{i}}^{-}}\leq 8\frC^{2}/h^{2}$ and 
	\begin{equation}\label{eq:int_part_appl}
		c\leq\frac{1}{4}\frac{N_{\mathrm{oi}}^{+}}{N_{\mathrm{i}}^{-}}
		\leq \frac{\absv{\caJ_{\mathrm{oi}}^{+}(j)}}{\absv{\caJ_{\mathrm{i}}^{-}(j)}}\leq 4\frac{N_{\mathrm{oi}}^{+}}{N_{\mathrm{i}}^{-}}\leq \frac{1}{c},\qquad j=1,\ldots,m.
	\end{equation}
	Each $\caJ_{\mathrm{oi}}^{+}(j)$ is contained in some box $Q_{k\ell}$ since $\wt{\caS}_{\mathrm{oi}}^{+}$ is a refinement of $\caS_{\mathrm{oi}}^{+}$, and the same applies to $\caJ_{\mathrm{i}}^{-}(j)$.
	
	We then apply Lemma \ref{lem:box} to the pair of vectors
	\begin{equation}
		(B_{ii})_{i\in\caJ_{\mathrm{oi}}^{+}(j)},\qquad (B_{ii})_{i\in\caJ_{\mathrm{i}}^{-}(j)}
	\end{equation}
	for each $j=1,\ldots,m$. To check the assumptions of Lemma \ref{lem:box}, notice that the lengths of the two vectors are comparable by \eqref{eq:int_part_appl}, and that, if $Q_{1}$ and $Q_{2}$ are the boxes containing $\caJ_{\mathrm{oi}}^{+}(j)$ and $\caJ_{\mathrm{i}}^{-}(j)$, then the centers $z_{1},z_{2}$ of $Q_{1},Q_{2}$ satisfy
	\begin{equation}
		\re z_{1}>c_{0}>c,\qquad \re z_{2}<0, \qquad \absv{z_{i}}\leq \frC+\frac{h}{2}<\frac{1}{c},\qquad \absv{z_{i}}>\frac{1}{\frC}-h>c.
	\end{equation}
	Consequently, we have found a vector-valued path $(\caB_{ii}(t))_{i\in\caI_{\mathrm{oi}}^{+}\cup\caI_{\mathrm{i}}^{-}}$ satisfying all properties (0$'$)--(v$''$) listed in Lemma~\ref{lem:box}, and in particular (0$'$) implies
	\begin{equation}
		\absv{\{\caB_{ii}(1):i\in\caI_{\mathrm{oi}}^{+}\cup\caI_{\mathrm{i}}^{-}\}}\leq 2m\leq \frac{16\frC^{2}}{h^{2}}.
	\end{equation}
	
	We repeat the procedure, i.e., applying Lemma \ref{lem:set_part}
	and then Lemma \ref{lem:box}, for the pairs of index sets $(\caI_{\mathrm{oi}}^{-},\caI_{\mathrm{i}}^{+})$ and $(\caI_{\mathrm{oo}}^{-},\caI_{\mathrm{oo}}^{+})$. Evolving the resulting three vectors  (from these three pairs of index sets) simultaneously, we now have found $\caB_{t}$ with all desired properties. For example, property (v$'$) of Proposition \ref{prop:path_finite} can be checked from \eqref{eq:box_cond_chi}, since
	\begin{equation}\begin{aligned}
			&N(\brkt{\caB_{t}^{3}\caB_{t}\adj}-\chi\brkt{\absv{\caB_{t}}^{4}})	\\
			=& \bigg[\sum_{j}\left(\sum_{i\in \caJ_{\mathrm{oi}}^{+}(j)}+\sum_{i\in\caJ_{\mathrm{i}}^{-}(j)}\right)
			+\sum_{j}\left(\sum_{i\in \caJ_{\mathrm{oi}}^{-}(j)}+\sum_{i\in\caJ_{\mathrm{i}}^{+}(j)}\right)	
			+\sum_{j}\left(\sum_{i\in \caJ_{\mathrm{oo}}^{-}(j)}+\sum_{i\in\caJ_{\mathrm{oo}}^{+}(j)}\right)
			\bigg]
			\\
			&\times \left((\caB_{t})_{ii}^{3}(\ol{\caB}_{t})_{ii}-\chi\absv{(\caB_{t})_{ii}}^{4}\right)	\\
			=&N(\brkt{\caB_{0}^{3}\caB_{0}\adj}-\chi\brkt{\absv{\caB_{0}}^{4}})=0.
	\end{aligned}\end{equation}
	We omit  the remaining details but only remark the constant $M$ in item (0) of Proposition \ref{prop:path_finite} can be chosen as
	\begin{equation}
		M=\frac{100\frC^{2}}{h^{2}},
	\end{equation}
	which depends only on $\frC$ since $h$ does. This completes the proof of Proposition \ref{prop:path_finite}.
\end{proof}

\subsection{Step 2: From finitely supported $A$ to $N$-independent $\rho_{A}$}\label{sec:path_fix}

As a consequence of Step 1, we have reduced a general given $A$ to a matrix whose spectrum has finite cardinality. 
Still the eigenvalues and their relative multiplicities of the resulting matrix are $N$-dependent, unlike in \cite[Eqs. (1.3)--(1.5)]{Liu-Zhang2023}.
In Step 2, we construct a further flow making the outcome of Step 1 $N$-independent. Once the $N$-dependent generalization of \cite{Liu-Zhang2023} indicated in Remark \ref{rem:Liu_gen} is established, Step 2 would not be needed.

\begin{prop}\label{prop:path_fix}
	For each $M\in\N$ and $\frC>0$ there exist small constant
	\footnote{The subscript ${\mathrm{TV}}$ indicates that $\delta_{\mathrm{TV}}$ controls the total variation of the difference of spectral measures $\rho_{B_{0}}-\rho_{B_{1}}$.} $\delta_{\mathrm{TV}}\equiv \delta_{\mathrm{TV}}(M,\frC)>0$ and large constant $\frC'\equiv\frC'(M,\frC)$ that satisfies the following: For $j=0,1$, let $\bsz_{j}=(z_{j,1},\ldots,z_{j,M})\in \C^{M}$ and let $\bsn=(n_{j,1},\ldots,n_{j,M})\in \bbrktt{1,N}^{M}$ satisfy $\sum_{i=1}^{M}n_{j,i}=N.$
	Suppose that 
	\begin{equation}\label{eq:close_zn}
		\norm{\bsz_{0}-\bsz_{1}}_{\infty}\leq \delta_{\mathrm{TV}},\qquad \norm{\bsn_{0}-\bsn_{1}}_{\infty}\leq \delta_{\mathrm{TV}}N,
	\end{equation}
	and for both $j=0,1$ 
	that matrices  $B_{j}\deq \bigoplus_{i=1}^{M}z_{j,i} I_{n_{j,i}}$ satisfy (a)--(c) of Proposition \ref{prop:path_finite} with the common constant $\frC$. Then, there exists a $C^{1}$ path $[0,1]\ni t\mapsto \caB_{t}\in\C^{N\times N}$ such that
	\begin{itemize}
		\item[(0$''$)] $\caB_{0}=B_{0}$, $\caB_{1}=B_{1}$.
		\item[(i)--(iv)] Items (i)--(iv) of Lemma \ref{lem:path_inv}  hold for all $t\in[0,1]$, with the constant $\frC'$ in place of $\frC$.
		\item[(v$'''$)] For all $t\in[0,1]$ we have
		\begin{equation}
			\Absv{\frac{\dd\chi(B_{t})}{\dd t}}\leq \frC'\absv{\chi(B_{1})-\chi(B_{0})}.
		\end{equation}
	\end{itemize}
\end{prop}
\begin{proof}
	Applying Lemma \ref{lem:box_ref} to $B_{0}$, we can find $i_{1},i_{2}\in\bbrktt{1,M}$ and a $\frC$-dependent constant $c>0$ such that
	\begin{equation}
		\re z_{0,i_{1}}<-c,\qquad \re z_{0,i_{2}}>c, \qquad n_{0,i_{1}}>cN,\qquad n_{0,i_{2}}>cN.
	\end{equation}
	Taking a common permutation of the entries of each vector $\bsz_{0},\bsz_{1},\bsn_{0},\bsn_{1}$, we may assume without loss of generality that $i_{1}=1$ and $i_{2}=2$. Take 
	\begin{equation}
		\wh{n}_{i}=n_{0,i}\wedge n_{1,i}, \qquad i=1,\ldots,M,
	\end{equation}
	and we consider the following decompositions of $B_{t}$ for $j=0,1$;
	\begin{equation}\begin{aligned}
			B_{j}&=B_{j,1}\oplus B_{j,2}\oplus B_{j,3},	\\
			B_{j,1}&\deq z_{j,1}I_{\wh{n}_{1}}\oplus z_{j,2}I_{\wh{n}_{2}},\\
			B_{j,2}&\deq \bigoplus_{i=3}^{M} z_{j,i}I_{\wh{n}_{i}},	\\
			B_{j,3}&\deq \bigoplus_{i=1}^{M}z_{j,i}I_{n_{j,i}-\wh{n}_{i}}.
	\end{aligned}\end{equation}
	Notice that when $\delta_{\mathrm{TV}}<c/2$ we have $\wh{n}_{1},\wh{n}_{2}>cN/2$ by \eqref{eq:close_zn}  so that
	\begin{equation}
		\frac{c}{2}<\frac{\wh{n}_{1}}{\wh{n}_{2}}<\frac{2}{c}.
	\end{equation}
	
	In what follows, we construct a path $\caB_{t,k}$, $t\in [0,1]$,  connecting $B_{0,k}$ to $B_{1,k}$ for each $k=1,2,3$, and then take their direct sum to get $\caB_{t}$. More precisely, we crudely construct $\caB_{t,2}$ and $\caB_{t,3}$, and carefully design $\caB_{t,1}$ in order to fine--tune the parameters $\brkt{B^{2}B\adj}$ and $\chi(B)$. The crude flows are defined just by linear interpolations as follows;
	\begin{equation}\begin{aligned}
			\caB_{t,2}\deq&\bigoplus_{i=3}^{M} ((1-t)z_{0,i}+tz_{1,i})I_{\wh{n}_{i}}, \\
			\caB_{t,3}\deq&((1-t)\absv{B_{0,3}}+t\absv{B_{1,3}})\exp\left[\ii( (1-t)\Theta_{0}+t\Theta_{1})\right],
	\end{aligned}\end{equation}
	where in the second line we used the polar decomposition $B_{j,3}=\absv{B_{j,3}}\exp (\ii \Theta_{j})$ for $j=0,1$, with the diagonal matrix $\Theta_{j}$ whose entries are in $[0,2\pi)$. Then, for $\delta_{\mathrm{TV}}<\frC/100$, we easily find that the sub-matrix $\caB_{t,2}\oplus \caB_{t,3}$ satisfies items (i)--(iii) of Lemma \ref{lem:path_inv} with a $\frC$-dependent constant $\frC'$. 
	
	In order to design $\caB_{t,1}$, we now consider the time dependent quantities
	\begin{equation}\begin{aligned}
			q_{1}(t)\deq &\sum_{k=2,3}\frac{1}{\wh{n}_{1}+\wh{n}_{2}}\Tr \caB_{t,k}^{2}\caB_{t,k}\adj,	\\
			q_{2}(t)\deq &\sum_{k=2,3}\frac{1}{\wh{n}_{1}+\wh{n}_{2}}\Tr\left[\caB_{t,k}^{3}\caB_{t,k}\adj -\chi_{t}\absv{\caB_{t,k}}^{4}\right],
	\end{aligned}\end{equation}
	where we abbreviated
	\begin{equation}
		\chi_{t}\deq (1-t)\chi(B_{0})+t\chi(B_{1}).
	\end{equation}
	We will apply the implicit function theorem, Lemma \ref{lem:Imp}, to the map $\wh{F}:[0,1]\times \C^{2}\to\C^{2}$ (as a function from $[0,1]\times \R^{4}$ to $\R^{4}$, hence taking $n=1$ and $m=4$ in Lemma \ref{lem:Imp}) defined by
	\begin{equation}
		\wh{F}(t,w_{1},w_{2})\deq F_{\chi_{t},p}(z_{0,1}+w_{1},z_{0,2}+w_{2})+(q_{1}(t),q_{2}(t)), \qquad p\deq\frac{\wh{n}_{1}}{\wh{n}_{1}+\wh{n}_{2}},
	\end{equation}
	where recall that $F_{\chi, p}$ was defined in Lemma~\ref{lem:Jac}. Since the domain is finite dimensional, we may choose the Euclidean norm for both variables. We can easily verify the first condition in \eqref{eq:Imp_a1} from Lemma~\ref{lem:Jac}, with the constant $C_{1}$ in \eqref{eq:Imp_a1} depending only on $\frC$. The remaining two conditions can be easily checked for small enough $h_{x},h_{y}$ depending only on $\frC$ and $M$. In fact, by choosing small enough $\delta_{\mathrm{TV}}$ in \eqref{eq:close_zn} (with the threshold in terms of $\frC$ and $M$), one may take\footnote{Here we use a minor variant of Lemma \ref{lem:Imp} for $n=1$, where the ball $V_{h}=[-h_{x},h_{x}]$ is replaced with $[0,h_{x}$].} $h_{x}=1$ and $h_{y}>0$ depending on $\frC,M$ but not on $h$. To see this, we note that 
	\begin{equation}
		\Absv{\frac{\dd \chi_{t}}{\dd t}}=\absv{\chi_{1}-\chi_{0}}\lesssim h, \qquad \sup_{t\in[0,1]}\Absv{\frac{\dd q_{1}(t)}{\dd t}}+\Absv{\frac{\dd q_{2}(t)}{\dd t}}\lesssim \delta_{\mathrm{TV}},
	\end{equation}
	with the implicit constants depending only on $\frC$ and $M$, so that uniformly over $t\in[0,1]$ and $\absv{\bsw}\leq 1$ 
	\begin{equation}\begin{aligned}
		&\norm{D_{t}\wh{F}(t,\bsw)}=O(\delta_{\mathrm{TV}}),	&
		&\norm{D_{\bsw}\wh{F}(t,\bsw)}=O(1),	&
		&\\
		&\norm{D_{t}^{2}\wh{F}(t,\bsw)}=O(\delta_{\mathrm{TV}}),	& &\norm{D_{t}D_{\bsw}\wh{F}(t,\bsw)}=O(\delta_{\mathrm{TV}}),	& &\norm{D_{\bsw}^{2}\wh{F}(t,\bsw)}=O(1).
	\end{aligned}\end{equation}
	On the other hand, it also follows that the constant in the last condition \eqref{eq:Imp_a2} of Lemma \ref{lem:Imp} is $C_{2}=O(\delta_{\mathrm{TV}})$, so that 
	\begin{equation}
		\frac{h_{y}}{2C_{1}C_{2}}\gtrsim \frac{1}{\delta_{\mathrm{TV}}},
	\end{equation}
	with the implicit constant depending only on $\frC$ and $M$. Therefore by choosing small enough $\delta_{\mathrm{TV}}$ we get $\wt{h}_{x}=h_{x}\wedge(h_{y}/(2C_{1}C_{2}))=h_{x}=1$, hence the implicit function in Lemma \ref{lem:Imp} is defined on all of $[0,1]$.
	
	So far we have proved that, by taking small enough $\delta_{\mathrm{TV}}$ in \eqref{eq:close_zn}, there exists a unique $C^{1}$ map
	\begin{equation}
		[0,1]\ni t\mapsto (w_{1}(t),w_{2}(t))\in D(0,h_{y})^{2},
	\end{equation}
	for some $h_{y}$ depending only on $\frC$ and $M$, such that 
	\begin{equation}\label{eq:Imp_N_indep}
		\wh{F}(t,w_{1}(t),w_{2}(t))=F_{\chi_{t},p}(z_{0,1}+w_{1},z_{0,2}+w_{2})+(q_{1}(t),q_{2}(t))=0,
	\end{equation}
	and $\absv{\dd w_{i}(t)/\dd t}=O(1)$. We  define 
	\begin{equation}
		\caB_{t}:=\caB_{t,1}\oplus \caB_{t,2}\oplus \caB_{t,3},\qquad \caB_{t,1}\deq (z_{0,1}+w_{1}(t))I_{\wh{n}_{1}}\oplus (z_{0,2}+w_{2}(t))I_{\wh{n}_{2}}.
	\end{equation}
	
	We finally check the required properties in Proposition \ref{prop:path_fix}. The first property (0$''$) is due to the uniqueness of the implicit solution $(w_{1},w_{2})$ of \eqref{eq:Imp_N_indep} and the fact that
	\begin{equation}
		\wh{F}(1,z_{1,1}-z_{0,1},z_{1,2}-z_{0,2})=0.
	\end{equation}
	Items (i)--(iii) of Lemma \ref{lem:path_inv} follow from the construction. The last two properties (iv) and (v$'''$) can be checked from the definition of $\wh{F}$. In particular, we have a stronger form of (v$'''$) as
	\begin{equation}
		\chi(B_{t})=\chi_{t}=(1-t)\chi(B_{0})+t\chi(B_{1}).
	\end{equation}
	This completes the proof of Proposition \ref{prop:path_fix}.
\end{proof}

\section{Proof of Theorem \ref{theo:mainthm}}\label{sec:pf_main}
Recall the assumptions of Theorem \ref{theo:mainthm}: The sequences of matrices $X\equiv X^{(N)}$ and $A\equiv A^{(N)}$ are given respectively by Definitions \ref{defn:X} and \ref{defn:A}, and $A$ further satisfies
\begin{equation}
	\absv{\alpha(A)-\alpha_{\infty}}\leq N^{-\epsilon},
\end{equation}
for an $N$-independent constant $\alpha_{\infty}>-1/3$.

By a usual compactness argument, it suffices to prove the following statement:
There is a subsequence $N_{n}$ of $\N$ such that 
\begin{equation}\label{eq:goal}
	\lim_{n\to\infty}\int_{\C^{k}}F(\bsz)(p_{k}^{(N_{n})}-p_{k,\alpha_\infty})\dd^{2k}\bsz=0,
\end{equation}
where $p_{k}^{(N)}$ is the $k$-point correlation function of $N^{1/4}\gamma(A)(A+X)$. 

In the first step, for each sufficiently large $N\geq N_{0}(\frC)$, we apply Proposition \ref{prop:path_finite} for $B_{0}=A^{-1}$ to construct a matrix-valued flow $\caA_{t}=\caB_{t}^{-1}$ 
such that $\caA_{1}$ can be written as
\begin{equation}
	\caA_{1}=\bigoplus_{i=1}^{M}z_{i}I_{n_{i}},
\end{equation}
where the number $M$ depends only on $\frC$. Then by Theorem \ref{theo:mainthmflow} we find that
\begin{equation}\label{eq:conv1}
	\int_{\C^{k}}F(\bsz)(p_{k,0}^{(N)}(\bsz)-p_{k,1}^{(N)}(\bsz))\dd^{2k}\bsz=O(N^{-c}),
\end{equation}
where $p_{k,t}^{(N)}$ stands for the $k$-point correlation function of $N^{1/4}\gamma(\caA_{t})(\caA_{t}+X)$.

Now we have found the two $N$-dependent sequence of vectors
\begin{equation}
	\bsz\equiv\bsz^{(N)}=(z_{1},\cdots,z_{M})\in \C^{M}, \qquad \bsc\equiv\bsc^{(N)}\deq (n_{1}/N,\cdots,n_{M}/N)\in (0,1)^{M},
\end{equation}
satisfying for some $\frC_{1}>0$ that
\begin{equation}\label{eq:norm_fin}
	\frC_{1}^{-1}\leq \absv{z_{i}}\leq \frC_{1},\qquad \sum_{i=1}^{M}c_{i}=1,\qquad \sum_{i=1}^{M}c_{i}\frac{1}{z_{i}^{2}\ol{z}_{i}}=0,
\end{equation}
as well as that
\begin{equation}\label{eq:alpha_fin}
	\sum_{i=1}^{M}c_{i}\left(\frac{1}{z_{i}^{3}\ol{z}_{i}}-\chi(A)\frac{1}{\absv{z_{i}}^{4}}\right)=0
\end{equation}
By the first condition in \eqref{eq:norm_fin}, we may take a subsequence $N_{n}$ of $\N$ so that $\bsc^{(N_{n})},\bsz^{(N_{n})}$ converges, i.e.
\begin{equation}\label{eq:cz_conv}
	\lim_{n\to\infty}\norm{\bsz^{(N_{n})}-\bsz^{(\infty)}}_{\infty}=\lim_{n\to\infty}\norm{\bsc^{(N_{n})}-\bsc^{(\infty)}}_{\infty}=0.
\end{equation}
Notice that the limits $\bsc^{(\infty)}$ and $\bsz^{(\infty)}$ automatically satisfy \eqref{eq:norm_fin}, and also that
\begin{equation}\label{eq:alpha_Nindep}
	\sum_{i=1}^{M}c^{(\infty)}_{i}\left(\frac{1}{(z^{(\infty)}_{i})^{3}\ol{z}^{(\infty)}_{i}}-\chi_{\infty}\frac{1}{\absv{z^{(\infty)}_{i}}^{4}}\right)=0,
\end{equation}
where we recall that
\begin{equation}
	\chi_{\infty}=\lim_{N\to\infty}\chi((A^{(N)})^{-1})=\lim_{n\to\infty}\chi((\caA^{(N_{n})})^{-1}). 
\end{equation}

Now we take sufficiently large $n\in\N$ so that 
\begin{equation}
	\Norm{\frac{1}{\bsz^{(N_{n})}}-\frac{1}{\bsz^{(\infty)}}}_{\infty}+\norm{\bsc^{(N_{n})}-\bsc^{(\infty)}}_{\infty}\leq \frac{\delta_{\mathrm{TV}}(M,\frC_{1})}{100},
\end{equation}
where $1/\bsz\deq(1/z_{i})_{1\leq i\leq M}$ and the threshold $h$ on the right-hand side is from Proposition \ref{prop:path_fix}, hence is determined by $\frC$ via $\frC_{1}$ in \eqref{eq:cz_conv}. Thus we may apply Proposition \ref{prop:path_fix} to get a flow $(\caA_{t}^{-1})_{t\in[1,2]}$ of $(N_{n}\times N_{n})$ diagonal matrices, starting from $\caA_{1}^{-1}$ and ending at $\caA_{2}^{-1}$ defined by
\begin{equation}\label{eq:A2}
	\caA_{2}=\bigoplus_{i=1}^{M}z_{i}I_{n^{(N_{n},\infty)}_{i}},
\end{equation}
	where $\bsn^{(N_{n},\infty)}\in\bbrktt{N_{n}}^{M}$ sums up to $N_{n}$ and
	\begin{equation}\label{eq:Nn}
		\max_{1\leq i\leq M}\absv{n^{(N_{n},\infty)}_{i}-c^{(\infty)}_{i}N_{n}}\leq M.
	\end{equation}
Moreover, from $\absv{\alpha(\caA_{2})-\alpha_{\infty}}=O(N^{-1})$ due to \eqref{eq:alpha_Nindep} and \eqref{eq:Nn}, we have
\begin{equation}
	\left|\frac{\dd\alpha(\caA_{t})}{\dd t}\right|\lesssim \absv{\alpha(\caA_{1})-\alpha(\caA_{2})}\leq N^{-\frc}.
\end{equation}
Then we apply Theorem \ref{theo:mainthmflow} once more
to find that
\begin{equation}\label{eq:conv2}
	\lim_{n\to\infty}\int_{\C^{k}}F(\bsz)(p_{k,1}^{(N_{n})}(\bsz)-p_{k,2}^{(N_{n})}(\bsz))\dd^{2k}\bsz=0,
\end{equation}
where now $p_{k,2}^{(N_{n})}$ is the $k$-point correlation function the ensemble
\begin{equation}
	N_{n}^{1/4}\gamma(\caA_{2})(\caA_{2}+X)\in\C^{N_{n}\times N_{n}}.
\end{equation}

In the last step, we aim at applying \cite[Theorem 1.2]{Liu-Zhang2023}. To this end we naturally extend $\bsn^{(N_{n},\infty)}$ for general $N$ beyond $N_{n}$ while keeping \eqref{eq:Nn} intact, defining a full ensemble $\caA_{2}+X$ of $(N\times N)$ matrices. Then, as in the beginning of the proof of Theorem \ref{theo:mainthmflow}, we use a two-moment-matching GFT to replace the IID matrix $X$ with the Ginibre ensemble. This incurs only an error of $O(N^{-c})$ for the eigenvalue correlation function $p_{k,2}^{(N)}$ of $\caA_{2}+X$. We now have the ensemble $\caA_{2}+X^{\mathrm{Gin}(\C)}$ with \eqref{eq:A2} and its rescaled $k$-point correlation function $p_{k,2}^{(N)}$, for which \cite[Theorem 1.2]{Liu-Zhang2023} directly applies and leads to
\begin{equation}\label{eq:conv3}
	\lim_{N\to\infty}\int_{\C^{k}}F(\bsz)(p_{k,2}^{(N)}(\bsz)-p_{k,\alpha_{\infty}}(\bsz))\dd^{2k}\bsz=0.
\end{equation}
Combining \eqref{eq:conv1}, \eqref{eq:conv2}, and \eqref{eq:conv3}, we conclude \eqref{eq:goal} as desired. This concludes the proof of Theorem~\ref{theo:mainthm}.
\appendix

\section{Proof of Lemma \ref{lem:alpha}}\label{sec:aux}
\begin{proof}[Proof of Lemma \ref{lem:alpha}]
	We start with the general identities 
	\begin{equation}\begin{aligned}\label{eq:Hess_comput}
			\caH_{11}=(\partial_{z}+\partial_{\ol{z}})^{2}\vert_{z=0}\brkt{\absv{A-z}^{-2}}=4\re \Brkt{\frac{1}{A^{3}A\adj}}+2\Brkt{\frac{1}{\absv{A^{2}}^{2}}},	\\
			\caH_{22}=-(\partial_{z}-\partial_{\ol{z}})^{2}\vert_{z=0}\brkt{\absv{A-z}^{-2}}=-4\re \Brkt{\frac{1}{A^{3}A\adj}}+2\Brkt{\frac{1}{\absv{A^{2}}^{2}}},	\\
			\caH_{12}=\ii (\partial_{z}^{2}-\partial_{\ol{z}}^{2})\vert_{z=0}\brkt{\absv{A-z}^{-2}}=-4\im \Brkt{\frac{1}{A^{3}A\adj}}
	\end{aligned}\end{equation}
	
	To prove of (i), we first show that $\alpha(A)\in(-1,1]$ for any $A$ having a criticality. The facts that $\lambda_{1}>0$ and $\alpha(A)\geq -1$ follow from
	\begin{equation}	\lambda_{1}+\lambda_{2}=\Delta_{z}\vert_{z=0}\brkt{\absv{A-z}^{-2}}=4\partial_{z}\partial_{\ol{z}}\vert_{z=0}\brkt{\absv{A-z}^{-2}}=\brkt{A^{-2}(A\adj)^{-2}}>0.
	\end{equation}
	Since $\lambda_{1}>0$, by definition we have $\alpha(A)\leq 1$.
		
	Next, we construct an $A$ for each value $\alpha(A)\in[-1/3,1]$. This case we consider
	\begin{equation}\label{eq:Ac}
		D\deq\diag (\pm1\pm\ii c)^{\oplus N/4},\qquad A_{c}\deq\brkt{\absv{D}^{-2}}^{1/2}D, \qquad c\in[0,1].
	\end{equation}
	Then we easily find that $A_{c}$ has a criticality at the origin. Also, since $\absv{D}^{2}=\brkt{\absv{D}^{2}}I$, it trivially follows that
	\begin{equation}
		\Brkt{\frac{1}{\absv{D^{2}}^{2}}}=\frac{1}{(1+c^{2})^{2}},\qquad \Brkt{\frac{1}{D^{3}D\adj}}=\brkt{\absv{D}^{2}}^{-1}\brkt{D^{-2}}=\frac{1}{1+c^{2}}\re\frac{1}{(1+\ii c)^{2}}=\frac{1-c^{2}}{(1+c^{2})^3}.
	\end{equation}
	We thus get (recall $c\leq 1$; for $c>1$ the last two quantities must be inverted.)
	\begin{equation}
		\alpha(A_{c})=\frac{-2(1-c^{2})+1+c^{2}}{2(1-c^{2})+1+c^{2}}=\frac{-1+3c^{2}}{3-c^{2}},
	\end{equation}
	which maps $c\in[0,1]$ to $[-1/3,1]$. Since we already covered the regime $\alpha(A)\in(-1,-1/3)$ in \eqref{eq:shapeparex1}--\eqref{eq:shapeparex}, this completes the proof of (i).
	
	We move on to prove (ii). The inclusion 
	\begin{equation}
		[-1/3,1]\subset\bigcup_{N\in\N,\frC>0}\{\alpha(A):A\in\mathrm{Crit}_{N,\frC},AA\adj=A\adj A\}
	\end{equation}
	is already proved, since the matrix $A_{c}$ defined in \eqref{eq:Ac} is normal. For the converse, we take a normal $A$, and assume without loss of generality that
	$\brkt{1/(A^{3}A\adj)}\geq0$ by rotating $A$ if necessary. Then $\caH$ becomes a diagonal matrix, for which we have
	\begin{equation}\begin{aligned}\label{eq:CS}
			\lambda_{1}+3\lambda_{2}=\caH_{11}+3\caH_{22}=8\Brkt{\frac{1}{A^{2}(A\adj)^{2}}}-8\Brkt{\frac{1}{A^{3}A\adj}}
			=8\Brkt{\frac{1}{A\adj}\left(\frac{1}{AA\adj}-\frac{1}{A^{2}}\right)\frac{1}{A}}
			\geq 0,
	\end{aligned}\end{equation}
	where the last inequality follows from Cauchy-Schwarz applied to the matrix inner product $[B_{1},B_{2}]\deq \Brkt{(A\adj)^{-1}B_{1}B_{2}\adj A^{-1}}$. This completes the proof of (ii). The last item (iii) also follows from \eqref{eq:CS}, since the inequality saturates if and only if $A$ satisfies $A=\e{\ii\varphi}A\adj$ for some $\varphi\in[0,2\pi)$. 
\end{proof}
\begin{rem}
	If $A$ is not normal, then the last inequality of \eqref{eq:CS} is false. While one can still apply Cauchy--Schwarz inequality to the inner product $[\cdot ,\cdot]$, what we obtain is
	\begin{equation}\begin{aligned}\label{eq:nonnormal_CS}
			\Brkt{\frac{1}{A\adj} \frac{1}{A^{2}}\frac{1}{A}}=\left[\frac{1}{A\adj},\frac{1}{A}\right]
			\leq \left[\frac{1}{A\adj},\frac{1}{A\adj}\right]^{1/2}\left[\frac{1}{A},\frac{1}{A}\right]^{1/2}
			=\Brkt{\frac{1}{A^{2}(A\adj)^{2}}}^{1/2}\Brkt{\frac{1}{\absv{A}^{4}}}^{1/2}.
	\end{aligned}\end{equation}
	Notice that $\brkt{1/\absv{A}^{4}}$ is typically bigger than $\brkt{1/(A^{2}(A\adj)^{2})}$ when $A$ is not normal.
\end{rem} 

\section{Hermitian $A$}\label{sec:Herm_A}
In this appendix, we construct the path when $B$ is Hermitian (i.e. real diagonal). This will correspond to the case when $A$ is Hermitian (up to rotation), i.e. $\alpha(A)=-1/3$.
\begin{prop}\label{prop:path_real} 
	Let $\frC>1$ be a constant and let $B\in\C^{N\times N}$ be a diagonal matrix with 
	\begin{itemize}
		\item[(a)] $\norm{B},\norm{B^{-1}}\leq \frC$,
		\item[(b)] $\brkt{B^{2}B\adj}=0$,
		\item[(c)] $B_{ii}$'s are real, or equivalently, $\chi(B)=\brkt{B^{3}B\adj}\brkt{\absv{B}^{4}}^{-1}=1$.
	\end{itemize}
	Then there exist constants $\frC_{1}>0$ depending only on $\frC$ and a piecewise $C^{1}$ map $[0,1]\ni t\mapsto\caB_{t}\in\C^{N\times N}$ satisfying 
	items (i)--(v$'$) of Proposition \ref{prop:path_finite} and, instead of (0), that $\caB_{0}=B$ and $\absv{\supp\rho_{\caB_{1}}}=2$.
\end{prop}
\begin{proof}
	Consider the two functions $f_{\pm}:[0,\infty)\to(0,\infty)$
	\begin{equation}
		f_{\pm}(s)\deq \frac{1}{N}\sum_{i:\pm B_{ii}>0}(\pm B_{ii}+s)^{3}.
	\end{equation}
	Then we have $f_{+}(0)=f_{-}(0)$, and both $f_{\pm}$ are smooth, strictly increasing, and 
	\begin{equation}
		\Absv{\frac{f_{\pm}(s)}{(1+s)^{3}}-\frac{\#(i:\pm B_{ii}>0)}{N}}\leq \frac{2 \frC^{3}}{s^{3}}.
	\end{equation}
	Thus, by Lemma \ref{lem:box_ref} we can find a constant $c>0$ such that
	\begin{equation}\label{eq:path_Herm}
		c(1+s)^{3}\leq f_{\pm}(s)\leq \frac{1}{c}(1+s)^{3},\qquad s\in[0,\infty).
	\end{equation}
	For simplicity we abbreviate $g\deq f_{-}^{-1}\circ f_{+}:[0,\infty)\to[0,\infty)$. Then $g(s)$ is almost linear for large $s$, i.e.
	\begin{equation}\label{eq:path_Herm_1}
		\Absv{\frac{g(s)}{s}-\left(\frac{\#(i:B_{ii}>0)}{\#(i:B_{ii}<0)}\right)^{1/3}}\leq C\frac{1}{(1+s)^{3}}
	\end{equation}
	for some constant $C$ depending only on $\frC$. Notice from Lemma \ref{lem:box_ref} that the second term on the left-hand side is bounded from above and below by constants.

	Define $s_{t}\deq t/(1-t)$ for $t\in[0,1)$ and
	\begin{equation}
		\caB_{t}\deq \frac{1}{1+s_{t}}\wt{\caB}_{s_{t}},
		\qquad
		(\wt{\caB}_{s})_{ii}=\begin{cases}
			B_{ii}+s,& B_{ii}>0,\\
			B_{ii}-f_{-}^{-1}(f_{+}(s)), & B_{ii}<0.
		\end{cases}
	\end{equation}
	The final point $\caB_{1}$ is simply defined by the limit as $t\to 1$, i.e.
	\begin{equation}
		(\caB_{1})_{ii}=\begin{cases}
			1, & B_{ii}>0, \\
			-(\#(B_{ii}>0)/\#(B_{ii}<0))^{1/3}, & B_{ii}<0.
		\end{cases}
	\end{equation}
	
	Now we check that $\caB_{t}$ satisfies the desired properties. Item (i) of Proposition \ref{prop:path_finite}, $\caB_{0}=B$, and $\absv{\supp\rho_{\caB_{1}}}=2$ can be checked immediately. The norm bound $\norm{\caB_{t}}\leq \frC_{1}$ in item (ii) is a consequence of
	\begin{equation}
		\norm{\wt{\caB}_{s}}\leq \norm{B}+s+f_{-}^{-1}(f_{+}(s))\leq C(1+s),
	\end{equation}
	where the last inequality is due to \eqref{eq:path_Herm}. The other bound $\norm{\caB_{t}^{-1}}\leq \frC_{1}$ can be proved similarly. Item (iv) follows from
	\begin{equation}
		\brkt{\caB_{t}^{2}\caB_{t}\adj}=s_{t}^{-3}(f_{+}(s_{t})-f_{-}(f_{-}^{-1}(f_{+}(s_{t}))))=0,
	\end{equation}
	and item (v$'$) holds since $\caB_{t}$ is Hermitian for all $t\in[0,1]$. Finally, for item (iii) on the regularity of $\caB_{t}$, firstly when $B_{ii}>0$ the flow $(\caB_{t})_{ii}=t+(1-t)B_{ii}$ is simply the convex combination so that we trivially have $\absv{\dd (\caB_{t})_{ii}/\dd t}\leq \frC$. Next for $B_{ii}<0$ we have
	\begin{equation}\label{eq:path_Herm_2}
		\Absv{\frac{\dd}{\dd t}(\caB_{t})_{ii}}=\Absv{\frac{\dd}{\dd t}\frac{B_{ii}+g(s)}{1+s_{t}}}\leq  \Absv{g'(s_{t})s_{t}-g(s_{t})}+\absv{g'(s_{t})-B_{ii}},
	\end{equation}
	where we abbreviated $g\deq (f_{-}^{-1}\circ f_{+})$. Notice that, by the same reasoning as in \eqref{eq:path_Herm_1},
	\begin{equation}
		g'(s)= \frac{\sum_{i:B_{ii}>0}(B_{ii}+s)^{2}}{\sum_{i:B_{ii}<0}(-B_{ii}+g(s))^{2}}=\frac{s^{2}}{g(s)^{2}}\frac{\#(i:B_{ii}>0)}{\#(i:B_{ii}<0)}(1+O((1+s)^{-2})),
	\end{equation}
	where the implicit constant depends only on $\frC$. Combining with \eqref{eq:path_Herm_1}, we find that 
	\begin{equation}
		\absv{g(s)-sg'(s)}+\absv{g'(s_{t})-B_{ii}}\lesssim 1,
	\end{equation}
	concluding the proof of Proposition \ref{prop:path_real}.
\end{proof}

\section{Proof of Lemma \ref{lem:Imp}}\label{sec:Imp}
\begin{proof}[Proof of Lemma \ref{lem:Imp}]
	We aim at applying a fixed point theorem to the $x$-dependent map $f_{x}:U_{h_{y}}\to \R^{m}$
	\begin{equation}
		f_{x}(y)\deq y-(D_{y}F(0,0))^{-1}F(x,y),\qquad x\in V_{\wt{h}_{x}},
	\end{equation}
	so that solving $F(x,y)=0$ amounts to finding a fixed point $y$ with $f_{x}(y)=y$. Note that the map $x\mapsto f_{x}(y)$ is Lipschitz in $V_{h_{x}}$ uniformly over $y\in U_{h_{y}}$ as 
	\begin{equation}\label{eq:Imp1}
		\norm{f_{x_{1}}(y)-f_{x_{2}}(y)}_{2}\leq \norm{D_{y}F(0,0)}^{-1}_{2,2} \sup_{(x,y)\in V_{h_{x}}\times U_{h_{y}}}\norm{D_{x}F(x,y)}_{2,1}\norm{x_{1}-x_{2}}_{1}\leq C_{1}C_{2}\norm{x_{1}-x_{2}}_{1}.
	\end{equation}
	Note also that $f_{x}$ is Lipschitz in $U_{h_{y}}$ with constant 1/2 uniformly over $x\in V_{h_{x}}$:
	\begin{equation}\label{eq:Imp2}
		\norm{f_{x}(y_{1})-f_{x}(y_{2})}_{2}\leq \sup_{(x,y)\in V_{h_{x}}\times U_{h_{y}}}\norm{I-(D_{y}F(0,0))^{-1} D_{y}F(x,y)}_{2,2}\norm{y_{1}-y_{2}}_{2}\leq \frac{1}{2}\norm{y_{1}-y_{2}}_{2}.
	\end{equation}
	Combining \eqref{eq:Imp1} and \eqref{eq:Imp2} immediately proves that $f_{x}$ is a self-map of $U_{h_{y}}$ for each $x\in V_{\wt{h}_{x}}$:
	\begin{equation}
		\norm{f_{x}(y)}_{2}\leq \norm{f_{x}(y)-f_{0}(y)}_{2}+\norm{f_{0}(y)-f_{0}(0)}_{2}\leq C_{1}C_{2}\norm{x}_{1}+\frac{\norm{y}_{2}}{2}\leq h_{y},\quad x\in V_{\wt{h}_{x}}, y\in U_{h_{y}}.
	\end{equation}
	
	Now we recursively define the sequence of functions
	\begin{equation}
		g_{0}(x)\deq f_{x}(0),\qquad g_{n+1}(x)\deq f_{x}(g_{n}(x)).
	\end{equation}
	One can easily find by induction that $g_{n}(x):V_{\wt{h}_{x}}\to \R^{m}$ is $C^{1}$ for each $n$. By \eqref{eq:Imp2}, we have
	\begin{equation}
		\norm{g_{n+1}(x)-g_{n}(x)}_{2}\leq\frac{1}{2}\norm{g_{n}(x)-g_{n-1}(x)}_{2}\leq\cdots\leq \frac{1}{2^{n}}\norm{f_{x}(0)}\leq \frac{C_{1}C_{2}\wt{h}_{x}}{2^{n}}\leq\frac{h_{y}}{2^{n+1}},
	\end{equation}
	so that $g_{n}(x)$ converges uniformly over $x\in V_{\wt{h}_{x}}$. Denoting the limit by $g$, we have $f_{x}(g(x))=g(x)$, i.e. $F(x,g(x))=0$. The Lipschitz estimate \eqref{eq:Imp2} also proves the uniqueness, for if $y\adj\in U_{h}$ were another solution we would have
	\begin{equation}
		\norm{g(x)-y\adj}_{2}=\norm{f_{x}(g(x))-f_{x}(y\adj)}_{2}\leq \frac{1}{2}\norm{g(x)-y\adj}_{2}.
	\end{equation}
	
	Finally, to show the regularity of $g$, we use \eqref{eq:Imp1} and \eqref{eq:Imp2} to get  
	\begin{equation}
		\norm{g(x_{1})-g(x_{2})}_{2}=\norm{f_{x_{1}}(g(x_{1}))-f_{x_{2}}(g(x_{2}))}\leq C_{1}C_{2}\norm{x_{1}-x_{2}}_{1}+\frac{1}{2}\norm{g(x_{1})-g(x_{2})}_{2},
	\end{equation}
	so that $g$ is Lipschitz with constant $2C_{1}C_{2}$. Given that $g$ is Lipschitz, one can easily derive \eqref{eq:Imp_d} by differentiating the implicit equation $F(x,g(x))=0$.
\end{proof}

\end{document}